\documentclass[3p, preprint]{elsarticle}
\usepackage[utf8]{inputenc}
\usepackage{amsmath}
\usepackage{graphicx}
\usepackage{subcaption}
\usepackage{mathSymbols}
\usepackage{geometry}
\usepackage{url}
\usepackage{graphicx}
\usepackage{color}
\usepackage[colorlinks=true, citecolor = blue]{hyperref}
\usepackage{array}
\usepackage{breqn}
\usepackage[defaultcolor=red, final]{changes}

\definechangesauthor[name=R2, color=red]{R2}
\definechangesauthor[name=R3, color=orange]{R3}
\definechangesauthor[name=martin, color=green]{martin}

\newcolumntype{x}[1]{>{\centering\arraybackslash\hspace{0pt}}p{#1}}

\newcommand{\scaleFig}{.45}
\newcommand{\visc}[1]{\text{m}^{#1}\text{s}^{-1}}

\title{Parallel-in-time integration of the shallow water equations on the rotating sphere using Parareal and MGRIT}

\author[1]{João Guilherme Caldas Steinstraesser\corref{cor1}}
\ead{joao.steinstraesser@usp.br}

\author[1]{Pedro da Silva Peixoto}
\ead{ppeixoto@usp.br}

\author[2,3,4]{Martin Schreiber}
\ead{martin.schreiber@univ-grenoble-alpes.fr}

\affiliation[1]{organization={Universidade de São Paulo},
addressline={Rua do Matão, 1010},
city={São Paulo},
postcode={05508-090},
country={Brazil}}

\affiliation[2]{organization={Université Grenoble Alpes},
addressline={621 avenue Centrale},
city={Saint-Martin-d'Hères},
postcode={38400},
country={France}}

\affiliation[3]{organization={Inria AIRSEA team
},
addressline={700 Avenue Centrale},
city={Grenoble},
postcode={38058},
country={France}}

\affiliation[4]{organization={Technical University of Munich},
addressline={Boltzmannstrasse 3},
city={Garching b. Muenchen},
postcode={85748},
country={Germany}}

\cortext[cor1]{Corresponding author}

\begin{document}

\begin{frontmatter}
    \begin{abstract}

    Despite the growing interest in parallel-in-time methods as an approach to accelerate numerical simulations in atmospheric modeling, improving their stability and convergence remains a substantial challenge for their application to operational models. In this work, we study the temporal parallelization of the shallow water equations on the rotating sphere combined with time-stepping schemes commonly used in atmospheric modeling due to their stability properties, namely an Eulerian implicit-explicit (IMEX) method and a semi-Lagrangian semi-implicit method (SL-SI-SETTLS). The main goal is to investigate the performance of parallel-in-time methods, namely Parareal and Multigrid Reduction in Time (MGRIT) when these well-established schemes are used on the coarse discretization levels and provide insights on how they can be improved for better performance. We begin by performing an analytical stability study of Parareal and MGRIT applied to a linearized ordinary differential equation depending on some temporal parallelization parameters, including the choice of a coarse scheme. Next, we perform numerical simulations of two standard tests in atmospheric modeling to evaluate the stability, convergence, and speedup provided by the parallel-in-time methods compared to a fine reference solution computed serially.  We also conduct a detailed investigation on the influence of artificial viscosity and hyperviscosity approaches, applied on the coarse discretization levels, on the performance of the temporal parallelization. Both the analytical stability study and the numerical simulations indicate a poorer stability behavior when SL-SI-SETTLS is used on the coarse levels, compared to the IMEX scheme. With the IMEX scheme, a better trade-off between convergence, stability, and speedup compared to serial simulations can be obtained under proper parameters and artificial viscosity choices, opening the perspective of the potential competitiveness for realistic models.
    
    % As a consequence, large second-order artificial viscosity coefficients are required when SL-SI-SETTLS is used, which, however, reduces the accuracy of the solution provided by the temporal parallelization; on the other hand, higher-order viscosity approaches can be applied when IMEX is used, leading to a better trade-off between convergence, stability and speedup compared to serial simulations.

    % the former case requires the use of large second-order artificial viscosity coefficients, which, however, compromises the accuracy of the solution, whereas higher-order viscosity approaches can be used in the latter case, leading to a better trade-off between convergence, stability and speedup compared to serial simulations.

\end{abstract}

\begin{keyword}
    Parallel-in-time methods \sep shallow water equations on the sphere \sep Parareal method \sep MGRIT method \sep atmospheric modeling
\end{keyword}
\end{frontmatter}

\section{Introduction}

\indent The numerical simulation of atmospheric circulation models, in the context of climate modeling and numerical weather prediction, is a challenge that motivates constant research efforts. Since the middle of the last century, several spatial and temporal discretization schemes have been proposed, seeking a balance between accuracy, numerical stability, and computational costs for simulations in large domains both in space and time \cite{williamson:2007}. Over the decades, advances in terms of computing technology, with increasing processing and memory resources available, as well as the advent of massively parallel high-performance computing (HPC) systems, opened new possibilities and motivated the development of new approaches in atmospheric modeling \cite{washington:2008}.

\indent Among the several temporal discretization approaches, parallel-in-time (PinT) methods have recently raised an increasing research interest. The term \quotes{PinT} refers to a large variety of numerical methods \cite{gander:2015}, of which the most popular, \eg Parareal \cite{lions_al:2001}, Multigrid Reduction in Time (MGRIT) \cite{friedhoff_al:2013} and Parallel Full Approximation Scheme in Space and Time (PFASST) \cite{emmett_minion:2012}, are iterative algorithms that, by using a fine, computationally expensive and one or more coarser, less expensive discretizations of the problem, allow to compute several time steps simultaneously, thus replacing the classical approach of serial time-stepping.

\indent The interest in PinT in various application domains has grown mainly in the past two decades, with these methods being seen as an alternative for overcoming the saturation of already well-established spatial parallelism approaches and taking more advantage of HPC systems. This increasing popularity is also explained by the non-intrusive character of several PinT methods, allowing the user to combine them with arbitrary temporal and spatial discretizations, and the availability of open-source libraries implementing them, \eg Xbraid (MGRIT) \cite{xbraid-package} and LibPFASST (PFASST) \cite{libpfasst:2018}, which makes the temporal parallelization of an operational code a relatively easy task. However, PinT is still making its first steps in atmospheric circulation models. Indeed, examples of efficient application of PinT methods, in the sense of effectively allowing to reduce the time-to-solution when compared to serial time-stepping, mainly include problems of parabolic and diffusive nature, \eg \cite{trindade_pereira:2004, geiser_guttel:2012, astorino_al:2012, gotschel:2018}. \replaced[id=R2]{Also, temporal parallelization suffers from stability and convergence issues when applied to simple hyperbolic problems, \eg the linear advection equation. This issue has been identified and studied by several works, \eg \cite{bal:2005, gander_vanderwalle:2007, dai_maday:2011,de_sterck_al:2019}, and has in some way discouraged further application and investigation of PinT methods to more complex advection-dominated problems, such as those arising in atmospheric modeling.}{On the other hand, temporal parallelization suffers from stability and convergence issues when applied to hyperbolic and advection-dominated problems, such as those arising in atmospheric modeling. Several works have identified and studied these issues,} Notably, it is known that the lack of convergence and stability arise on high wavenumbers of the solution due to the mismatch of phase representations on the discretization levels \cite{ruprecht:2018}, and several approaches trying to overcome it have been proposed, \eg \cite{ruprecht_krause:2012,chen_al:2014,haut_wingate:2014}, but most of them have limited application or introduce algorithmic complications that reduce the possible gains in terms of time-to-solution \added[id=R2]{and diminish the non-intrusive character of the PinT methods, making them less attractive to complex applications such as atmospheric circulation problems}.

\indent \replaced[id=R2]{However, even if most studies of PinT methods for hyperbolic problems still focus on relatively simple one-dimensional models in order to develop a better understanding,}{Therefore, most studies of PinT methods for hyperbolic problems still focus on relatively simple one-dimensional models to develop a better understanding; however,} some recent works already propose studies towards the temporal parallelization of more complete atmospheric models. These works focus on the shallow water equations (SWE) on the rotating sphere, a two-dimensional model commonly used as a starting point in atmospheric modeling since it contains most of the numerical and implementation challenges related to the horizontal discretization found in more complete, tridimensional models \cite{williamson_al:1992}. \cite{abel_al:2020} studies an asymptotic MGRIT method based on the asymptotic Parareal proposed and applied to the one-dimensional rotating SWE by \cite{haut_wingate:2014}; \cite{hamon_schreiber:2020} implemented the PFASST method combined with a spherical harmonics discretization of the rotating SWE; and \cite{schreiber_al:2019} developed a parallel-in-time method based on a temporal discretization using rational approximation for exponential integrators.
% \added[id=R2]{It should also be noticed that research efforts have been made to apply PinT methods for solving more complex and operational models, with some results pointing out good results compared to serial simulations and others indicating that the considered PinT methods present slow convergence and are not able to provide speedup. For instance, the Parareal algorithm or proposed modifications to it have been applied, in the context of ocean circulation models, to the Princeton Ocean Model \cite{liu_hu:2008} and the Finite-volumE Sea ice-Ocean circulation Model (FESOM2) \cite{philippi_slawig:2022, philippi_slawig:2023} and in the context of magnetic plasma simulation by \cite{berry_al:2012, samaddar_al:2019}.}
% \added[id=martin]{martin@Joao: Please don't just write that there are speedups/downs and we list some papers, but cite these papers by also explaining their results. E.g., this method claimed speedups for special blabla. This one has been applied in a simplified context, but still not showing real speedups. Also, it's very important to differentiate between real speedups (wallclocktime) and pseudo speedups (mathematicians abusing speedup meaning)}
\added[id=R2]{It should also be noted that research efforts have been made to apply PinT methods to solve more complex and operational models. For instance, in the context of ocean circulation and sea-ice model, Parareal has been applied to the Finite-volumE Sea ice-Ocean circulation Model (FESOM2) by \cite{philippi_slawig:2022}, presenting slow convergence and stability issues that prevented theoretical speedups, with better convergence results, in the same context, being obtained by \cite{philippi_slawig:2023} with a variant of Parareal using spatial coarsening, but with still present instabilities in long-term simulations and lack of real wall-clock speedup due to parallel overheads. In the context of magnetically confined plasma simulation in tokamaks, in the Integrated Plasma Simulator (IPS) \cite{berry_al:2012} and the nuclear fusion research project ITER \cite{samaddar_al:2019}, an event-based modification of Parareal allows for well-measured speedup, parallel efficiency, and parallel scaling results.}

\indent In this work, we also study the application of PinT methods to the SWE on the rotating sphere. Here, we focus on two- and multilevel temporal parallelization using Parareal and MGRIT. We combine them with popular temporal and spatial discretization schemes for atmospheric circulation, used operationally and/or for research purposes in weather and climate forecast models, to evaluate if well-established numerical methods in the atmospheric modeling community would suit temporal parallelization. Namely, we consider a spectral discretization in space using spherical harmonics; in time, we consider two schemes, an Eulerian implicit-explicit (IMEX) one and the semi-Lagrangian semi-implicit SL-SI-SETTLS proposed by \cite{hortal:2002}. We highlight that there is a recent interest in semi-Lagrangian methods in the PinT framework, \eg with studies on the Burgers \cite{schmitt_al:2018} and advection \cite{sterck_al:2021, sterck_al:2022} equations indicating stability and convergence improvements of Parareal and/or MGRIT by using semi-Lagrangian coarse discretizations; thus, we seek to study if this behavior is also verified in the context of a more complex problem.

\indent This investigation is conducted following two approaches. First, we perform an analytical stability study of Parareal and MGRIT applied to a linearized ODE as a function of some parameters for the temporal parallelization, including the choice of a coarse time-stepping scheme (IMEX or SL-SI-SETTLS). This analysis is based on the work developed in the Parareal framework by \cite{staff_ronquist:2005}, and we extend it to the two-level MGRIT with arbitrary relaxation. Second, we perform numerical simulations to evaluate the PinT methods in terms of stability, convergence to a reference solution, and computational time compared to the reference simulation. We consider two test cases commonly used for studying the numerical simulation of the SWE on the rotating sphere. Moreover, we investigate the use of artificial viscosity and hyperviscosity on the coarse discretization levels and how they influence the stability and convergence of the temporal parallelization.

\indent The proposed study \replaced[id=R3]{is conducted having in sight possible practical applications in atmospheric modeling.}{is of great interest in terms of practical application in atmospheric modeling.} Indeed, the high complexity of operational models makes the non-intrusive character of PinT methods such as Parareal and MGRIT a very attractive feature. By combining these methods with time-stepping schemes, which are effectively used in these models, \added[id=R3]{and by considering discretization and parametric choices which are coherent with practical applications,} we may provide indications of the feasibility of their temporal parallelization and how the PinT performance can be improved by properly parametrizing it. 

\indent This paper is organized as follows: in Section \ref{sec:swe}, we present the shallow water equations on the rotating sphere and their discretization using IMEX and SL-SI-SETTLS in time and spherical harmonics in space; in Section \ref{sec:PinT} we describe the Parareal and MGRIT methods; the analytical stability study is developed in Section \ref{subsec:stability}; the numerical simulations for evaluating the performance of the methods and the influence of the artificial viscosity and hyperviscosity are presented in Section \ref{sec:numerical_tests}; finally, conclusions are presented in Section \ref{sec:conclusion}.

\section{The shallow water equations on the rotating sphere}
\label{sec:swe}

\subsection{Governing equations}
\label{subsec:governing_equations}

\indent The SWE on the rotating sphere read

\begin{equation}
    \label{eq:swe_sphere}
    \begin{aligned}
    \pdertt{\vecU} &= \vecL_{\vecG}(\vecU) + \vecL_{\vecC}(\vecU) + \vecN_{\vecA}(\vecU) + \vecN_{\vecR}(\vecU)\\
    & = \vecL \vecU + \vecN(\vecU)
    \end{aligned}
\end{equation}

\noindent where $\vecU = \vecthreeT{\Phi}{\xi}{\delta}$, $\Phi = \ol{\Phi} + \Phi' = gh$ is the geopotential field (with $g$ the gravitational acceleration, $h$ the fluid depth, and $\ol{\Phi}$ and $\Phi'$ the mean geopotential and the geopotential perturbation, respectively), $\xi := \vecz \cdot (\nabla \times \vecV)$ is the vorticity ($\vecz$ is the unit vector in the vertical direction), $\delta := \nabla \cdot \vecV$ is the divergence, $\vecV := \vectwoT{u}{v}$ is the horizontal velocity, $\vecL := \vecL_{\vecG} + \vecL_{\vecC}$, $\vecN := \vecN_{\vecA} + \vecN_{\vecR}$ and $\vecL_{\vecG}$, $\vecL_{\vecC}$, $\vecN_{\vecA}$ and $\vecN_{\vecR}$ are respectively the linear gravity, linear Coriolis, nonlinear advection and nonlinear rest terms, given respectively by

\begin{equation*}
    \begin{gathered}
        \vecL_{\vecG}(\vecU) = \vecthree{-\ol{\Phi}\delta}{0}{-\nabla^2\Phi},\qquad
        \vecL_{\vecC}(\vecU) = \vecthree{0}{-\nabla \cdot(f\vecV)}{\vecz \cdot \nabla \times (f\vecV)}\\
        \vecN_{\vecA}(\vecU) = \vecthree{-\vecV \cdot \nabla \Phi}{-\nabla \cdot (\xi \vecV)}{-\nabla^2 \left( \frac{\vecV \cdot \vecV}{2} \right) + \vecz \cdot \nabla \times (\xi \vecV)},\qquad
        \vecN_{\vecR}(\vecU) = \vecthree{-\Phi'\delta}{0}{0}
    \end{gathered}
\end{equation*}

\indent Moreover, in the numerical simulations performed in this work, we consider an artificial (hyper)viscosity approach \cite{lauritzen_al:2011}. We then include the linear term

\begin{equation*}
\vecL_{\nu}(\vecU) = (-1)^{\frac{q}{2} + 1}\nu 
 \vecthree{\nabla^q \Phi'}{ \nabla^q \xi}{\nabla^q \delta }\end{equation*}

\noindent where $\nu \geq 0$ is the viscosity coefficient and $q \geq 2$ is the viscosity order, with $q$ even.

\subsection{Spatial discretization}

\indent In this work, the SWE equations on the rotating sphere are discretized in space using a spectral approach based on spherical harmonics, which is briefly presented below. We refer the reader to \cite{durran:2010} for details. We remark that spherical harmonics discretization is used in important operational applications in atmospheric modeling, \eg by the Integrated Forecasting System (IFS) at the European Centre for Medium-Range Weather Forecast (ECMWF) \cite{ECMWF:2003}, the Global Spectral Model (GSM) at the National Centers for Environmental Prediction (NCEP) of the U.S. National Oceanic and Atmospheric Administration (NOAA) \cite{NOAA_NCEP:2016} and the Global Spectral Model (GSM) of the Japan Meteorological Agency (JMA) \cite{kanamitsu_al:1983}.

\indent A given time-dependent smooth field $\psi(\lambda, \mu, t)$ defined on the sphere, where $\mu := \sin(\theta)$, $\lambda$ is the longitude, $\theta$ is the latitude, and $t$ is the time, can be written as a spherical harmonics expansion

\begin{equation}
    \label{eq:spherical_harmonics_expansion}
    \psi(\lambda,\mu, t) = \sum_{m = -\infty}^{\infty} \sum_{n = |m|}^{\infty} \psi_{m,n}(t)Y_{m,n}(\lambda,\mu)
\end{equation}

\noindent where $Y_{m,n}$ is the spherical harmonic function of zonal and total wavenumbers $m$ and $n$, respectively, and $\psi_{m,n}$ is the respective expansion coefficient. The spherical harmonics are defined as the product of associated Legendre functions and Fourier modes, so the direct and inverse spherical transforms can be performed via Fast Fourier and Fast Legendre transforms, respectively, in the zonal and meridional directions. In the implementation considered here, we consider a triangular truncation for the expansion \eqref{eq:spherical_harmonics_expansion}:

\begin{equation*}
    \psi(\lambda,\mu, t) = \sum_{m = -M}^{M} \sum_{n = |m|}^{M} \psi_{m,n}(t)Y_{m,n}(\lambda,\mu)
\end{equation*}

\indent Two main reasons motivate using spectral discretization for the SWE. First, spherical harmonics are eigenfunctions of the spherical Laplacian:

\begin{equation*}
    \nabla^2 Y_{m,n} = \frac{-n(n+1)}{a^2}Y_{m,n}
\end{equation*}

\noindent where $a$ is the sphere's radius. This property is interesting since the spherical Laplacian arises on the temporal discretization of the SWE. We remark that a pseudospectral approach is adopted, \replaced[id=R3]{relying on collocated (A-grid) gridpoint variables, with the nonlinear terms as well as the Coriolis term being computed in the physical space }{with the nonlinear terms being computed in the physical space }(the discretization size of the physical grid being determined by the well-known anti-aliasing \quotes{3/2-rule} \cite{durran:2010}). \replaced[id=R3]{Second, the spherical harmonics expansions and their properties related to the spherical Laplacian allow to easily formulate semi-implicit methods, thus avoiding the so-called \quotes{pole problem} found in physical grid discretizations used \eg in explicit finite difference methods: when uniform latitude-longitude meshes are used, cells near the poles have small longitudinal length, requiring the use of too small time step size to fulfill CFL stability constraints \cite{hack_ruediger:1992}.}{Second, the spherical harmonics expansion avoids the so-called \quotes{pole problem} found in physical grid discretizations used \eg in finite difference methods: when uniform latitude-longitude meshes are used, cells near the poles have small longitudinal length, requiring the use of too small time step size to fulfill CFL stability constraints in explicit time-stepping schemes.}

\subsection{Temporal discretization}
\label{sec:temporal_discretization}

\indent In this work, we consider two time-stepping schemes for discretizing \eqref{eq:swe_sphere}, namely an Eulerian, Strang-splitting implicit-explicit (IMEX) scheme and the semi-Lagrangian semi-implicit SL-SI-SETTLS proposed by \cite{hortal:2002}, which are briefly presented below. In both schemes, using a backward Euler method, the (hyper)viscosity term $\vecL_{\nu}$ is solved at the end of each time step.

\subsubsection{IMEX}
\label{sec:IMEX}

\indent In the implicit-explicit scheme considered here, the stiff, linear terms of the governing equations are treated implicitly, and the nonlinear ones are treated explicitly, allowing to overcome stability constraints imposed by the former, which makes it a popular time-stepping approach in fluid dynamics problems \cite{carpenter_al:2005}. \replaced[id=R3]{Although not as popular in operational numerical weather prediction models as semi-Lagrangian (SL) methods, IMEX schemes are also a relevant class of methods in atmospheric modeling, being used in research models such as the Model for Prediction Across Scales (MPAS) at the U.S.\,National Center for Atmospheric Research (NCAR) \cite{skamarock_al:2012} and the Nonhydrostatic Unified Model of
the Atmosphere (NUMA) at the Naval Research Laboratory (NPS/NRL) \cite{kelly_giraldo:2012}, including implicit-explicit approaches for the horizontal terms of the governing equations, as considered here, and also explicit and implicit treatments respectively for the horizontal and vertical terms \cite{mengaldo_al:2018}.}{Although not as popular in operational numerical weather prediction models as semi-Lagrangian (SL) methods, IMEX schemes are based on fixed spatial grids, being conceptually simpler and having better conservation properties. They are used in research models \eg at the U.S.\,National Center for Atmospheric Research (NCAR) and the Naval Research Laboratory (NPS/NRL), including implicit-explicit approaches for the horizontal terms of the governing equations, as considered here, and also explicit and implicit treatments respectively for the horizontal and vertical terms \cite{mengaldo_al:2018}.}

\indent We consider a second-order Strang-splitting IMEX scheme, with a half timestep of the implicit solver, followed by a full explicit timestep and a second implicit half timestep:

\begin{equation}
    \label{eq:IMEX}
    \vecU^{n + 1} = \vecF_I^{\Dt/2}\left(\vecF_E^{\Dt}\left(\vecF_I^{\Dt/2}\left(\vecU^n\right)\right)\right)
\end{equation}

\indent The implicit term $\vecF_I^{\Dt/2}$ consists of a Crank-Nicolson discretization:

\begin{equation}
    \label{eq:IMEX_implicit}
    \begin{aligned}
    \vecU^{*} - \frac{\Dt}{4}\vecL(\vecU^{*}) = \vecU^{n} + \frac{\Dt}{4}\vecL(\vecU^{n})
    \end{aligned}
\end{equation}

\noindent which is solved for $U^{*}$ using the \quotes{semi-implicit treatment of the Coriolis term} described by \cite{temperton:1997}. The explicit term $\vecF_E^{\Dt}$ consists of a second-order Runge-Kutta scheme:

\begin{equation}
    \label{eq:IMEX_explicit}
    \begin{aligned}
        \ol{\vecU} &= \Dt \vecN(\vecU^*)\\
        \ol{\ol{\vecU}} &= \Dt\vecN(\vecU^* + \ol{\vecU}) \\
        \vecU^{**} &= \frac{1}{2} \left( \ol{\vecU} + \ol{\ol{\vecU}} \right)
    \end{aligned}
\end{equation}

\subsubsection{SL-SI-SETTLS}
\label{sec:SL_SI_SETTLS}

\indent Semi-Lagrangian semi-implicit schemes are popular numerical methods in atmospheric circulation modeling, being used by several operational numerical weather prediction models, \eg the IFS-ECMWF, the GSM-JMA and the Global Forecast System (GFS) at NCEP/NOAA, to cite only a few \cite{mengaldo_al:2018}. They use a semi-implicit (Crank-Nicolson) discretization of the linear terms and a semi-Lagrangian (SL) treatment of the nonlinear advection term. The principle of SL schemes is to combine the Eulerian and Lagrangian approaches for spatiotemporal PDEs, the former relying on a fixed spatial grid (\eg the IMEX method presented above), being conceptually simpler but usually restricted by Courant-Friedrichs-Lewy (CFL) stability conditions limiting the time step size. In contrast, the latter follows the trajectories of the fluid particles along time instead of using a fixed grid, which allows larger time steps but is much more complex in implementation. SL methods avoid these issues by using a Lagrangian approach at each time step $[t_n, t_{n+1}]$: it traces the trajectories of the particles arriving at each point of a fixed spatial grid at time $t_{n+1}$, performing a spatial interpolation at time $t_n$ to retrieve the departure values. It simplifies SL compared to purely Lagrangian schemes but still allows the use of larger time steps than Eulerian methods \cite{staniforth:1991}.

\indent SL methods require to perform an estimation of the Lagrangian trajectories $(t,\vecx(t))$ along each time step $[t_n, t_{n+1}]$, with several approaches being proposed in the literature. We consider here the Stable Extrapolation Two-Time-Level Scheme (SETTLS) proposed by \cite{hortal:2002}, in which the trajectories are computed iteratively based on an estimation of the velocity field on the intermediate timestep: 

\begin{equation}
    \label{eq:SETTLS_trajectories}
    \vecx_d^{k+1} = \vecx_a - \frac{\Dt}{2} \vecV\left(t_n + \frac{\Dt}{2}\right) \approx \vecx_a - \frac{\Dt}{2} \left( [2\vecV(t_n, \vecx_d) - \vecV(t_{n-1}, \vecx_d) + \vecV(t_n, \vecx_a) ] \right)
\end{equation}

\noindent where $\Dt := t_{n+1} - t_n$, $\vecx_a$ denotes the arrival point of the trajectories $(t, \vecx(t))$ along $[t_n, t_n+1]$, corresponding to the fixed spatial grid, and $\vecx_d$ denotes the departure point. The SL-SI-SETTLS thus reads

\begin{equation}
    \label{eq:SL_SI_SETTLS}
    \frac{\vecU^{n+1} - \vecU^n_*}{\Dt} = \frac{1}{2}\left(\vecL \vecU^{n+1} + \vecL \vecU^n_*\right) + \frac{1}{2} \left([2\vecN_R(\vecU^n) - \vecN_R(\vecU^{n-1}) ]_* + \vecN_R(\vecU^{n})\right)
\end{equation}

\noindent where the subscript $*$ denotes interpolation to the departure point $\vecx_d$. In \eqref{eq:SETTLS_trajectories} and  \eqref{eq:SL_SI_SETTLS}, the terms in brackets are linear extrapolations to $t_{n+1}$, which are then averaged with the known values at $t_n$ to obtain an estimation at ${t_{n+1/2}} := t_n + \Dt/2$.

\indent We notice that this two-step time-stepping scheme, with the solution at time $t_{n+1}$ depending on the two preview times $t_n$ and $t_{n-1}$, poses additional challenges in terms of practical implementation of the parallel-in-time methods. Therefore, we use a modified version of the method in the PinT framework as discussed in Section \ref{subsec:pint_implementation}.

\section{Parallel-in-time methods}
\label{sec:PinT}

\indent Following the classification proposed by \cite{bellen_zennaro:1989}, PinT methods comprise iterative schemes allowing to compute several time steps simultaneously, algorithms based on spatial domain decomposition methods, and direct methods parallelizing the scheme used for advancing each time step. Detailed reviews on PinT can be found in \cite{gander:2015, ong_schroder:2020}.

\indent In this work, we focus on two of the most popular PinT methods, Parareal and MGRIT. Both can be interpreted as iterative, predictor-corrector algorithms, based on the simultaneous use of coarse (low expensive) and fine (expensive) time-stepping methods, the former being computed sequentially along the entire temporal domain, whereas the latter is computed in parallel, \ie with several time steps being computed simultaneously. Parareal is a two-level scheme, using fine and coarse discretizations of the problem; on the other hand, MGRIT is a multilevel scheme, using more than two discretization levels, besides other generalizations \wrt Parareal, as explained below. The popularity of these methods can be explained by their non-intrusive character, allowing the user to implement their own time-stepping scheme as coarse and fine methods, and by their simple implementation and formulation, mainly in the case of Parareal. In this section, we briefly describe these methods and introduce some notation used in this paper. For this description, we consider the time-dependent system of ODEs

\begin{equation}
    \label{eq:base_problem}
    \pdert{\vecU}(t) = \vecf(t, \vecU(t)), \qquad \vecU(0) = \vecU_0, \qquad t \in [0, T]
\end{equation}

\noindent which can be obtained \eg via a spatial discretization of a PDE.

\subsection{Parareal}
\label{subsec:parareal}

\indent First developed by \cite{lions_al:2001} and presented as a predictor-corrector algorithm by \cite{baffico_al:2002, bal_maday:2002}, the Parareal method iteratively computes approximations to the solution of \eqref{eq:base_problem} using simultaneously two numerical schemes, $\phi_f$ and $\phi_c$, named respectively \emph{fine} and \emph{coarse} propagators. The method aims to provide these approximations with a smaller computational cost compared to the serial simulation of the accurate (thus expensive) fine propagator. Let $T_0 = 0, T_1, \dots, T_{N} = T$ be a discretization of the temporal domain $[0,T]$. We denote by $\phi_f(\vecU, T_n, T_{n+1})$ and $\phi_c(\vecU, T_n, T_{n+1})$ the propagation of $\vecU$ from $T_n$ to $T_{n+1}$, using respectively $\phi_f$ and $\phi_c$. The initial guess (solution at iteration $k = 0$) is provided by the serial simulation of the coarse propagator along the entire temporal domain:

\begin{equation*}
    \vecU^0_{n+1} = \phi_c(\vecU^0_n, T_{n}, T_{n+1}), \qquad n = 0, \dots, N - 1
\end{equation*}

\indent In the following iterations, Parareal computes

\begin{equation}
    \label{eq:parareal_description}
    \vecU^{k+1}_{n+1} = \phi_c(\vecU^{k+1}_n, T_n, T_{n+1}) + \phi_f(\vecU^{k}_n, T_n, T_{n+1}) - \phi_c(\vecU^{k}_n, T_n, T_{n+1}), \qquad n = 0, \dots, N -1 ,\qquad k\geq 0
\end{equation}

\noindent where $\vecU^k_n$ is an approximation to the solution of \eqref{eq:base_problem} at time $T_n$ and iteration $k$. Note that the only term on the right-hand side of \eqref{eq:parareal_description} that needs to be computed serially at each iteration is the first one, using the coarse propagator (which is supposed to be relatively cheap). The remaining terms, including the expensive, fine one, can be computed in parallel (\ie the propagations along each time slice $[T_n, T_{n+1}]$ can be distributed to different parallel processors) since they depend only on the solution of the previous iteration $k$, which has already been computed for every $n = 0, \dots, N$.

\indent In general, the coarse propagator is defined with a larger timestep compared to the fine one. This, however, is not necessary. The fine and coarse schemes can be defined by different spatial resolutions, numerical schemes, or integration orders, for example. Finally, it can be easily shown that the Parareal solution converges exactly to the fine solution (\ie the solution obtained via a serial simulation of the fine propagator $\phi_f$) in at most $N$ iterations \cite{gander:2008}; however, a much faster convergence is required in practice.

\subsection{MGRIT}
\label{subsec:MGRIT}

\indent The MGRIT algorithm is a multilevel, predictor-corrector iterative parallel-in-time method introduced by \cite{friedhoff_al:2013} and based on spatial multigrid methods \cite{ries_trottenberg:1979}. In the following paragraphs, we describe the main ideas behind MGRIT and its parameters that are relevant to the study proposed in this work; we refer the reader to \cite{falgout_al:2014, xbraid-package} for details that are omitted here.

\subsubsection{Basic definitions}

\indent Let $\nlevels$ be the number of temporal discretization levels, indexed from $0$ to $\nlevels - 1$, the former being the finest one and the latter the coarsest one. The time-stepping scheme $\phi_l$ in each level is defined by a timestep $\Dt_l$. As in Parareal, in general (but not necessarily), the coarse levels are defined such as to have a coarser temporal discretization \wrt the fine levels. We then define a coarsening factor $\cfactor \geq 1$ such that $\Dt_{l+1} = \cfactor \Dt_l, l = 0, \dots, \nlevels - 2$. For the sake of simplicity, we consider the same coarsening factor $\cfactor$ for all pairs of consecutive levels and also that the temporal discretization in each level is homogeneous.

\indent In the following, we consider a pair of consecutive levels $(l, l+1)$, called respectively \quotes{fine} and \quotes{coarse} levels, to introduce some notation and vocabulary of MGRIT, which can be easily extended for any pair of levels.

\begin{figure}[!htbp]
    \centering
    \includegraphics[scale = .9]{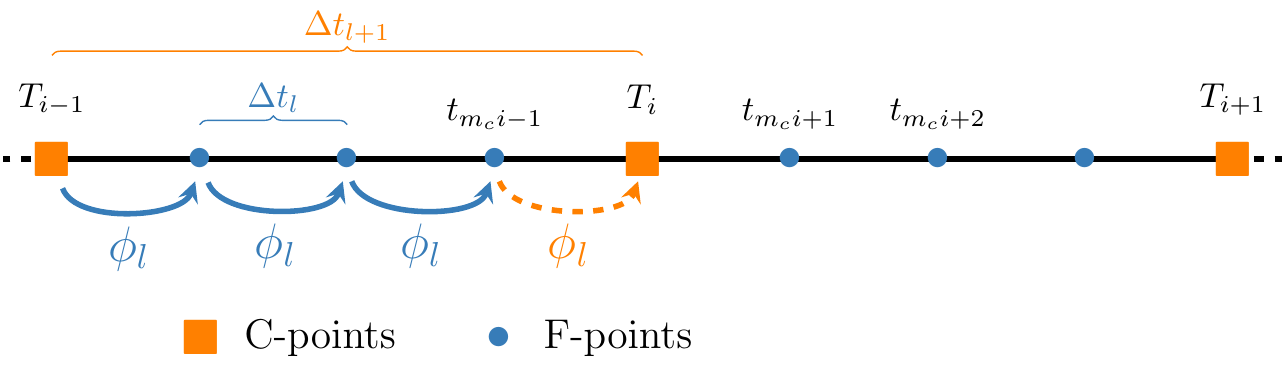}
    \caption{Basic definitions for the temporal discretization adopted in MGRIT. Bullets and squares represent, respectively, F-points and C-points. The sequence of solid arrows represents an F-relaxation. The dashed arrow represents a C-relaxation.}
    \label{fig:MGRIT_definitions}
\end{figure}

\indent The temporal domain is divided into $N_{l+1} := T / \Dt_{l+1}$ coarse time steps, defining the coarse time instances $T_i = i \Dt_{l+1}, i = 0, \dots, N_{l+1}$; and $N_l := T / \Dt_l = \cfactor N_{l+1}$ fine time steps, defining the time instances $t_j = j \Dt_l, j = 0, \dots, N_l$. Note that $T_i  = t_{\cfactor i}, i = 0, \dots, N_{l+1} $. The times $T_i, i = 0, \dots, N_{l+1}$ of the coarse temporal discretization are named \emph{C-points} and the fine time instances not present in the coarse discretization, \ie  $t_j = 1, \dots, N_l - 1, j \neq \cfactor i, i = 0, \dots, N_{l+1}$, are named \emph{F-points}.  These definitions are illustrated in Figure \ref{fig:MGRIT_definitions}.

\indent The numerical approximations for the solution of \eqref{eq:base_problem} are defined on the fine temporal grid, such that $\vecU^{(l)}_j$ and $\vecU^{(l)}_{\cfactor i}$ denote approximations at the F-points and C-points, respectively, with $i = 0, \dots, N_{l+1}$ and $j = 0, \dots, N_{l}, j \neq \cfactor i$.

\indent We call \emph{F-relaxation} (or \emph{F-relaxation sweep}) the update of the solution at all F-points between two C-points (\ie $t_j \in \ ]T_i, T_{i+1}[$), by using $\phi_l$ with the C-point value $\vecU^{(l)}_{\cfactor i}$ as an initial solution. Analogously, we call \emph{C-relaxation} the update of each C-point value $\vecU^{(l)}_{\cfactor i}$ by advancing one fine time step using $\phi_l$ with initial solution $\vecU^{(l)}_{\cfactor i -1}$, see Figure \ref{fig:MGRIT_definitions}. An update of the approximations on the fine level can be performed using a combination of F- and C-relaxations; in general, MGRIT is set to use \emph{FCF-relaxation}, \ie an F-relaxation followed by a C- and a second F- one, but one can define any $\text{F(CF)}^{\nrelax}$-relaxation, $\nrelax \geq 0$, including the simplest F-relaxation (which is used \eg by Parareal, corresponding to $\nrelax = 0$). In general, we can expect a faster MGRIT convergence with a more complex relaxation scheme, but with a higher numerical cost.

\subsubsection{Full approximation scheme}

\indent The idea of the MGRIT algorithm is to compute error corrections on the coarse grid that improve the solution obtained by relaxation on the fine grid. It is performed via the Full Approximation Scheme (FAS) \cite{brandt:1977}, which can be briefly described in five steps:

\begin{enumerate}
    \item Perform an $\text{F(CF)}^{\nrelax}$-relaxation on the fine temporal grid;
    \item Restrict the solution and its residual from the fine to the coarse grid at the C-points;
    \item Solve the so-called coarse grid equation;
    \item Use the solution of the coarse grid equation to compute a coarse grid error approximation;
    \item Correct the solution on the fine grid by injecting the computed error approximation from the coarse to the fine grid and performing a further F-relaxation sweep.
\end{enumerate}

\subsubsection{Multilevel scheme}

\indent In the two-level algorithm presented above, the third step is the only one that needs to be computed sequentially along the entire temporal domain. However, in the case of multilevel runs, this step is performed recursively, thus adding additional degrees of parallelism. In this case, the entire algorithm described above is applied for the pair of levels $(l+1, l+2)$ and so on. The order of recursive calls in each level is another parameter to be set in an MGRIT simulation, influencing the convergence speed and computational cost of the method. In this work, we consider F-cycles followed by one post-V-cycle for the order of recursive calls through levels, as defined in \cite{xbraid-package}.

\indent Finally, as in Parareal, an exact convergence towards the serial solution of the finest temporal discretization is obtained in a finite number of iterations. In a MGRIT simulation with $N_0$ time points on the finest level, coarsening factor $\cfactor$ and using an $\text{F(CF)}^{\nrelax}$-relaxation strategy, exact convergence is obtained in at most $N_0/((\nrelax+1)\cfactor)$ iterations \cite{falgout_al:2014}.

\subsubsection{Implementation}
\label{subsec:pint_implementation}

\indent The numerical results presented in Section \ref{sec:numerical_tests} were obtained using the XBraid library \cite{xbraid-package}, which implements MGRIT combined with user-defined time-stepping methods. The Parareal simulations were also performed using XBraid, by setting the parameters $\nlevels = 2$ and $\nrelax = 0$.

As a major feature, XBraid allows us to easily incorporate MGRIT into an operational code. Notably, every parallelization is made in a \quotes{black-box} fashion inside the library itself, with the user only needing to define \emph{what} is communicated between processors, but not \emph{when} or \emph{how}. This, however, may not be convenient in the cases of two-step time-stepping schemes (such as SL-SI-SETTLS, described in Section \ref{sec:SL_SI_SETTLS}), since two consecutive time steps $[t_{n-1}, t_{n}]$ and $[t_{n}, t_{n+1}]$, both contributing to the solution at $t_{n+1}$, may be treated by different parallel processors. In \cite{schmitt_al:2018}, a modified version of Parareal is proposed in order to deal with it. In XBraid, this issue can be solved on the coarsest level by including the solution at $t_{n-1}$ in the solution to be communicated, since the simulation is sequential at this level. In all other levels, however, time steps are computed in parallel, with no predefined order, making it impracticable to define this extra communication in a coherent way. An alternative approach would be to replace the two-step scheme with a one-step one, with the solution at $t_{n-1}$ replaced by the one at $t_{n}$. Preliminary studies using the Gaussian bumps test case introduced in Section \ref{sec:numerical_tests} were conducted comparing the MGRIT solution under two scenarios: (i) with the alternative one-step scheme applied in all time steps of all levels, except for the coarsest one, in which the original two-step scheme is applied; (ii) with the one-step scheme being applied in all time steps of all levels, including the coarsest one. These results show that the chosen procedure has little influence on the convergence and stability of MGRIT. Therefore, and to avoid extra communications on the coarsest level, the following procedure is considered for adapting SL-SI-SETTLS to MGRIT:

\begin{itemize}
    \item On the coarsest level, the original two-step SETTLS scheme is used in all time steps, except for the \emph{temporal boundaries} between processors, \ie if $[t_{n-1}, t_{n}]$ and $[t_{n}, t_{n+1}]$ are treated by different processors, in which case the one-step alternative is adopted;
    \item On all other levels using SL-SI-SETTLS, the one-step alternative is used in all time steps.
\end{itemize}

\indent We remark that, except for the simulations performed for evaluating wall-clock times, all the simulations presented in Section \ref{sec:numerical_tests} were run using a single processor in time, which implies that the original SL-SI-SETTLS scheme is applied on the entire coarsest level. \added[id=R3]{Notably, in the two-level configurations using SL-SI-SETTLS on the coarse level and another scheme on the fine one, the PinT simulation is entirely conducted without modification of SL-SI-SETTLS. Finally, we notice that the analytical stability study presented in the next section, which indicates stability issues of Parareal and MGRIT using SL-SI-SETTLS on the coarse level (which are confirmed in the numerical simulations presented in Section \ref{sec:numerical_tests}), are conducted considering the original SL-SI-SETTLS.}

\section{Stability study}
\label{subsec:stability}

\indent We perform in this section a stability study of Parareal and MGRIT depending on given parameters and given fine and coarse time-stepping methods. This study follows the procedure presented in \cite{cox_matthews:2002}. We consider the ODE

\begin{equation}
    \label{eq:ODE}
    \dot{u} = \lambdaL u + \vecN(u)
\end{equation}

\noindent which is obtained \eg by solving \eqref{eq:swe_sphere} with a spectral method, in which case $\lambdaL \in \complex$ is a wavenumber mode of $\vecL$ and $u$ is a spectral coefficient of $\vecU$. Problem \eqref{eq:ODE} is linearized by considering $\vecN(u) = \lambdaN u$, with $\lambdaN \in \complex$. Thus,

\begin{equation}
    \label{eq:linearized_ODE}
    \dot{u} = \lambdaL u + \lambdaN u
\end{equation}

\indent As pointed out by \cite{cox_matthews:2002}, a full stability study of a numerical scheme to \eqref{eq:linearized_ODE} should be performed in a four-dimensional space, \ie as a function of the real and imaginary parts of $\xiL := \lambdaL \Dt$ and $\xiN := \lambdaN \Delta t$. Since the SWE equations on the rotating sphere are characterized by the propagation of purely imaginary wavenumber modes, we perform this study in the function of $Re(\xiN)$ and $Im(\xiN)$ for fixed values of $\xiL$. Namely, we are interested in the cases $\xiL \approx 0 $ (Rossby modes) and positive and negative purely imaginary values of $\xiL$ (inertia-gravity wave modes).

\indent Let us first determine some physically meaningful values of $\xiL$ for this stability study, which can be easily obtained via a simplification of the linear operator $\vecL$ defined in \eqref{eq:swe_sphere}. By considering the Coriolis parameter $f$ to be constant (which corresponds to the so-called $f$-plane approximation), we can write 

\begin{equation*}
    \vecL(\vecU) = \vecL\vecU = \matthree{0}{0}{-\ol{\Phi}}{0}{0}{-f}{-\nabla^2}{f}{0}\vecthree{\Phi}{\xi}{\delta}
\end{equation*}

\noindent which, in the spherical harmonics spectral space, reads

\begin{equation}
    \label{eq:linear_operator_spectral}
    \mc{L}\vecU_{m,n} = \matthree{0}{0}{-\ol{\Phi}}{0}{0}{-f}{-\dfrac{n(n+1)}{a^2}}{f}{0}\vecthree{\Phi_{m,n}}{\xi_{m,n}}{\delta_{m,n}}
\end{equation}

\indent The three eigenvalues of $\mc{L}$ are

\begin{equation*}
    \lambda_{\vecL} = 0, \qquad \lambda_{\vecL} = \pm i\sqrt{f^2 + \frac{n(n+1)}{a^2}\ol{\Phi}}
\end{equation*}

\indent In the numerical simulations presented in Section \ref{sec:numerical_tests}, the mean geopotential field has an order of magnitude of $10^5 \text{m}^2\text{s}^{-2}$. We then consider here $\ol{\Phi} = 10^5 \text{m}^2\text{s}^{-2}$. The radius of the Earth is $a \approx 6371.22\times10^3 \text{m}$ and the Coriolis parameter in a $f$-plane approximation reads $f = 2\Omega \approx 2 \times 7.292 \times 10^{-5} \text{s}^{-1}$, where $\Omega$ is the Earth's angular velocity. Since $\ol{\Phi}/a^2 \gg f^2$, we neglect this term in the approximation of $\lambda_{\vecL}$ and we plot the stability regions for fixed integer multiples of  $\tilde{\xi}_{\vecL} := i\sqrt{\ol{\Phi}/a^2} \approx 2.5\times10^{-4}i$, where the scaling of this quantity can be interpreted as a consequence of the choice of wavenumber $n$ and/or the choice of timestep $\Dt$.

\indent For a given numerical scheme for \eqref{eq:ODE}, we compute the amplification factor $A := |u^{n+1}| / |u^n|$ and plot its stability region (\ie where $|A| \leq 1$) on the $Re(\xiN)-Im(\xiN)$ plane. We are mainly interested in the overlapping between the stability region and the imaginary axis.

\indent We begin by presenting the stability plots for IMEX and SL-SI-SETTLS individually (\ie in a serial framework), where we introduce some specific details of the stability analysis of each method. Then, we present the stability study of Parareal presented by \cite{staff_ronquist:2005} and we extend it to the MGRIT framework in a two-level configuration with an arbitrary relaxation strategy.

\subsection{Stability of the time-stepping schemes individually}

\subsubsection{Stability of IMEX}

\indent Applying \eqref{eq:IMEX} to the linearized problem \eqref{eq:ODE}, the implicit step \eqref{eq:IMEX_implicit} with a half timestep $\Dt/2$ has the amplification factor

\begin{equation*}
    A_{\text{IMEX-I}} = \frac{4 + \xiL}{4 - \xiL}
\end{equation*}

\noindent and, for the explicit step with a full timestep,

\begin{equation*}
    A_{\text{IMEX-E}} =  \frac{\xiN + \xiN(1+\xiN)}{2}  = \frac{\xiN(2 + \xiN)}{2}
\end{equation*}

\indent Thus, the amplification factor of IMEX reads

\begin{equation*}
    A_{\text{IMEX}} = A_{\text{IMEX-I}} A_{\text{IMEX-E}} A_{\text{IMEX-I}} = \left(\frac{4 + \xiL}{4 - \xiL}\right)^2  \frac{\xiN(2 + \xiN)}{2}
\end{equation*}

\subsubsection{Stability of SL-SI-SETTLS}

\indent In the stability study of semi-Lagrangian schemes, the amplification factors also depend on the spatial wavenumber (see \eg \cite{durran:2010}). Indeed, the amplification factor $A$ is more rigorously defined by considering a solution under the form $u(x,t_n) = A^n e^{i\kappa x}$, where $\kappa$ is the wavenumber. In the case of Eulerian schemes (\eg IMEX), all terms are evaluated in the same spatial point and the exponential term vanishes, leading to $A = |u^{n+1}|/|u^n|$, which is not the case in the semi-Lagrangian framework. We define $s := x_j - x_*$ as the distance between the grid point $x_j$ and the departure point $x_*$ (at $t_n$) of the trajectory arriving at $x_j$ at $t_{n+1}$. Replacing in \eqref{eq:SL_SI_SETTLS}, we obtain

\begin{equation*}
    A^{n+1}e^{i\kappa x_j} - A^{n}e^{i\kappa x_*} = \frac{\xiL}{2}\left( A^{n+1}e^{i\kappa x_j} + A^{n}e^{i\kappa x_*} \right)  + \frac{\xiN}{2} \left( \left[2 A^{n} - A^{n-1} \right]e^{i\kappa x_*} + A^n e^{i\kappa x_j}\right)
\end{equation*}

\noindent which can be simplified and rearranged to 

\begin{equation}
    \label{eq:stability_SL_SI_SETTLS}
    A_{\text{SL-SI-SETTLS}}^2 \left( 1 - \frac{\xiL}{2} \right) - A_{\text{SL-SI-SETTLS}} \left( e^{-i\kappa s} \left(1 +  \frac{\xiL}{2} + \xiN \right) + \frac{\xiN}{2} \right) + \frac{\xiN}{2}e^{-i \kappa s} = 0
\end{equation}

\indent Eq. \eqref{eq:stability_SL_SI_SETTLS} has two roots $A_{\text{SL-SI-SETTLS}}^+$ and $A_{\text{SL-SI-SETTLS}}^-$, and the stability region of SL-SI-SETTLS is the intersection of their stability regions. Moreover, the stability region depends on the spatial wavenumber $\kappa s$. As shown by \cite{hortal:2002}, SL-SI-SETTLS has improved stability properties since the intersection of the stability regions for all $\kappa s$ is not empty, which is not the case of previous SL-schemes used in IFS-EMCWF. In this work, we compute the stability region of SL-SI-SETTLS by intersecting the stability regions of $\kappa s \in [0, 2\pi]$ with a step $\pi/10$.

\begin{figure}[!htbp]
    \centering
    \includegraphics[scale = .5]{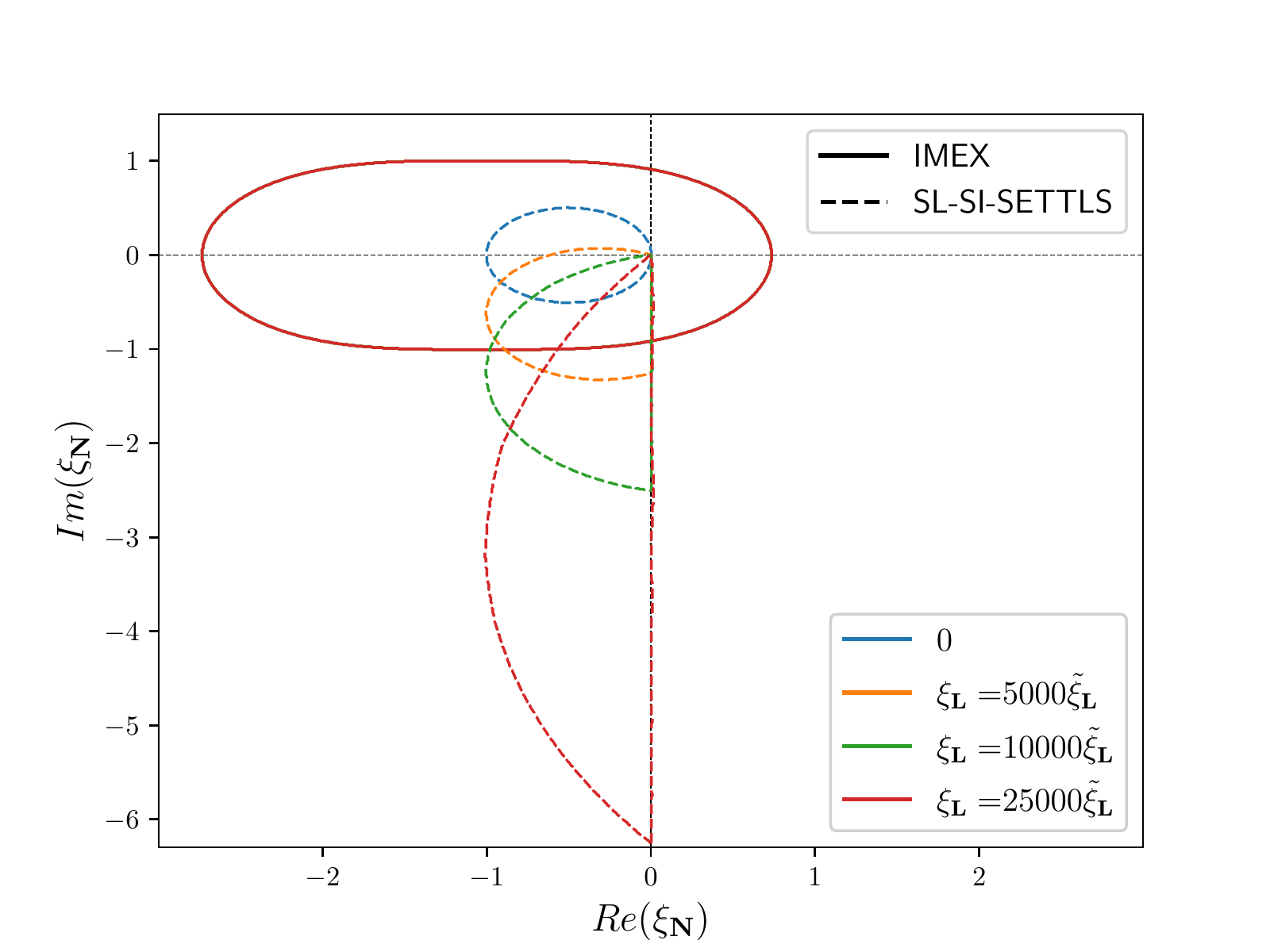}
    \caption{Stability regions of IMEX and SL-SI-SETTLS applied to \eqref{eq:ODE} as a function of $\xiN$ for fixed values of $\xiL$. The contours for IMEX visually coincide.}
    \label{fig:stability_schemes}
\end{figure}

\indent Figure \ref{fig:stability_schemes} compares the stability regions of both schemes for $\xiL \in \{0, 5\times10^3 \tilde{\xi}_{\vecL}, 10^4 \tilde{\xi}_{\vecL}, 2.5\times10^4 \tilde{\xi}_{\vecL} \}$. For negative $\xiL$, the plots are symmetric \wrt the $Re(\xiN)$-axis. We note that the stability regions of IMEX do not vary with $\xiL \in i\reals$ since, in this case, $|(4+\xiL)/(4-\xiL)| = 1$, indicating that $A_{\text{IMEX}} = \xiN ( 2+ \xiN)/2$. On the other hand, the stability region of SL-SI-SETTLS depends on $\xiL$, and for larger values of the linear mode, the method outperforms IMEX, notably with a larger intersection between the stability region and the imaginary axis. 

\subsection{Stability of Parareal}

\indent We briefly present the stability analysis developed by \cite{staff_ronquist:2005} for the parareal iteration

\begin{equation}
    \label{eq:parareal2}
    u^{k}_n = \phi_c(u^k_{n-1}) + \phi_f(u^{k-1}_{n-1}) - \phi_c(u^{k-1}_{n-1})
\end{equation}

\noindent using the fine and coarse timestepping methods $\phi_c$ and $\phi_f$ with respective timesteps $\Dt_f$ and $\Dt_c = \cfactor \Dt_f$. The parareal solution $u^k_n$ at iteration $k$ is computed at time $t_n := n \DT$, with $\Dt/\Dt_f =: N_f$ and $\Dt/\Dt_c =: N_c$.

\indent Let $A_f := A_f(\Dt_f)$ and $A_c := A_c(\Dt_c)$ be the stability functions of $\phi_f$ and $\phi_c$, respectively. Then, \eqref{eq:parareal2} can be written as 

\begin{equation}
    \label{eq:parareal_stability}
    u^{k}_n = A_c^{N_c}u^k_{n-1} + A_f^{N_f}u^{k-1}_{n-1} - A_c^{N_c}u^{k-1}_{n-1} = A_c^{N_c}u^k_{n-1} + (A_f^{N_f} -  A_c^{N_c})u^{k-1}_{n-1}
\end{equation}

\begin{figure}[!htbp]
    \begin{subfigure}{.5\linewidth}
        \includegraphics{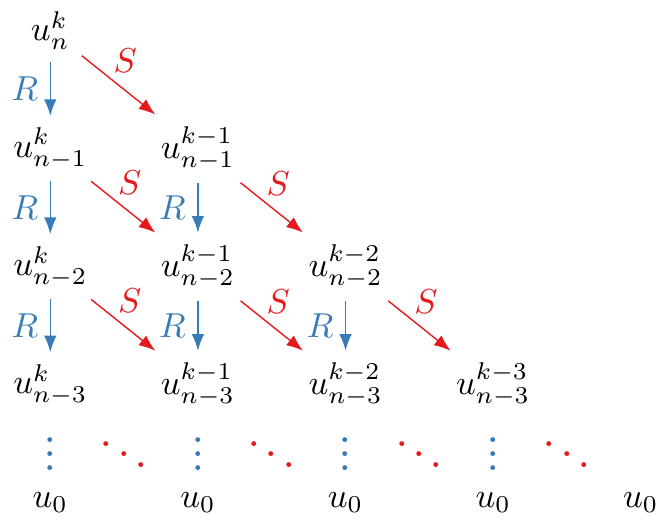}
        \caption{Parareal, adapted from \cite{staff_ronquist:2005}}
        \label{fig:recurrence_parareal}
    \end{subfigure}
    \begin{subfigure}{.5\linewidth}
        \includegraphics{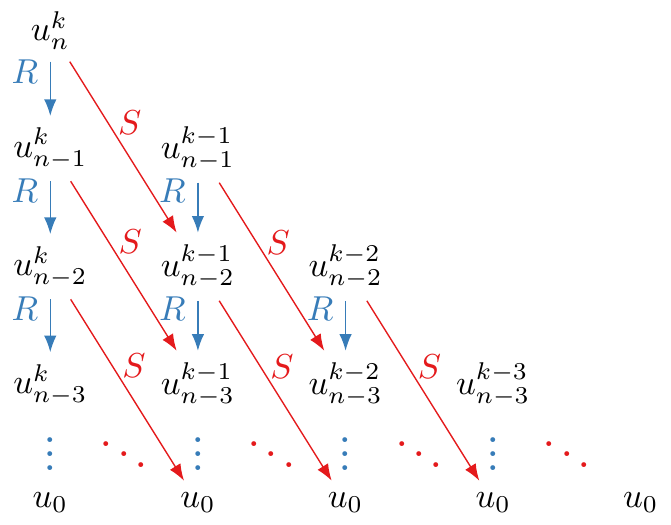}
        \caption{Two-level MGRIT with $\nrelax = 1$}
        \label{fig:recurrence_mgrit_nrelax1}
    \end{subfigure}
    \caption{Illustration of the recurrence relations of: (a) the Parareal iteration $u^k_n = Ru^k_{n-1} + Su^{k-1}_{n-1}$, with $R := A_c^{N_c}$ and $S := A_f^{N_f} -  A_c^{N_c}$; and (b) the two-level MGRIT iteration with $\nrelax = 1$, $u^k_n = Ru^k_{n-1} + Su^{k-1}_{n-2}$, with $R := A_c^{N_c}$ and $S := A_f^{2N_f} -  A_c^{2N_c}$.}
    \label{fig:recurrence}
\end{figure}

\indent Eq. \eqref{eq:parareal_stability} is a recurrence relation leading to

\begin{equation}
    \label{eq:recurrence_parareal}
    u^{k}_n = \left( \sum_{i = 0}^k \vectwo{n}{i} \left(A_f^{N_f} -  A_c^{N_c}\right)^i \left(A_c^{N_c}\right)^{n-i}\right)u_0
\end{equation}

\noindent where $u_0$ is the initial solution at $t = 0$. The coefficients in \eqref{eq:recurrence_parareal} are the binomial coefficients of Pascal's triangle, which can be easily identified in Figure \ref{fig:recurrence_parareal}, adapted from \cite{staff_ronquist:2005}.

% \begin{figure}[!htbp]
%     \centering
%     \includegraphics{pint/figures/stability/TikZfigures-recurrence_parareal.pdf}
%     \caption{Illustration of the recurrence relation of the Parareal iteration $u^k_n = Ru^k_{n-1} + Su^{k-1}_{n-1}$, with $R := A_c^{N_c}$ and $S := A_f^{N_f} -  A_c^{N_c}$. Adapted from \cite{staff_ronquist:2005}. }
%     \label{fig:recurrence_parareal}
% \end{figure}

\indent From \eqref{eq:recurrence_parareal}, we define the Parareal stability function

\begin{equation}
    \label{eq:parareal_stability_function}
    A_{\text{parareal}}(n, k, N_f, N_c, \phi_f, \phi_c) := \sum_{i = 0}^k \vectwo{n}{i} \left(A_f^{N_f} -  A_c^{N_c}\right)^i \left(A_c^{N_c}\right)^{n-i}, \qquad n \geq k
\end{equation}

\begin{figure}[!htbp]
    \begin{subfigure}{.5\linewidth}
    \centering
    \includegraphics[scale = \scaleFig]{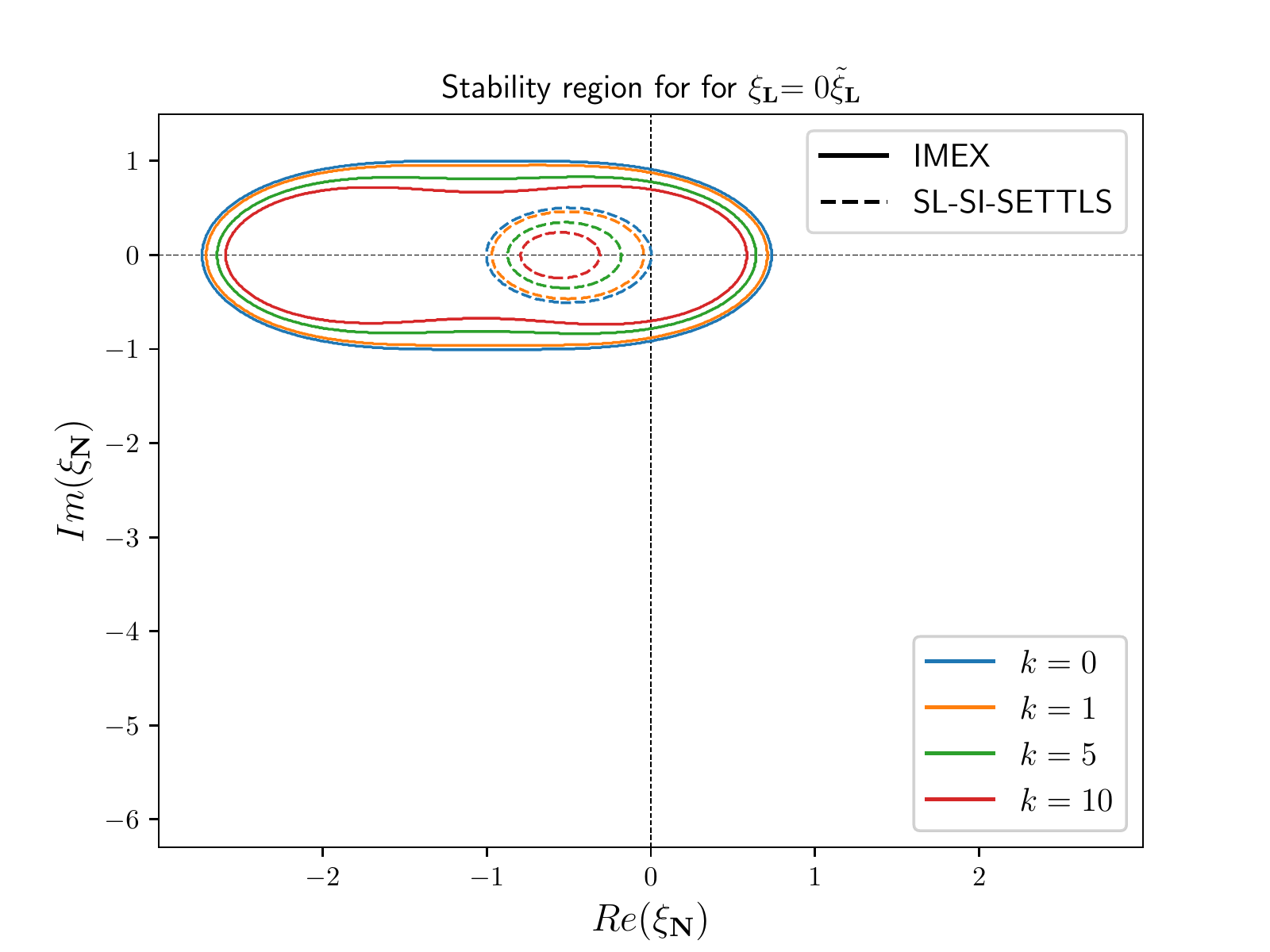}
    \caption{$\xiL = 0$}
    \end{subfigure}
    \begin{subfigure}{.5\linewidth}
    \centering
    \includegraphics[scale = \scaleFig]{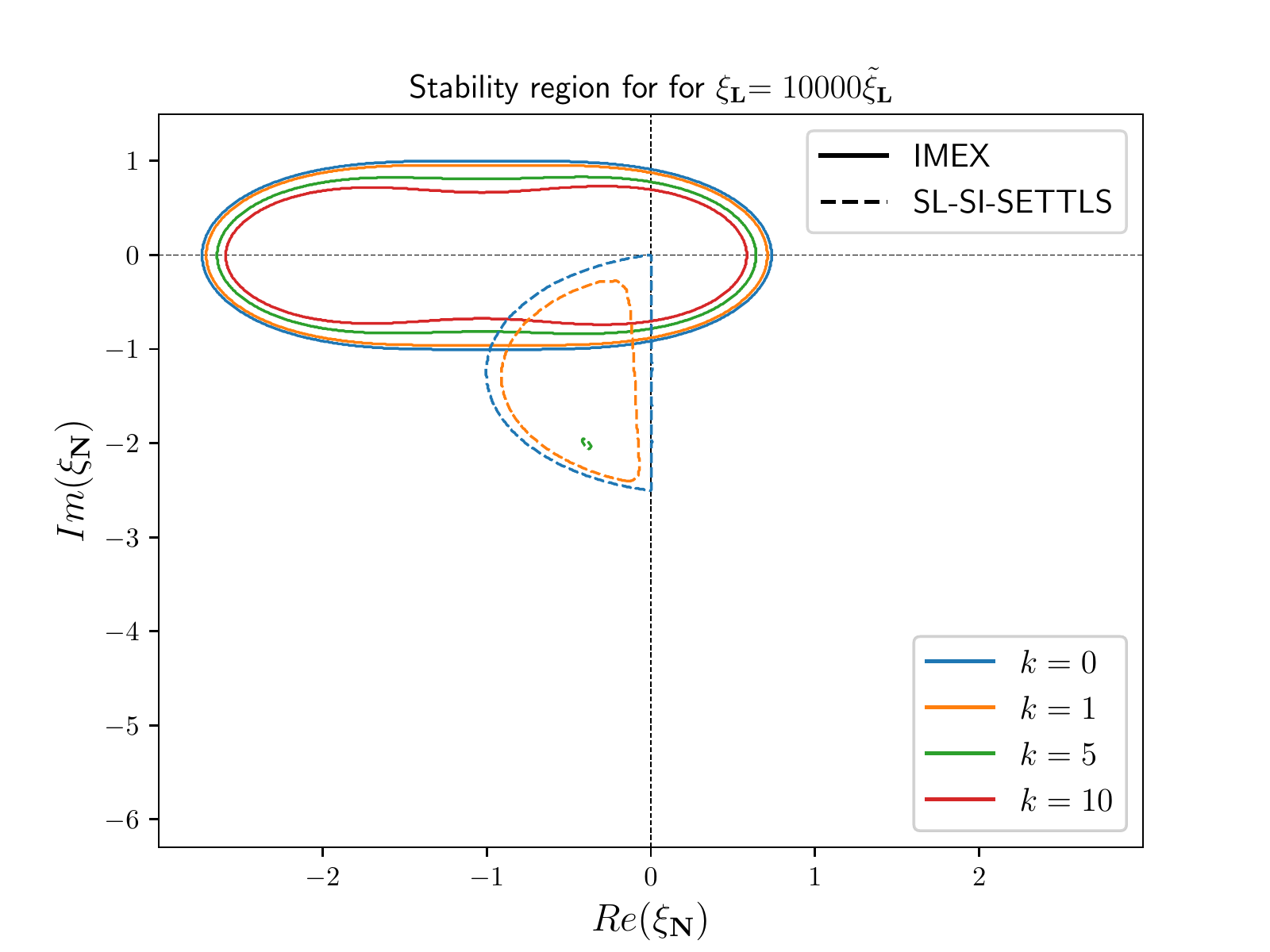}
    \caption{$\xiL = 10^4 \tilde{\xi}_{\vecL}$}
    \end{subfigure}
    \caption{Stability regions of Parareal (with $\cfactor = 2$, $n = 100$, $k \in \{0, 1, 5, 10\}$) applied to \eqref{eq:ODE} as a function of the coarse timestepping method and $\xiN$ for fixed values of $\xiL$. In all configurations, the fine time-stepping method is IMEX. In the case $\xiL = 10^4 \tilde{\xi}_{\vecL}$, the stability region of SL-SI-SETTLS is empty for $k = 10$.}
    \label{fig:stability_parareal_iterations}
\end{figure}

\indent Let us study the behavior of $A_{\text{parareal}}$ along iterations and as a function of the coarse time-stepping method. We fix $n = 100$, $\cfactor = 2$, $N_c = 1$, and IMEX as a fine scheme and we plot the stability contours at iterations $k = 0, 1, 5, 10$ in Figure \ref{fig:stability_parareal_iterations}. In all cases (for both choices of coarse time-stepping method and for all chosen values of $\xiL$), we observe a decrease of the stability region along iterations. Also, for both coarse schemes, the stability region is not larger than the ones using the schemes alone (which corresponds to $k = 0$). A much more notable stability loss is observed in the configurations using SL-SI-SETTLS as a coarse method, whose stability contours no longer intersect the $Im(\xi_N)$-axis for $k \geq 1$. Also, in the case $\xiL = 10^4 \tilde{\xi}_{\vecL}$, in which the serial SL-SI-SETTLS outperforms IMEX, the stability regions vanish (or almost) in the Parareal framework after few iterations. On the other hand, the configurations using IMEX as a coarse method are able to preserve quite well their shape and only suffer a small reduction in size.

\subsection{Stability of MGRIT}

\indent We generalize the Parareal stability analysis presented above to the MGRIT framework. We consider here the two-level case with $\text{F(CF)}^{\nrelax}$-relaxation and we study the influence of $\nrelax$ on the stability contours using the fact that the two-level MGRIT can be seen as an overlapping variant of Parareal, with the parallel steps performed at each iteration being computed along  $\nrelax + 1$ time slices \cite{gander_al:2018}:

\begin{equation}
    \label{eq:parareal_FCF_nrelax}
    u^k_n = \phi_c(u^k_{n-1}) + \phi_f^{\nrelax + 1}\left(u^{k-1}_{n-(\nrelax+1)}\right) - \phi_c\phi_f^{\nrelax}\left(u^{k-1}_{n-(\nrelax+1)}\right)
\end{equation}

\indent In terms of the fine and coarse stability functions:

\begin{equation}
    \label{eq:parareal_FCF_nrelax_B}
    u^k_n = \left[A_c^{N_c}u^k_{n-1} + \left(A_f^{(\nrelax+1)N_f} - A_c^{N_c} A_f^{\nrelax N_f} \right) u^{k-1}_{n-(\nrelax+1)} \right]
\end{equation}

\noindent which can be rewritten as

\begin{equation*}
    u^k_n = A_{\text{MGRIT}}(n, k, \phi_f, \phi_c, \nlevels = 2, \nrelax) u_0   
\end{equation*}

\noindent with

\begin{equation}
    \label{eq:MGRIT_stability_function_FCFnrelax_twolevels}
    \begin{gathered}
    A_{\text{MGRIT}}(n, k, N_f, N_c, \phi_f, \phi_c, \nlevels = 2, \nrelax) := \\ \sum_{i = 0}^{\lfloor k/(\nrelax+1) \rfloor} \vectwo{n-i\nrelax}{i} \left(A_f^{N_f\nrelax} \left(A_f^{N_f} - A_c^{N_c}  \right)\right)^i \left(A_c^{N_c}\right)^{n-i(\nrelax+1)}, \qquad n \geq \frac{k}{\nrelax+1}
    \end{gathered}
\end{equation}

\noindent being the stability function of MGRIT, where $\lfloor\cdot \rfloor$ is the floor function. It is easy to check that \eqref{eq:MGRIT_stability_function_FCFnrelax_twolevels} reduces to \eqref{eq:parareal_stability_function} in the case $\nrelax = 0$ (Parareal).

\indent The derivation of \eqref{eq:MGRIT_stability_function_FCFnrelax_twolevels} is analogous to the Parareal case. The two-level MGRIT iteration \eqref{eq:parareal_FCF_nrelax_B} can be identified with a Pascal's triangle, but with a \quotes{jump} of $\nrelax$ rows in all \quotes{diagonal relations}. Figure \ref{fig:recurrence_mgrit_nrelax1} illustrates it in the case $\nrelax = 1$.

\begin{figure}[!htbp]
    \begin{subfigure}{.5\linewidth}
    \centering
    \includegraphics[scale = \scaleFig]{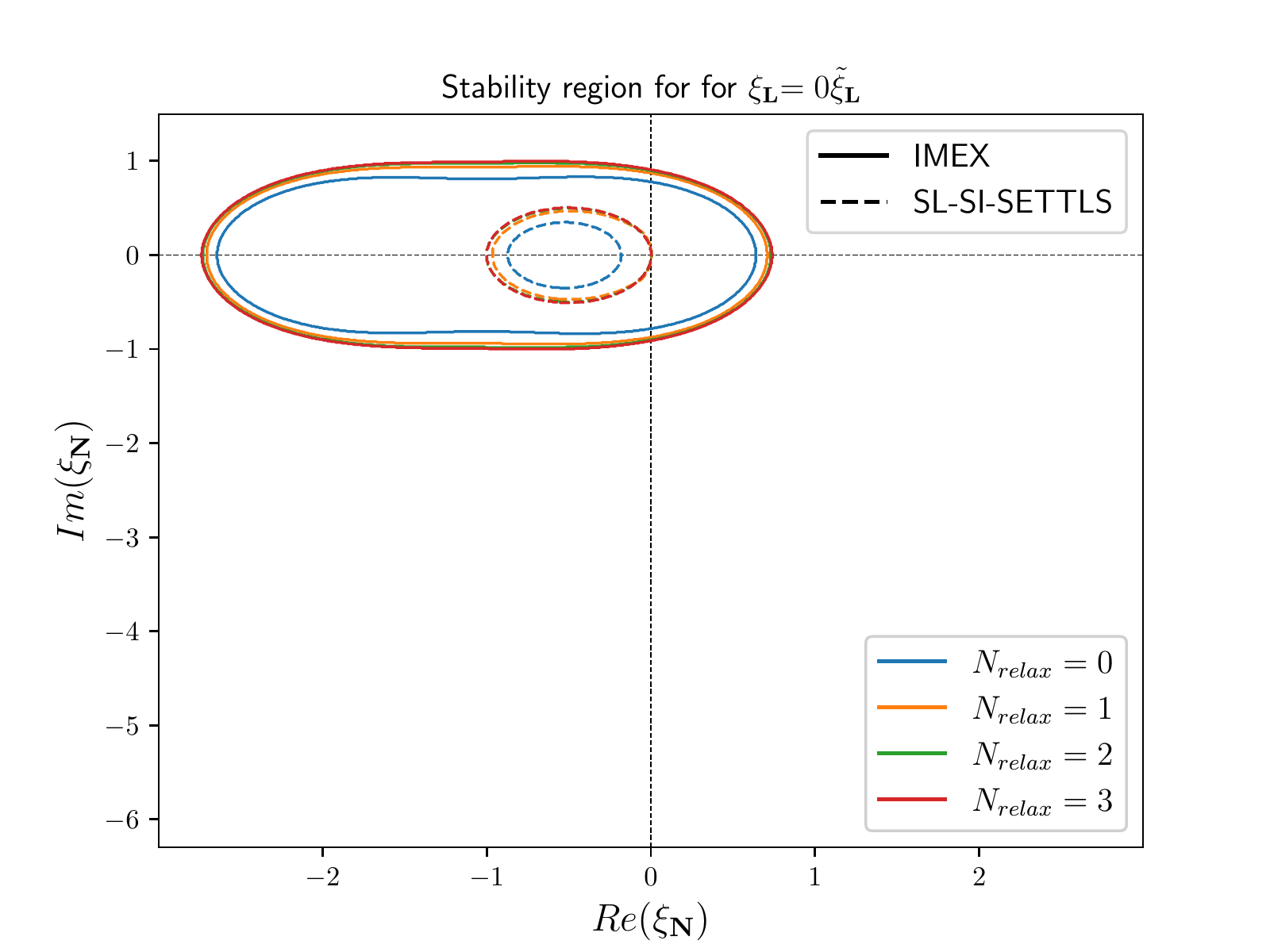}
    \caption{$\xiL = 0$}
    \end{subfigure}
    \begin{subfigure}{.5\linewidth}
    \centering
    \includegraphics[scale = \scaleFig]{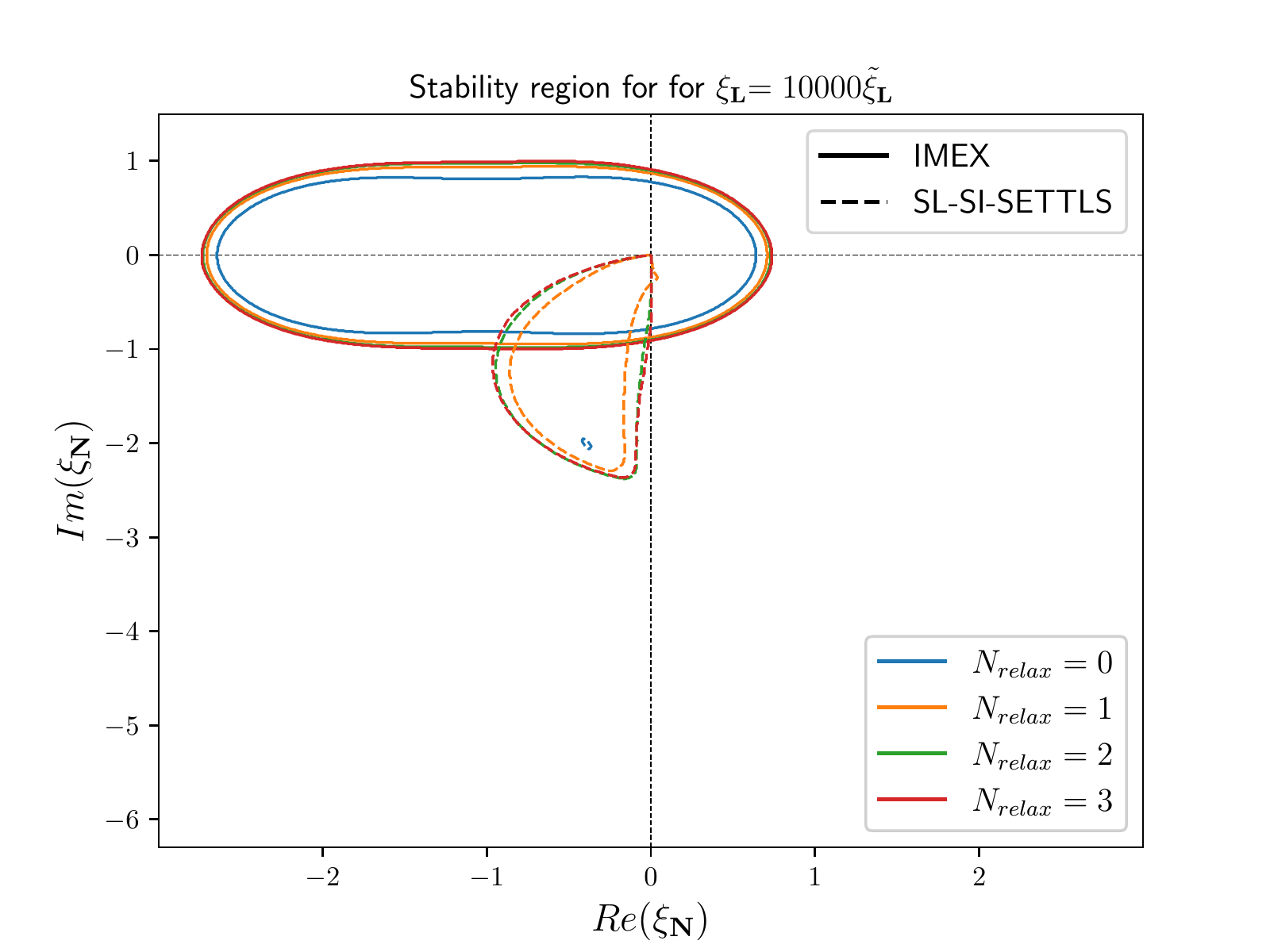}
    \caption{$\xiL = 10^4 \tilde{\xi}_{\vecL}$}
    \end{subfigure}
    \caption{Stability regions of MGRIT (with $k = 5$, $\cfactor = 2$, $n = 100$ and $\nrelax \in \{0, 1, 2, 3\}$) applied to \eqref{eq:ODE} as a function of the coarse timestepping method and $\xiN$ for fixed values of $\xiL$. In all configurations, the fine time-stepping method is IMEX. In almost all cases, the contour plots for $\nrelax = 2$ and $\nrelax = 3$ visually coincide.}
    \label{fig:stability_mgrit_nrelax}
\end{figure}

\indent In Figure \ref{fig:stability_mgrit_nrelax}, we plot the stability regions for fixed MGRIT parameters ($k = 5$, $n = 100$, $\cfactor = 2$) and $\nrelax \in \{0, 1, 2, 3\}$. We observe that the stability regions increase when $\nrelax$ increases, mainly from $\nrelax = 0$ to $\nrelax = 1$. In most cases, the stability regions of the serial schemes are almost or totally recovered with $\nrelax \geq 2$; however, for large $\xiL$ ($\xiL = 10^4 \tilde{\xi}_{\vecL}$), it is not enough to improve the stability region of MGRIT using SL-SI-SETTLS as coarse method, remaining smaller than the stability region of the serial scheme.

\added[id = R2]{We highlight that the small stability regions in the configurations using SL-SI-SETTLS may be linked to the fact that the stability region is defined as the intersection between the stability regions of all wavenumbers $\kappa$. In the Parareal and MGRIT framework, each of these regions is reduced depending on the parametric choice; thus, their intersection is drastically reduced compared to the stability region of the serial scheme. It does not imply that SL schemes will necessarily suffer from this issue; if the amplification factor can be written in the form $A(\xiL, \xiN, e^{-i \kappa s}) = e^{-i\kappa s}\tilde{A}(\xiL, \xiN)$, then the region where $|A| \leq 1$ will be the same for all wavenumbers, and its reduction in the PinT framework may be less remarkable. This is the case, for instance, of semi-Lagrangian exponential methods \cite{peixoto_schreiber:2019}, which also present a better stability behavior in PinT simulations, as will be presented in future work.}

\section{Numerical tests}
\label{sec:numerical_tests}

\indent We consider here two test cases for evaluating the application of Parareal and MGRIT combined with IMEX or SL-SI-SETTLS for the parallel-in-time integration of the SWE on the rotating sphere. For each test, we perform a convergence study, in which we study the convergence speed and stability as a function of selected Parareal and MGRIT parameters (namely the number of levels $\nlevels$, the coarsening factor $\cfactor$, the relaxation strategy $\nrelax$ and the spectral resolution $\Mcoarse$ of the coarse levels), and we evaluate the computational times and speedups \wrt the parallelization. Moreover, in the first test case, we study the influence of the artificial (hyper)viscosity on the convergence and stability of the temporal parallelization, and we use the obtained conclusions for setting up the viscosity parameters in the second test case.

\indent We are also interested in comparing the performance of the PinT methods using IMEX or SL-SI-SETTLS. Since the finest level is commonly discretized with a small time step, it provides an accurate solution to the problem and is not critical in terms of stability, which is not the case for the coarser levels, in which larger time steps are used. Therefore, we simplify this study by fixing the time-stepping scheme on the finest level ($l = 0$), namely IMEX, and considering the same time-stepping scheme (IMEX or SL-SI-SETTLS) on all coarse levels ($l > 0$).

\indent \added[id=R3]{The choices of parameters and test cases presented in this section represent an effort towards the operational application of PinT methods in atmospheric modeling, even though their application to more complex and realistic problems would evidently require more detailed studies. As done in the analytical stability study presented in Section \ref{subsec:stability}, in which the modes $\xiL$ were chosen based on realistic physical parameters, we consider here, for instance, meaningful temporal and spatial discretization sizes and viscosity coefficient values. Moreover, the two test cases represent quite complex dynamics, mainly the second one, which is a standard and challenging benchmark in atmospheric modeling research.}

\indent We first define the errors and computational time measures used for evaluating the PinT performance:

\paragraph{Error definition}

\indent As proposed by \cite{hamon_al:2020}, we evaluate the PinT errors (w.r.t. a given reference solution) in the spectral space. Indeed, it is known that PinT methods suffer from stability and convergence issues for the higher wavenumber modes of the solution\deleted[id=R2]{, which is especially critical in the case of hyperbolic or advection-dominated problems}. Therefore, we evaluate the ability of Parareal and MGRIT to converge on different regions of the wavenumber spectrum. Let $\psi_{m,n}$ be the spectral coefficient of a given function $\psi$ with modes $(m,n) \in [-M, \dots, M] \times [|m|, \dots, M]$ (corresponding to a triangular truncation of a spherical harmonics transform). We define a spectral resolution  $0 < R_{\text{norm}} \leq M$ for evaluating the errors. Then, the error of a numerical approximation to $\psi$, compared to a given reference solution $\psi_{\text{ref}}$, reads

\begin{equation}
    \label{eq:error_definition}
    E_{\psi,R_{\text{norm}}} := \frac{\norm{\infty, R_{\text{norm}}}{\psi - \psi_{\text{ref}}}}{\norm{\infty, R_{\text{norm}}}{\psi_{\text{ref}}}},\qquad \norm{\infty, R_{\text{norm}}}{\psi} := \max_{\substack{m \in \{ 0, \dots, R_{\text{norm}} \} \\ n \in \{m, \dots, R_{\text{norm}}\}}} |\psi_{m,n}|
\end{equation}

\indent In some situations, we also evaluate the more classical $L_2$ error computed in the physical space:

\begin{equation*}
    E_{\psi,L_2} := \frac{\norm{L_2}{\psi - \psi_{\text{ref}}}}{\norm{L_2}{\psi_{\text{ref}}}},\qquad \norm{L_2}{\psi} := \sqrt{ \frac{1}{M^{(x)}M^{(y)}}\sum_{i = 1}^{M^{(x)}}\sum_{j = 1}^{M^{(y)}} \psi_{i,j}^2} 
\end{equation*}

\noindent where $\psi_{i,j}$ are the discrete values of $\psi$ in a homogeneous physical grid with $M^{(x)}$ and $M^{(y)}$ points respectively along latitudes and longitudes.

\paragraph{Evaluation of computational time and speedup}

\indent The numerical simulations presented in the following paragraphs were executed in the GRICAD cluster from the University of Grenoble Alpes. The considered nodes are composed of 32 physical cores of Intel Xeon Gold 6130 @2.10GHz distributed between two sockets. The serial-in-time reference simulations were run considering a shared memory spatial parallelization using 16 OpenMP threads bound to cores in the same socket. The parallel-in-time simulations were run with a hybrid OpenMP-MPI parallelization, respectively in space and time, with $\nproc/2$ nodes being allocated to define $\nproc$ MPI tasks (which we name hereafter processors), each one using 16 OpenMP threads bound to cores on the same socket. \added[id=R3]{We remark that other space-time parallelization strategies could be adopted, \eg an MPI-MPI one relying on the splitting of the MPI communicator.}

\indent The computing times of the reference and the PinT simulations, the latter using $\nproc$ processors, are denoted $T_{\text{ref}}$ and $T_{\text{PinT}}(\nproc)$. The speedup provided by the temporal parallelization is defined by

\begin{equation*}
    S(\nproc) := \frac{T_{\text{ref}}}{T_{\text{PinT}}(\nproc)}
\end{equation*}

\subsection{Gaussian bumps test case}

\indent In this first test case, adapted from the single Gaussian bump benchmark presented by \cite{swarztrauber:2004}, a geopotential field with mean value $\ol{\Phi} = g\ol{h}$, with $\ol{h} = 29400\text{m}$ and $g = 9.80616 \text{ms}^{-2}$, is perturbed by three Gaussian bump centered at $(\lambda_1,\theta_1) = (\pi/5, \pi/3)$, $(\lambda_2,\theta_2) = (6\pi/5, \pi/5)$ and $(\lambda_3,\theta_3) = (8\pi/5, -\pi/4)$:

\begin{equation*}
    \Phi(\lambda, \theta) = \ol{\Phi} + A \sum_{i = 1}^3e^{-a_id_i(\lambda,\theta)^2}
\end{equation*}

\noindent where $A = 6000\text{m}$, $a_1 = 20$, $a_2 = 80$, $a_3 = 360$ and

\begin{equation*}
    d_i(\lambda,\theta) = \arccos(\sin(\theta)\sin(\theta_i) + \cos(\theta)\cos(\theta_i)\cos(\lambda - \lambda_i))
\end{equation*}

\indent The initial velocity field is identical to zero and the final time of simulation is $T = 36\text{h}$.

\begin{figure}[!htbp]
    \centering
    \includegraphics[scale=.5]{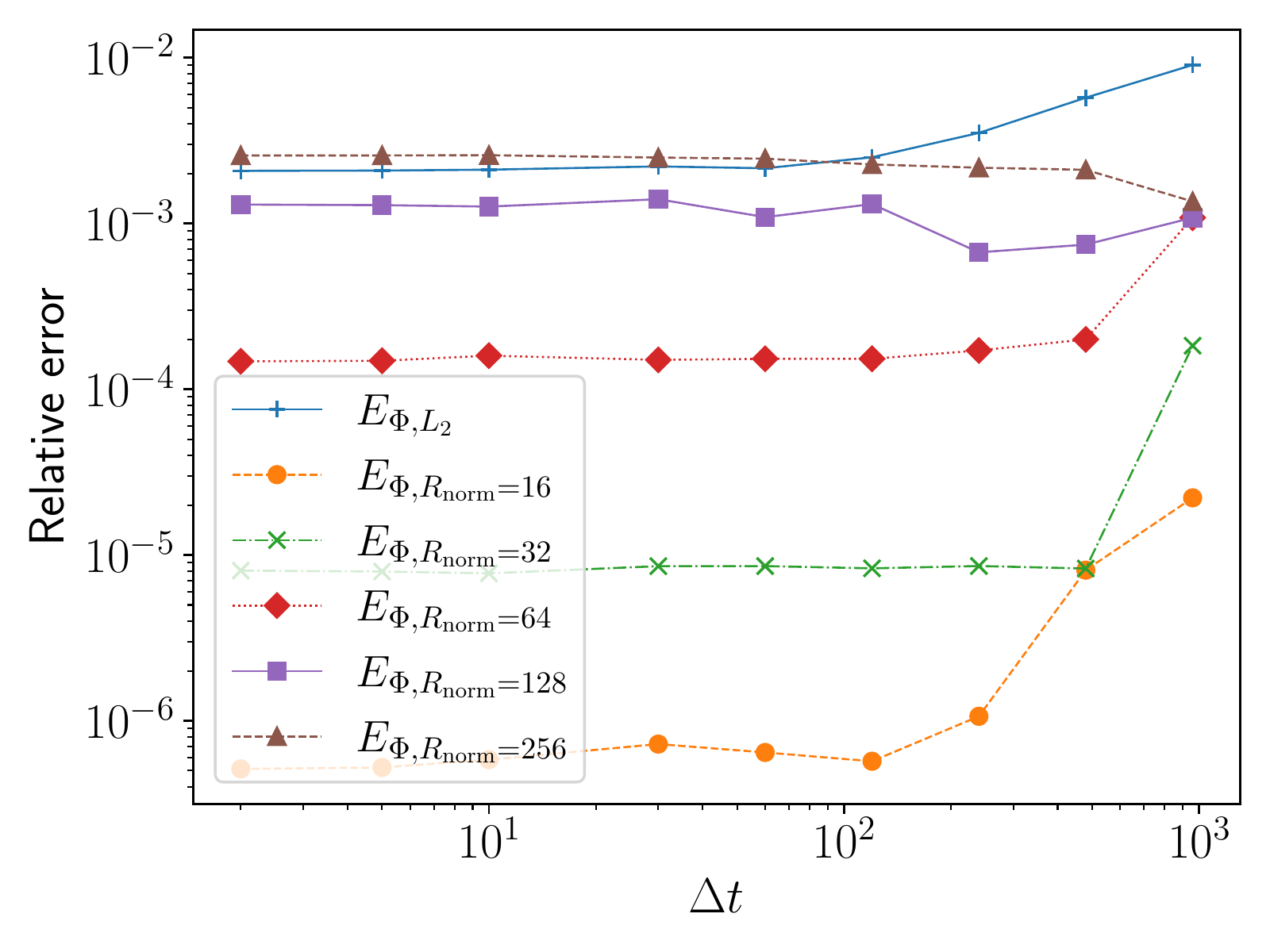}
    \caption{Gaussian bumps test case: relative $L_2$ error (computed in the physical space) and spectral errors (for various values of $\rnorm$) between a solution obtained with spectral resolution $M = 512$ and time step $\Dt = 2$ and solutions obtained with $M = M_0 = 256$ and various time steps. IMEX is used in all cases.}
    \label{fig:err_ref_gaussian_bumps}
\end{figure}

\begin{figure}[!htbp]
    \begin{subfigure}{.5\linewidth}
        \centering
        \includegraphics[scale=.3]{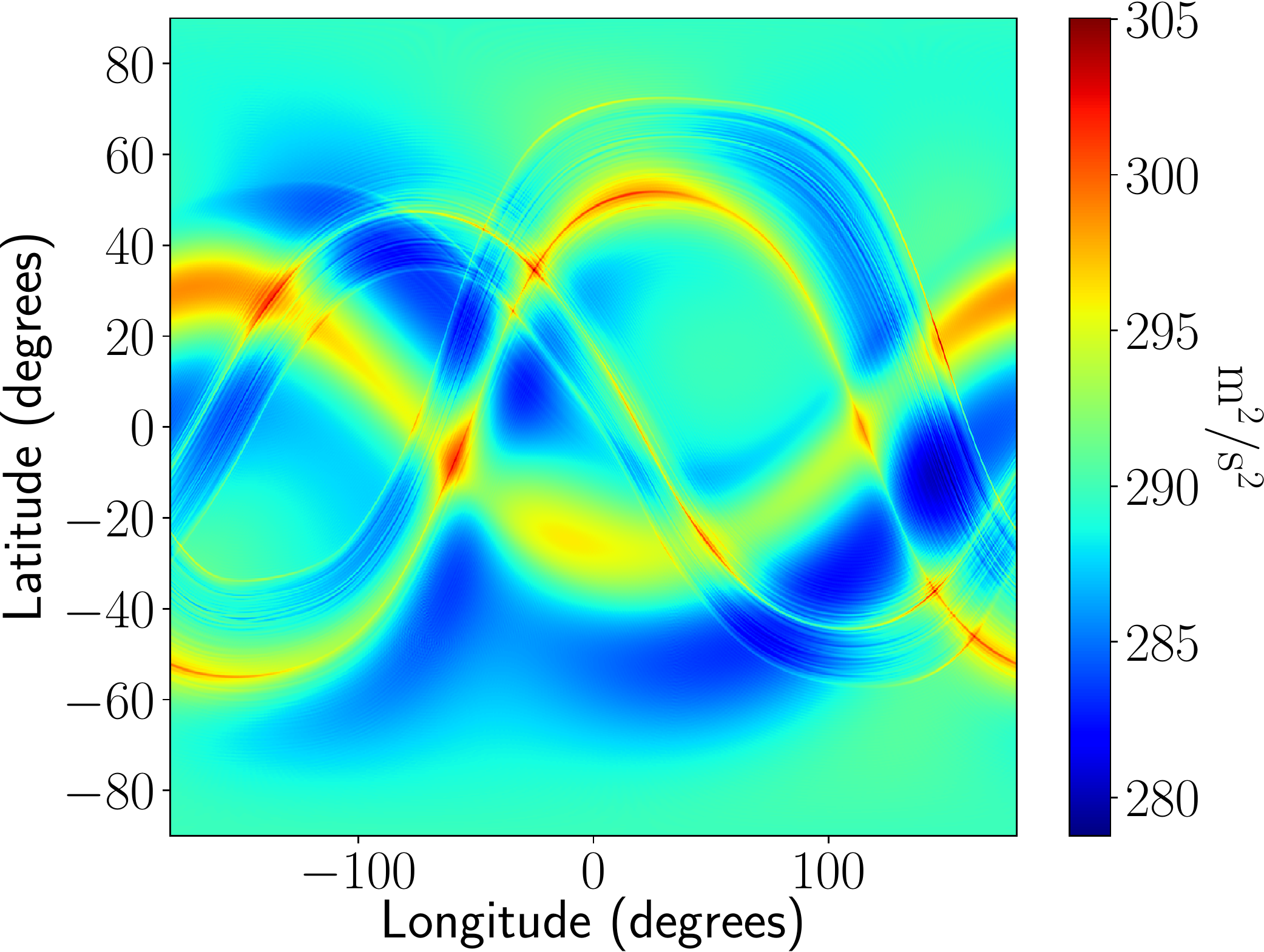}
        \caption{$M = 512$, $\Dt = 2$}
    \end{subfigure}
    \begin{subfigure}{.5\linewidth}
        \centering
        \includegraphics[scale=.3]{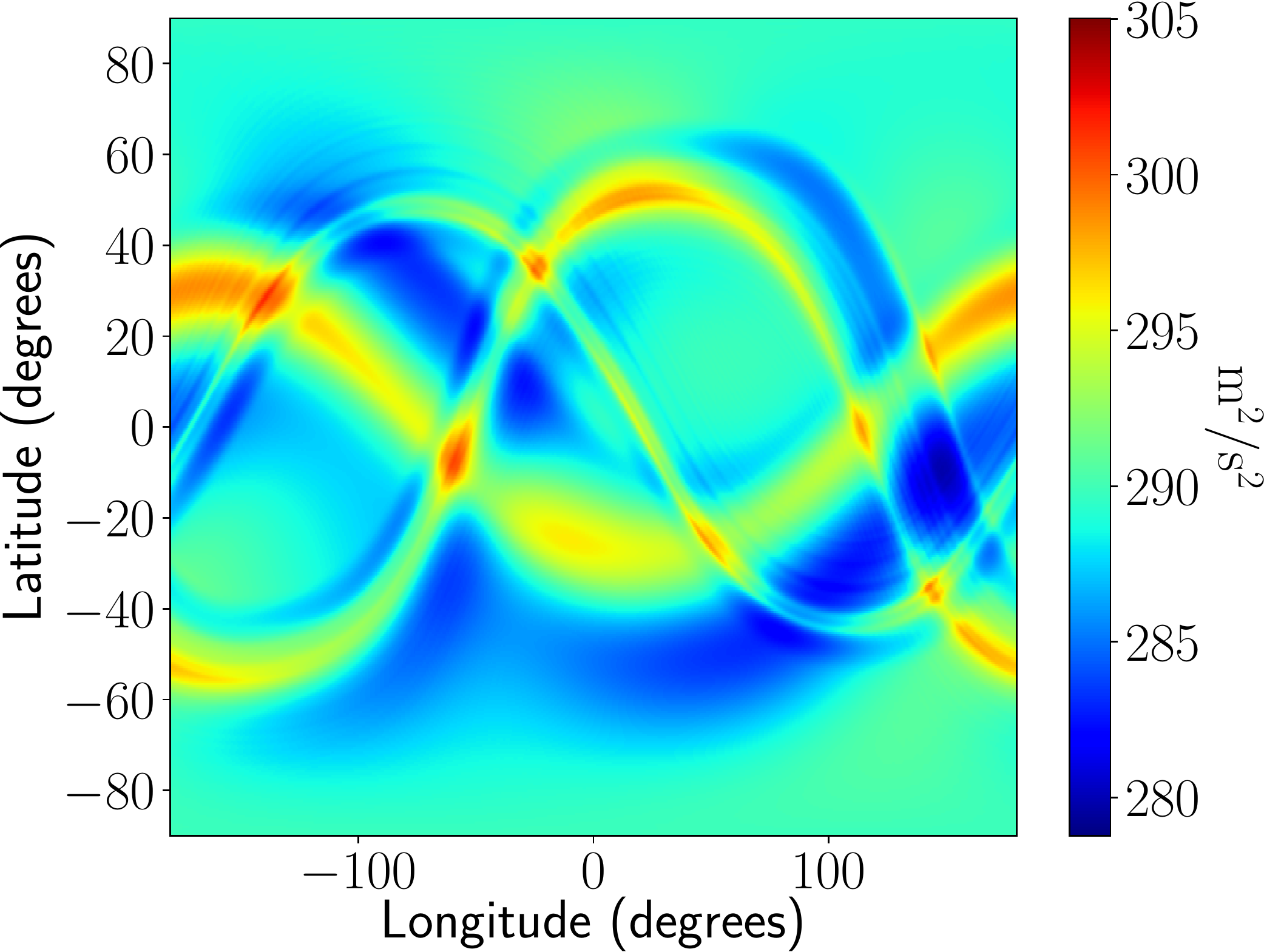}
        \caption{$M = 256$, $\Dt = 60$}
    \end{subfigure}
    \caption{Gaussian bumps test case: solution at $t = T = 36\text{h}$ computed with IMEX.}
    \label{fig:ref_gaussian_bumps}
\end{figure}

\indent We begin by defining a discretization to be used on the fine level of the parallel-in-time simulations (and also as a reference solution for evaluating the errors). Since the solution provided by Parareal and MGRIT converges to the one obtained via a serial simulation on the fine level, we want the latter to be a good enough approximation to the exact solution of the problem (or a very refined solution); on the other hand, we do not want the fine level to be overresolved (\eg with too fine discretizations in time and/or space), since it may lead to unrealistically overestimated speedups provided by the temporal parallelization \cite{gotschel_al:2021}. We fix the spectral discretization to $M_0 = 256$ and perform a serial simulation with time step sizes in $[2,960]$, which we compare to a solution obtained with spectral resolution $M = 512$ and time step $\Dt = 2$. In all cases, no artificial viscosity is used, and the solutions are obtained with IMEX. Figure \ref{fig:err_ref_gaussian_bumps} presents the $L_2$ errors, computed in the physical space, as well as the spectral errors for given $\rnorm$ values. Considering all the presented norms, no important error reduction is observed by using time steps smaller than approximately $100\text{s}$, with the errors being dominated by spatial discretization. We then choose $(M_0 = 256, \Dt_0 = 60)$ to be used on the fine levels of the temporal parallelization. The solution at $t = T$ obtained under this configuration is compared to the solution obtained with $(M = 512, \Dt = 2)$ in Figure \ref{fig:ref_gaussian_bumps}. We notice that for the chosen simulation length $T$, the turbulence regime is not yet fully developed and the $-5/3$ power law of the two-dimensional kinetic energy spectrum \wrt the wavenumber \cite{lindborg:1999} is not yet verified, as illustrated in Figure \ref{fig:gaussian_bumps_spectrum}, which compares the kinetic energy spectrum of the fine solution with those obtained with selected PinT configurations.

\subsubsection{Convergence study}

\indent  The finest level in the PinT configuration is discretized with the same time step size $\Dt_0 = \Dt_{\text{ref}} = 60$, spectral resolution $M_0 = M_{\text{ref}} = 256$ and artificial viscosity coefficient $\nu_0 = \nu_{\text{ref}} = 0$ used in the reference simulation. We perform MGRIT simulations considering $\nlevels \in \{2,3\}$, with a coarsening factor $\cfactor \in \{2, 4\}$. The same spectral resolution is considered on all coarse levels, being set to $M_{\text{coarse}} \in \{51, 128\}$, which correspond respectively to spatial coarsening factors of $1/5$ and $1/2$ \wrt the finest level. \added[id=R3]{The use of spatial coarsening on the coarse levels is motivated by stability issues observed in the PinT simulations, both when the coarse schemes are IMEX and SL-SI-SETTLS, if too large spectral resolutions are adopted since relatively large time step sizes are used on the coarse levels.} We also consider several relaxation strategies, with $\nrelax \in \{0,1,5\}$. We recall that MGRIT with $(\nlevels, \nrelax) = (2,0)$ is equivalent to Parareal. Finally, all coarse levels $l>0$ use an artificial second-order viscosity coefficient $\nu_{\text{coarse}} =  10^6 \visc{2}$, which is the same value adopted in the numerical simulations presented by \cite{schreiber_al:2019} for integrating the SWE on the rotating sphere with PFASST. A more detailed study on the influence of the artificial viscosity order and coefficient applied in each discretization level is presented in Section \ref{subsec:gaussian_bump_viscosity}.

\begin{figure}[!htbp]
    \begin{subfigure}{.5\linewidth}
        \centering
        \includegraphics[scale=.425]{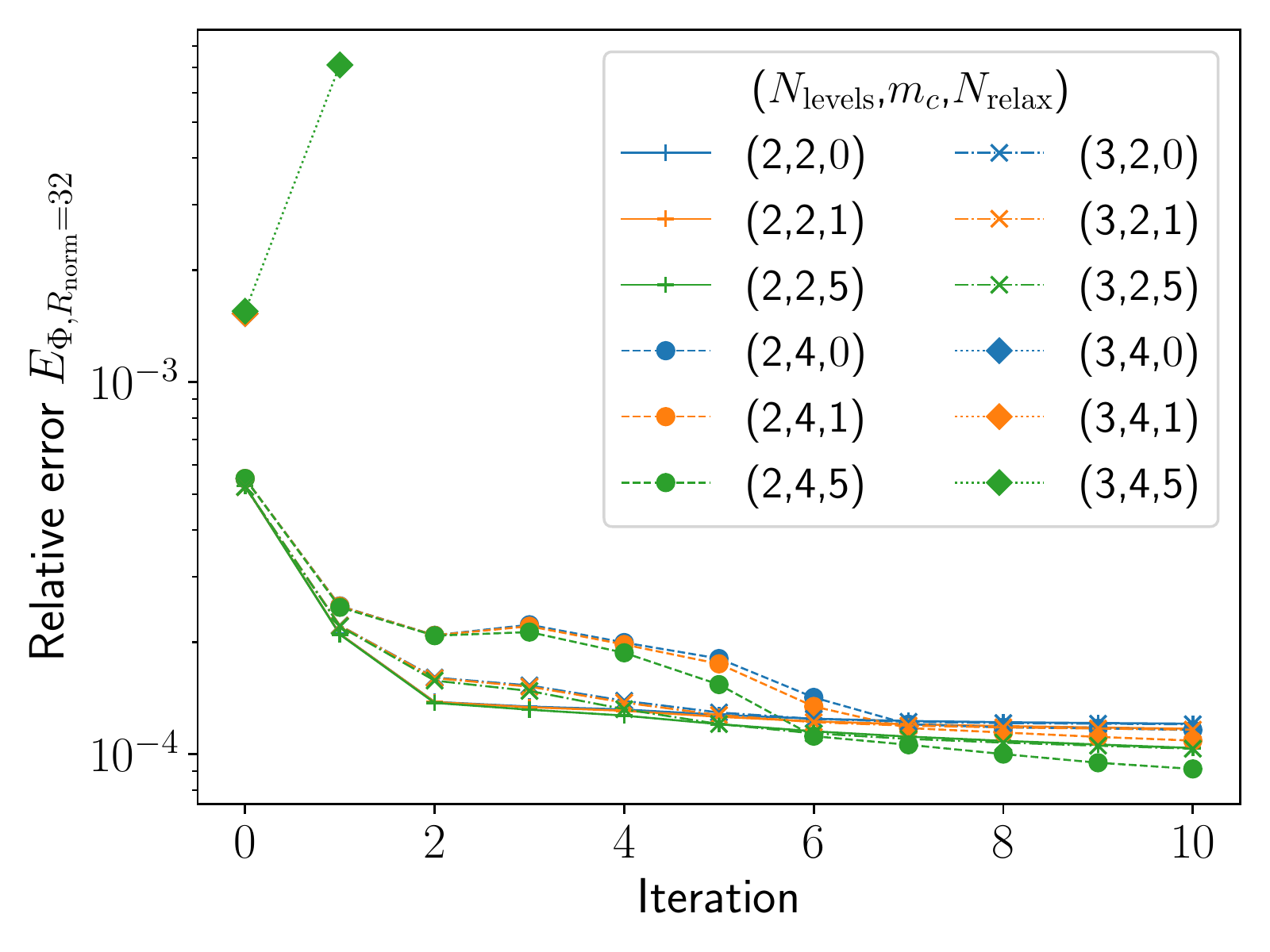}
        \caption{$M_{\text{coarse}} = 51$, $\rnorm = 32$}
    \end{subfigure}
    \begin{subfigure}{.5\linewidth}
        \centering
        \includegraphics[scale=.425]{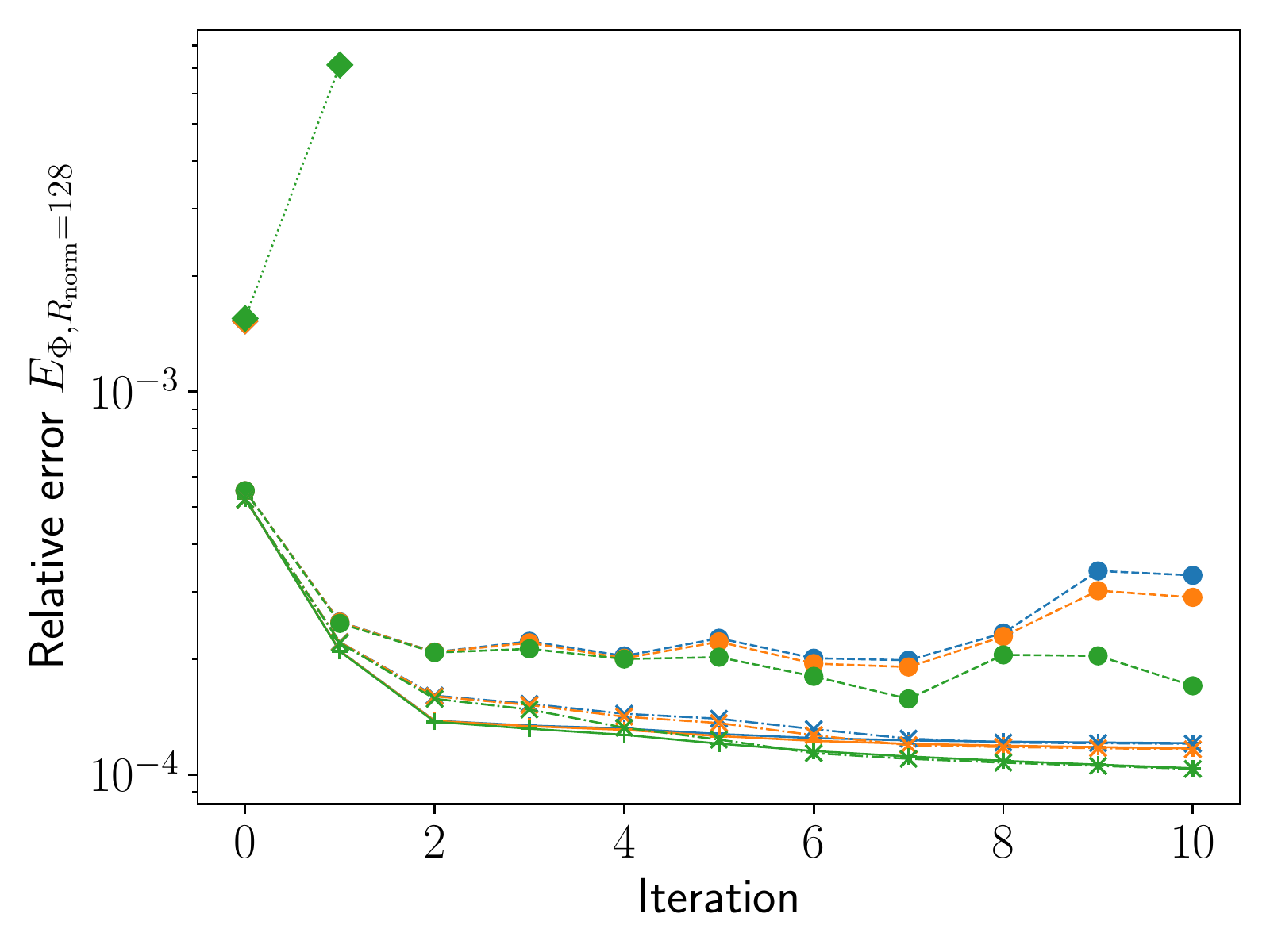}
        \caption{$M_{\text{coarse}} = 51$, $\rnorm = 128$}
    \end{subfigure}
    \begin{subfigure}{.5\linewidth}
        \centering
        \includegraphics[scale=.425]{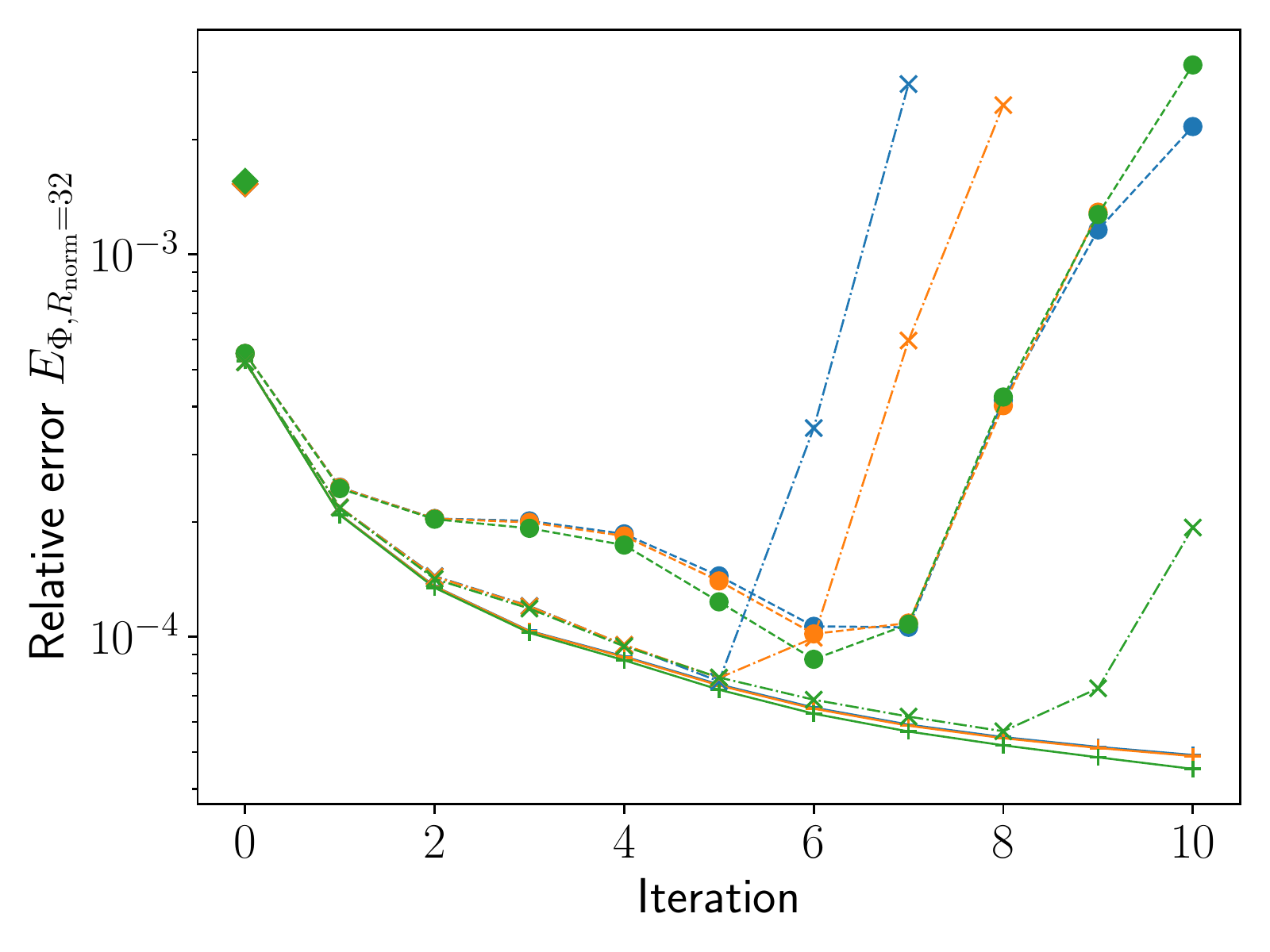}
        \caption{$M_{\text{coarse}} = 128$, $\rnorm = 32$}
    \end{subfigure}
    \begin{subfigure}{.5\linewidth}
        \centering
        \includegraphics[scale=.425]{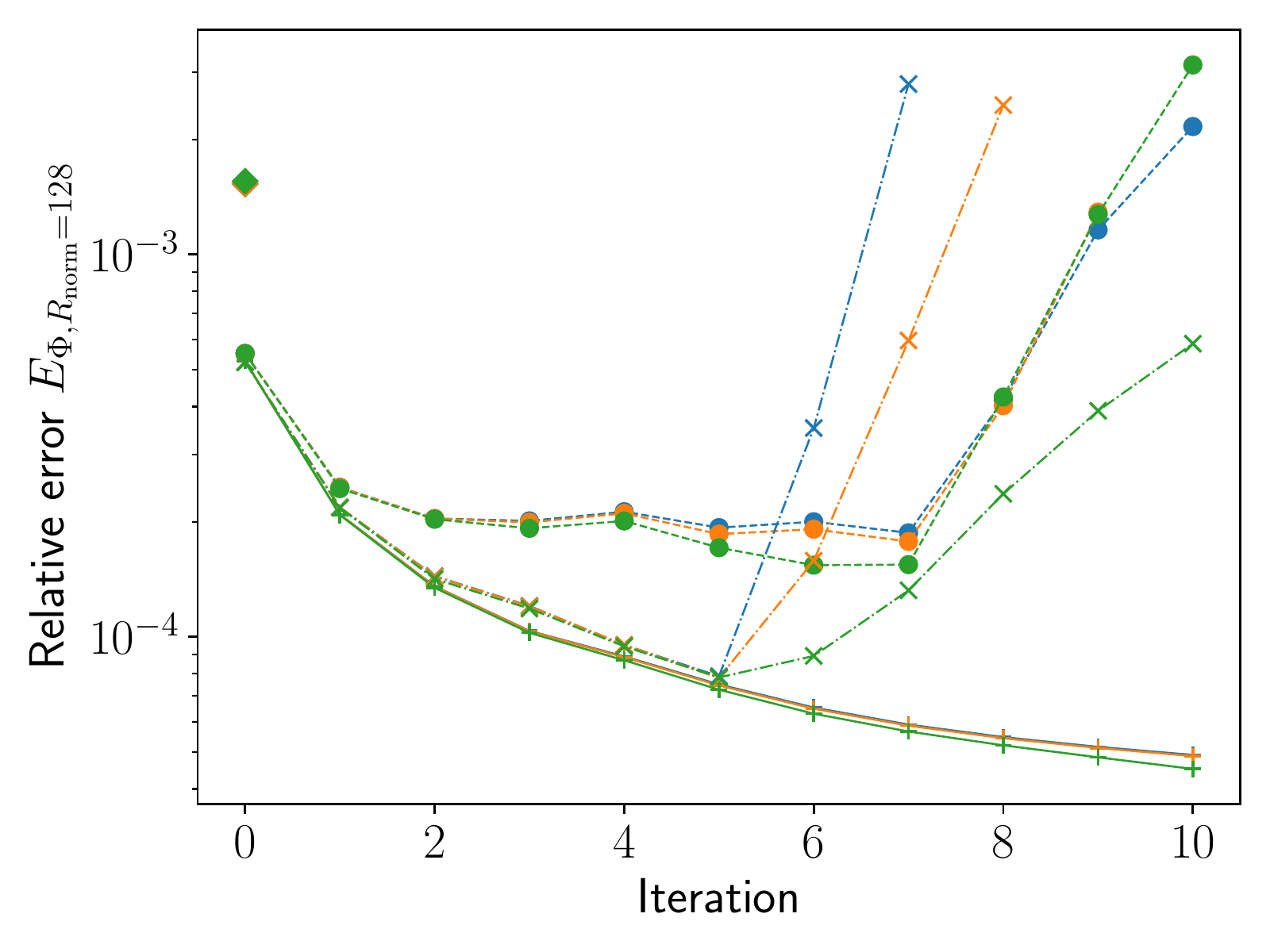}
        \caption{$M_{\text{coarse}} = 128$, $\rnorm = 128$}
    \end{subfigure}
    \caption{Gaussian bumps test case: relative error $E_{\Phi,\rnorm}$ between the PinT and fine solutions at $t = T$ along iterations, for $\rnorm = 32$ (left) and $\rnorm = 128$ (right); and $\Mcoarse = 51$ (top) and $\Mcoarse = 128$ (bottom), with IMEX used on the coarse levels. Simulations are identified by $(\nlevels, \cfactor, \nrelax)$.}
    \label{fig:gaussian_bumps_errors_params_mgrit_IMEX}
\end{figure}

\indent Figure \ref{fig:gaussian_bumps_errors_params_mgrit_IMEX} presents the evolution, along iterations, of the relative spectral error between the parallel-in-time and fine solutions, evaluated at the final time of simulation when IMEX is used on the coarse levels. The errors are shown for two chosen spectral resolutions, namely $\rnorm = 32$ and $\rnorm = 128$. We observe that only a few simulations are stable and able to converge (and, when this is the case, with a relatively fast error decrease in the first iterations followed by a slower convergence in the next ones), namely when the spectral resolution on the coarse levels, the number of levels and/or the coarsening factor are sufficiently small; otherwise, the time step on the coarsest level is too large and leads to unstable behaviors. With $\Mcoarse = 51$, all configurations can converge, more or less rapidly, at the beginning of the wavenumber spectrum ($\rnorm = 32$), except for the most aggressive configuration ($(\nlevels, \cfactor) = (3,4)$); when more modes are considered in the error analysis ($\rnorm = 128$), the simulations with $(\nlevels, \cfactor) = (2,4)$ also starts to develop instabilities that compromise convergence, with this unstable behavior being less remarkable when a more expensive relaxation strategy is adopted. On the other hand, when a finer spectral resolution is considered on the coarse levels ($\Mcoarse = 128$), almost none of the MGRIT configurations are stable, even under $\rnorm = 32$, the only exception being the less aggressive one ($(\nlevels, \cfactor) = (2,2)$). The configuration $(\nlevels, \cfactor) = (3,2)$ initially presents a convergent behavior; still, the error starts to increase after five iterations, which takes place sooner at the end of the wavenumber spectrum ($\rnorm = 128$), and when the relaxation is less expensive; under other configurations, the instabilities start to develop even sooner. Finally, in all convergent configurations, we notice a slight improvement when the relaxation strategy is more expensive, but the error behaviors are essentially the same.

\begin{figure}[!htbp]
    \begin{subfigure}{.5\linewidth}
        \centering
        \includegraphics[scale=.425]{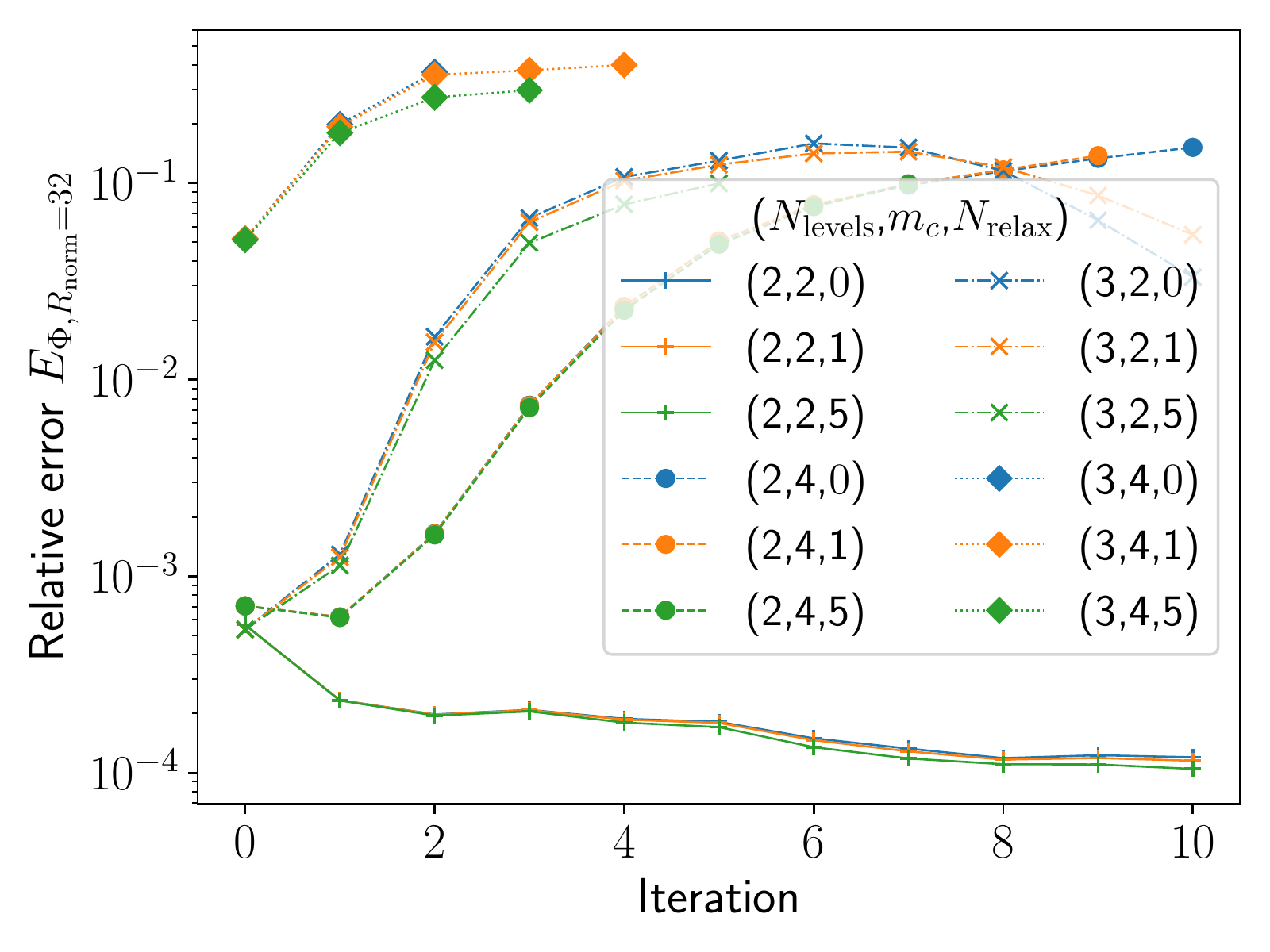}
        \caption{$M_{\text{coarse}} = 51$, $\rnorm = 32$}
    \end{subfigure}
    \begin{subfigure}{.5\linewidth}
        \centering
        \includegraphics[scale=.425]{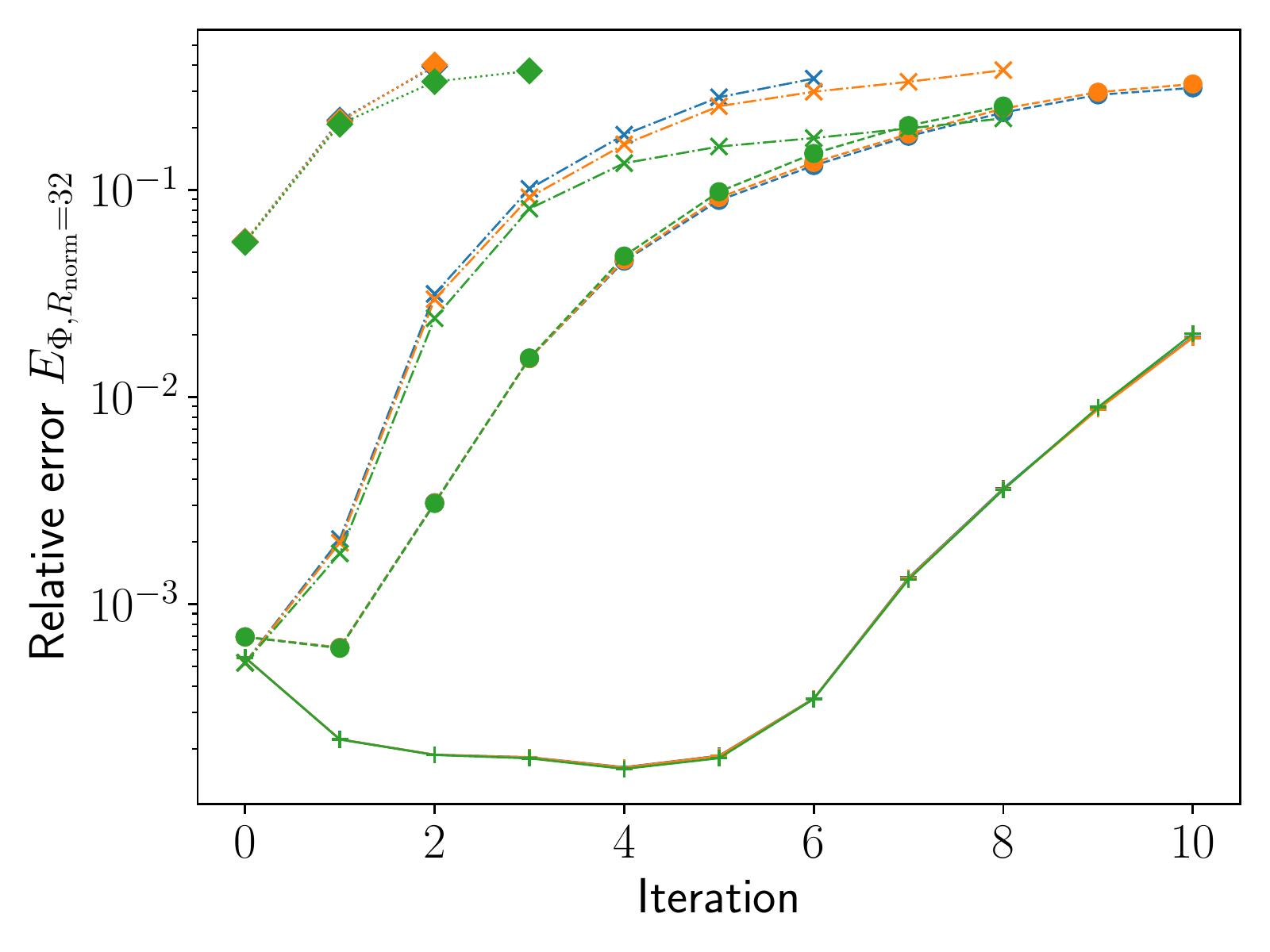}
        \caption{$M_{\text{coarse}} = 128$, $\rnorm = 32$}
    \end{subfigure}
    \caption{Gaussian bumps test case: relative error $E_{\Phi,\rnorm}$ between the PinT and fine solutions at $t = T$ along iterations, for $\rnorm = 32$ and $\Mcoarse = 51$ (left) and $\Mcoarse = 128$ (right), with SL-SI-SETTLS used on the coarse levels. Results are identical under $\rnorm = 128$. Simulations are identified by $(\nlevels, \cfactor, \nrelax)$.}
    \label{fig:gaussian_bumps_errors_params_mgrit_SL_SI_SETTLS}
\end{figure}

\indent The same results are presented in Figure \ref{fig:gaussian_bumps_errors_params_mgrit_SL_SI_SETTLS} for the simulations using SL-SI-SETTLS on the coarse levels. We notice a more important unstable behavior compared to the simulations using IMEX as a coarse scheme. Among all tested configurations, the only one presenting a relatively stable and convergent behavior is $(\nlevels,\cfactor,\Mcoarse) = (2,2,51)$, \ie the one with the smallest time step and spectral resolution on the coarsest level (however, even in this case the convergence is less clear and not monotonic as observed in simulations using IMEX). In all other configurations, the error increases after a few iterations, independently of the chosen norm for evaluating the error. Moreover, contrary to the results with IMEX, we do not observe significant differences between the errors computed with $\rnorm = 32$ and $\rnorm = 128$, meaning that simulations fail to converge in the large spatial scales, and only the former case is presented. Finally, we observe only very small improvements when using more expensive relaxation strategies, but it is not enough to improve the stability of the simulations and lead to convergent behavior.

\begin{figure}[!htbp]
    \centering
    \begin{minipage}{.3\linewidth}
        \subcaptionbox{Iteration 0}
        {
              \includegraphics[scale=.225]{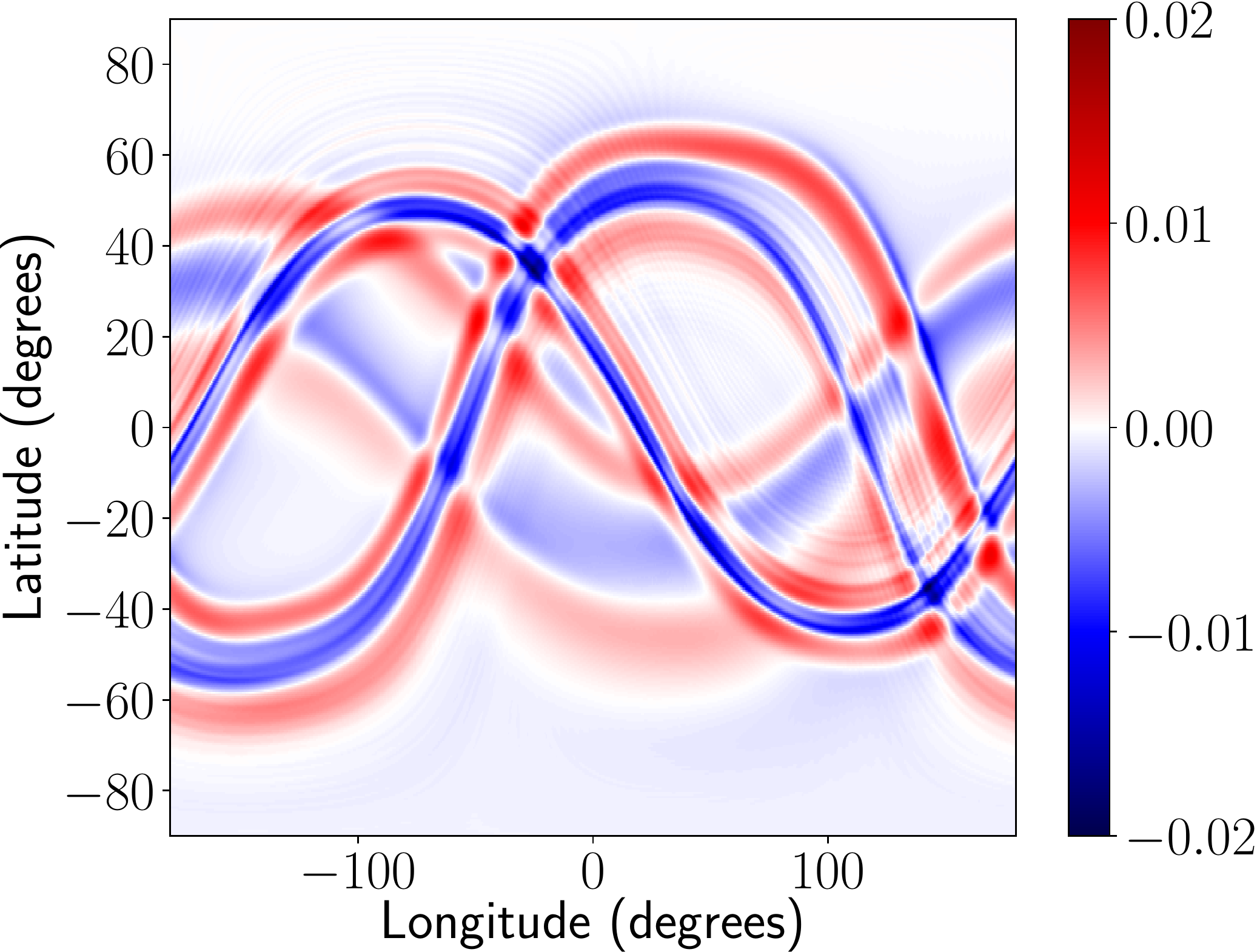} 
        }
    \end{minipage}
    \begin{minipage}[t]{.3\linewidth}
        \subcaptionbox{$(3,2,0,51,\text{IMEX})$}
        {
              \includegraphics[scale=.225]{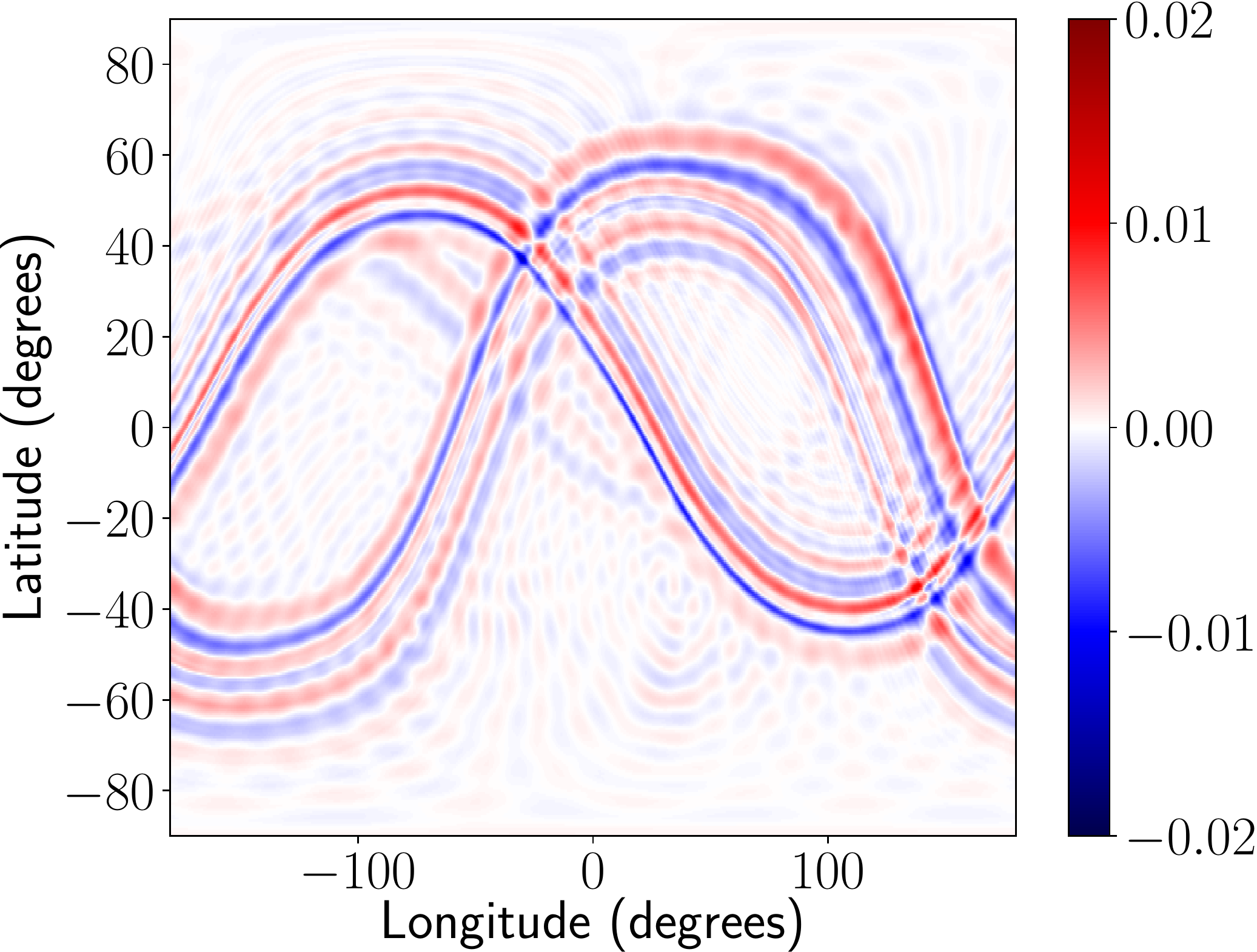}
        }
        \subcaptionbox{$(2,2,0,51,\text{SL-SI-SETTLS})$}
        {
            \includegraphics[scale=.225]{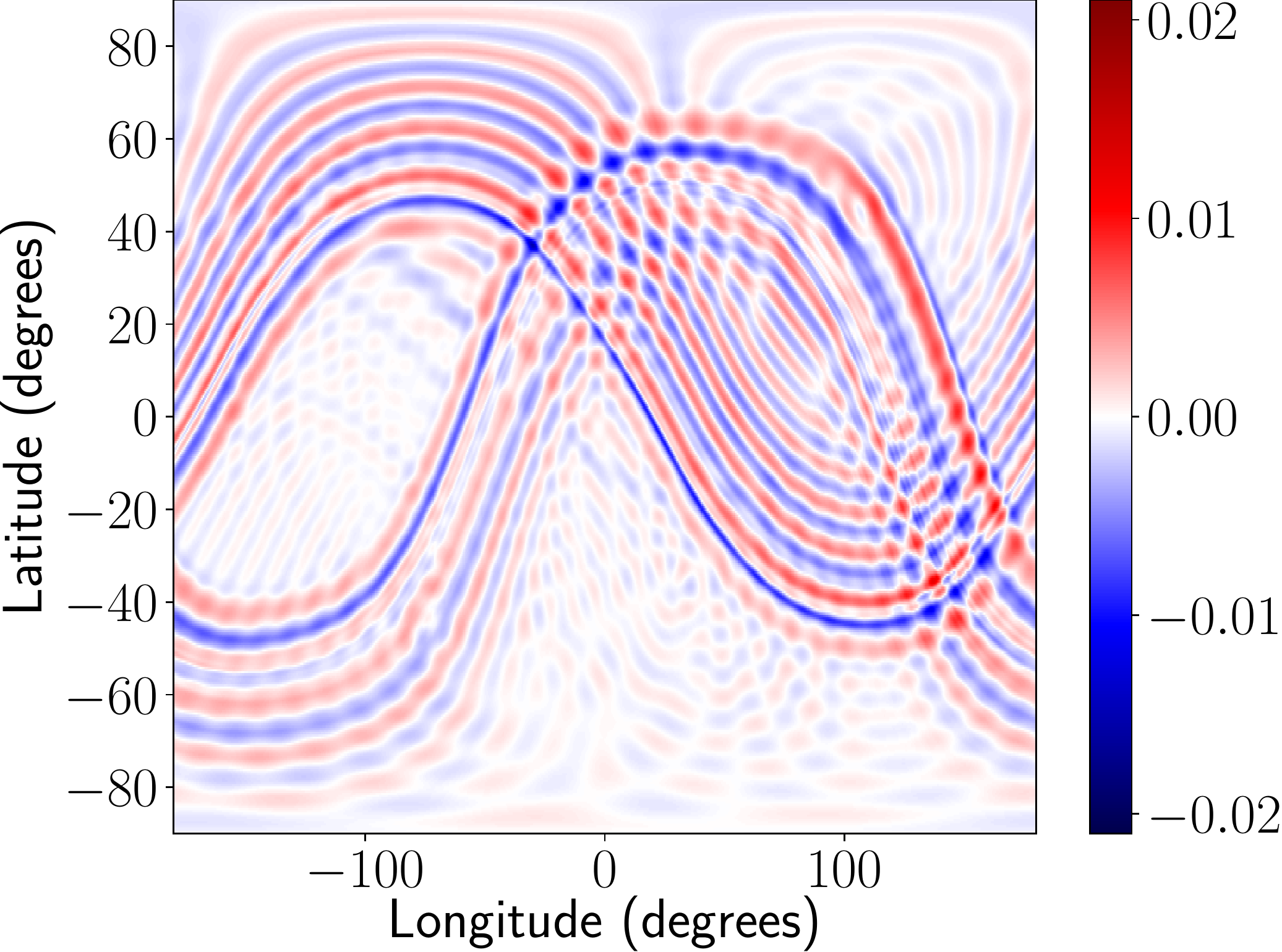}
        }
    \end{minipage}
    \begin{minipage}[t]{.3\linewidth}
        \subcaptionbox{$(3,2,5,128,\text{IMEX})$}
        {
            \includegraphics[scale=.225]{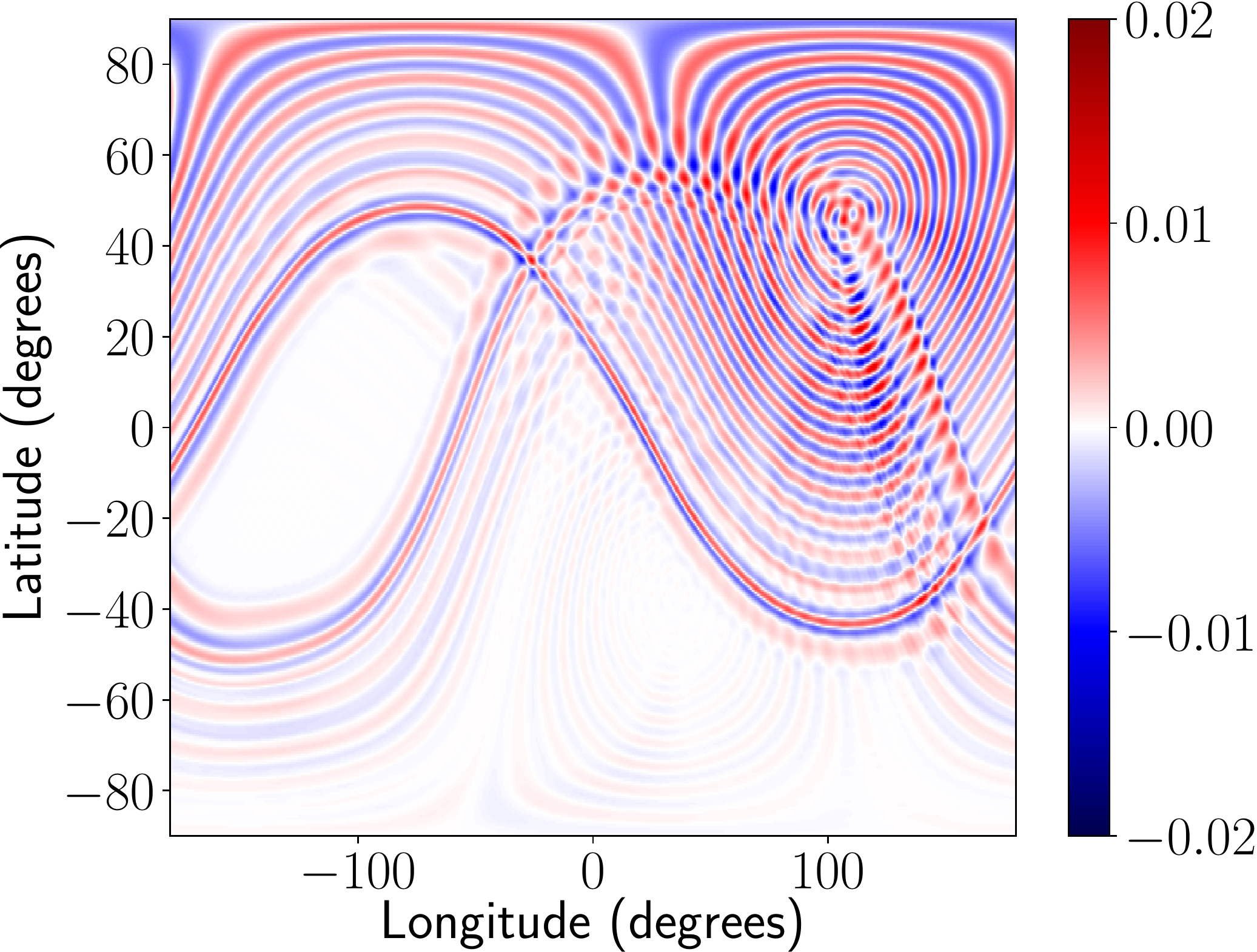}
        }
        \subcaptionbox{$(2,2,0,128,\text{SL-SI-SETTLS})$}
        {
            \includegraphics[scale=.225]{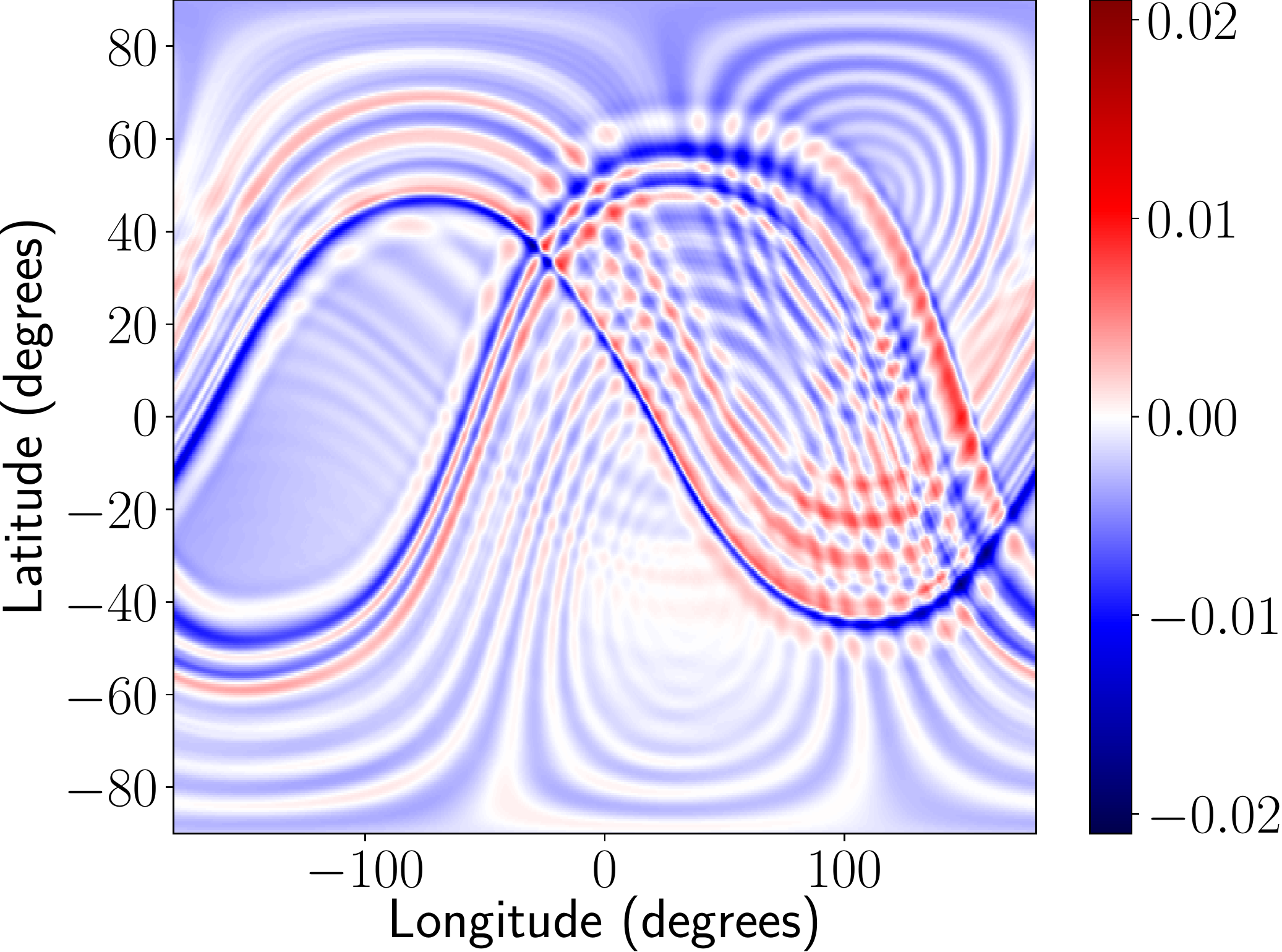}
        }
    \end{minipage}
    \caption{Gaussian bumps test case: the relative difference between the PinT and fine geopotential fields at $t = T$ after 7 iterations for chosen configurations $(\nlevels,\cfactor,\nrelax,\Mcoarse, \text{coarse time-stepping method})$. The left plot represents the solution at the initial iteration, which is visually similar for all configurations.}
    \label{fig:gaussian_bumps_solution_params_mgrit}
\end{figure}

\indent To illustrate the physical solutions obtained in the PinT simulations, we present in Figure \ref{fig:gaussian_bumps_solution_params_mgrit} the relative difference between the PinT and the reference geopotential fields at iteration $k = 7$ for chosen configurations, compared to the solution obtained at the initial iteration. With IMEX, we choose $(\nlevels, \cfactor, \nrelax, \Mcoarse) = (3,2,0,51)$ and $(\nlevels, \cfactor, \nrelax, \Mcoarse) = (3,2,5,128)$. The former case presents a monotonic convergence behavior, and we observe small-scale oscillations of the error. On the other hand, the latter case diverges after a few iterations (which takes place sooner in the largest wavenumbers, \ie for larger $\rnorm$ values): in the seventh iteration, the PinT already diverged under $\rnorm = 128$ but not under $\rnorm = 32$, and we indeed observe fine oscillations dominating the error plot. With SL-SI-SETTLS as coarse scheme, we present the results for $(\nlevels, \cfactor, \nrelax, \Mcoarse) = (2,2,0,51)$ and $(\nlevels, \cfactor, \nrelax, \Mcoarse) = (2,2,0,128)$. The former presents the best convergence behavior among the simulations depicted in Figure \ref{fig:gaussian_bumps_errors_params_mgrit_SL_SI_SETTLS}, but we observe an increasing magnitude of the fine-scale errors along iterations. The latter configuration starts to diverge after five iterations, under $\rnorm = 32$ and $\rnorm = 128$: indeed, the plot reveals small- and fine-scale errors at iteration 7.

\begin{figure}[!htbp]
    \begin{subfigure}{.5\linewidth}
        \centering
        \includegraphics[scale=.35]{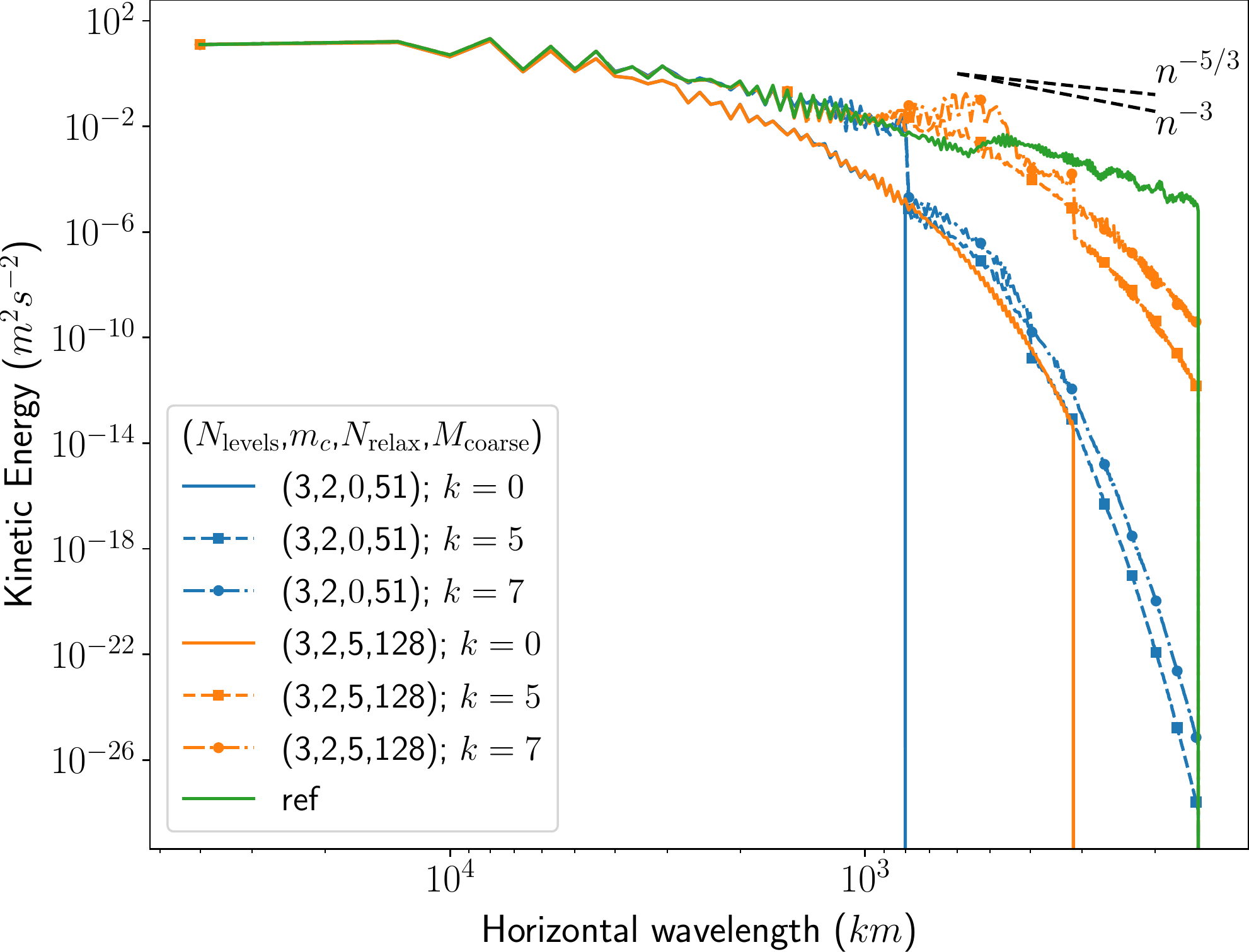}
        \caption{IMEX}
    \end{subfigure}
    \begin{subfigure}{.5\linewidth}
        \centering
        \includegraphics[scale=.35]{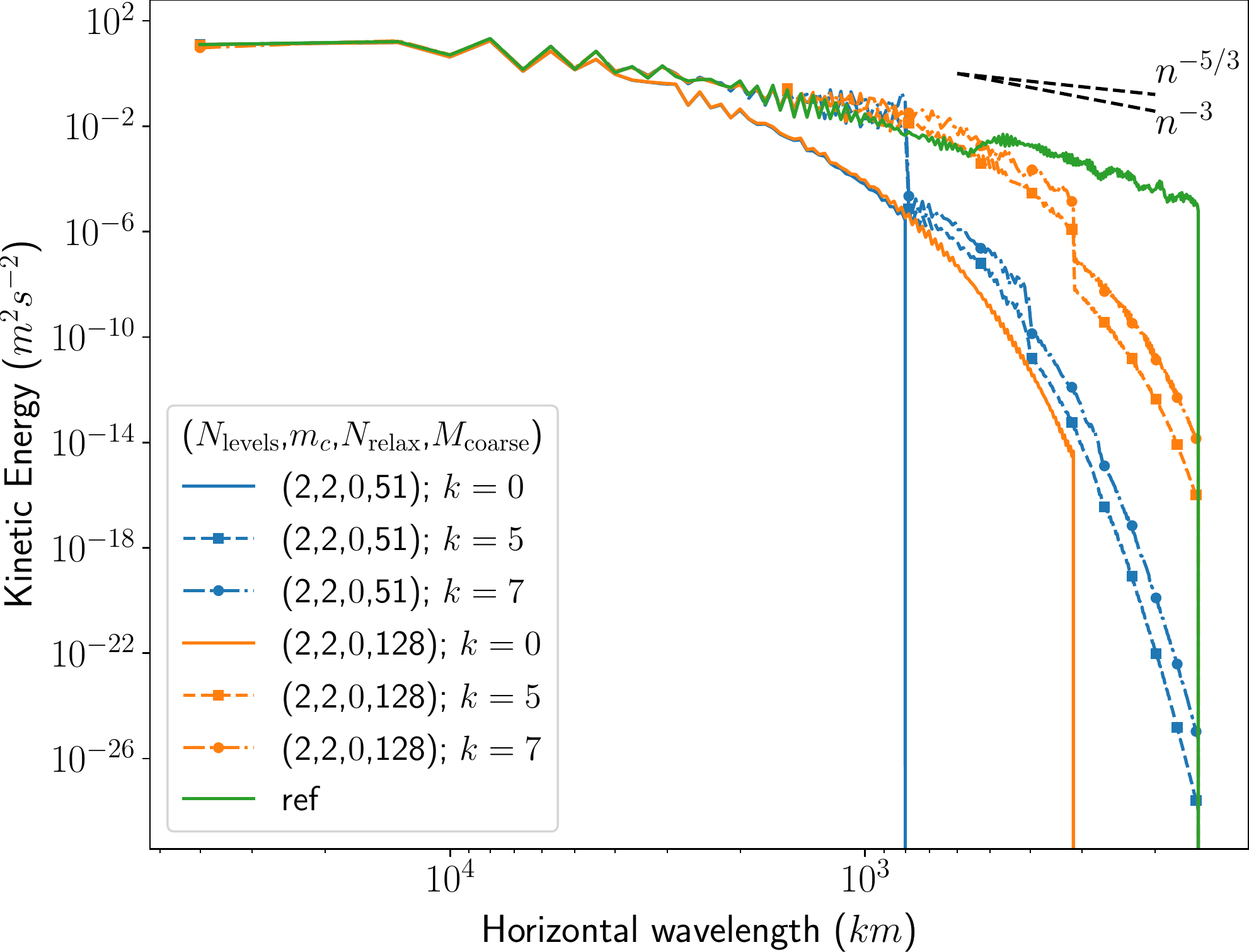}
        \caption{SL-SI-SETTLS}
    \end{subfigure}
    \caption{Gaussians bumps test case: kinetic energy spectra at $t = T$ of the reference (fine) and PinT solutions at iterations 0, 5 and 7 for given configurations. Left and right: IMEX and SL-SI-SETTLS as coarse time-stepping schemes, respectively.} 
    \label{fig:gaussian_bumps_spectrum}
\end{figure}

\indent The stability and convergence behaviors of the simulations presented in Figure \ref{fig:gaussian_bumps_solution_params_mgrit} are also illustrated in Figure \ref{fig:gaussian_bumps_spectrum}, in which we compare the kinetic energy (KE) spectrum at $t = T$ at given iterations with those corresponding to the fine, reference solution. For a triangular spectral truncation, the KE spectrum is defined by \cite{koshyk_hamilton:2001}

\begin{equation*}
    E(n, t) = \frac{1}{4} \frac{a^2}{n(n+1)} \sum_{m = -n}^n\left( |\xi_{m,n}(t)|^2 + |\delta_{m,n}(t)|^2 \right)
\end{equation*}

\noindent where $\xi_{m,n}(t)$ and $\delta_{m,n}(t)$ are the spherical harmonics coefficients of the vorticity and divergence, respectively. In Figure \ref{fig:gaussian_bumps_spectrum}, we present the spectrum as a function of the wavelength $\mc{L}_n := 2\pi a / n$. Concerning the simulations using IMEX, the unstable behavior under configuration $(\nlevels,\cfactor,\nrelax,\Mcoarse) = (3,2,5,128)$ is clear in the spectrum, with over-amplifications of medium to large wavenumbers along iterations; on the other hand, the spectrum of the simulation using $(\nlevels,\cfactor,\nrelax,\Mcoarse) = (3,2,0,51)$ remains below the fine one, but slowly converging to it, indicating a more difficult convergence of large wavenumber modes. In the simulations using SL-SI-SETTLS, we observe some over-amplification in the middle region of the wavenumber spectrum in both PinT configurations.

\subsubsection{Influence of artificial (hyper)viscosity}
\label{subsec:gaussian_bump_viscosity}

\indent We now study the influence of the artificial (hyper)viscosity coefficients applied on each coarse discretization level on the convergence and stability performance of the PinT methods. Since different temporal and spatial discretization sizes are used on each level, their respective stability constraints are not the same; therefore, different viscosity approaches should probably be applied on different levels. A brief overview of viscosity and hyperviscosity approaches in spectral methods and some guidelines for choosing the magnitude of the viscosity coefficients are presented in \ref{app:viscosity}. For the study proposed here, we consider the configuration $(\nlevels, \cfactor, \nrelax, \Mcoarse) = (3,2,0,128)$, both with IMEX and SL-SI-SETTLS as coarse time-stepping schemes. In the former case, this configuration leads to an initial convergent behavior followed by instabilities arising after the fifth iteration (Figure \ref{fig:gaussian_bumps_errors_params_mgrit_IMEX}); in the latter, the unstable behavior is observed from the initial iteration (Figure \ref{fig:gaussian_bumps_errors_params_mgrit_SL_SI_SETTLS}). As before, no viscosity is used on the fine level $(\nu_0 = 0)$; for the two coarse levels, we consider viscosity orders $q_1$ and $q_2$ from $\{2,4,6\}$, with the viscosity coefficients respectively in $\{10^5,10^6,10^7\}$ ($\text{m}^2\text{s}^{-1}$), $\{10^{15},10^{16},10^{17}\}$ ($\text{m}^4\text{s}^{-1}$) and $\{10^{25},10^{26},10^{27}\}$ ($\text{m}^6\text{s}^{-1}$). \added[id=R3]{These ranges of coefficients are coherent with those documented in the literature of atmospheric modeling and used in operational models \cite{jablonowski_williamson:2011}; for instance, reported values of fourth-order coefficient varying between orders of $10^{12}\visc{4}$ and $10^{16} \visc{4}$, with larger coefficients being adopted, in general, when the spectral resolution is smaller.}

\begin{figure}[!htbp]
    \begin{subfigure}{.5\linewidth}
        \centering
        \includegraphics[scale=.35]{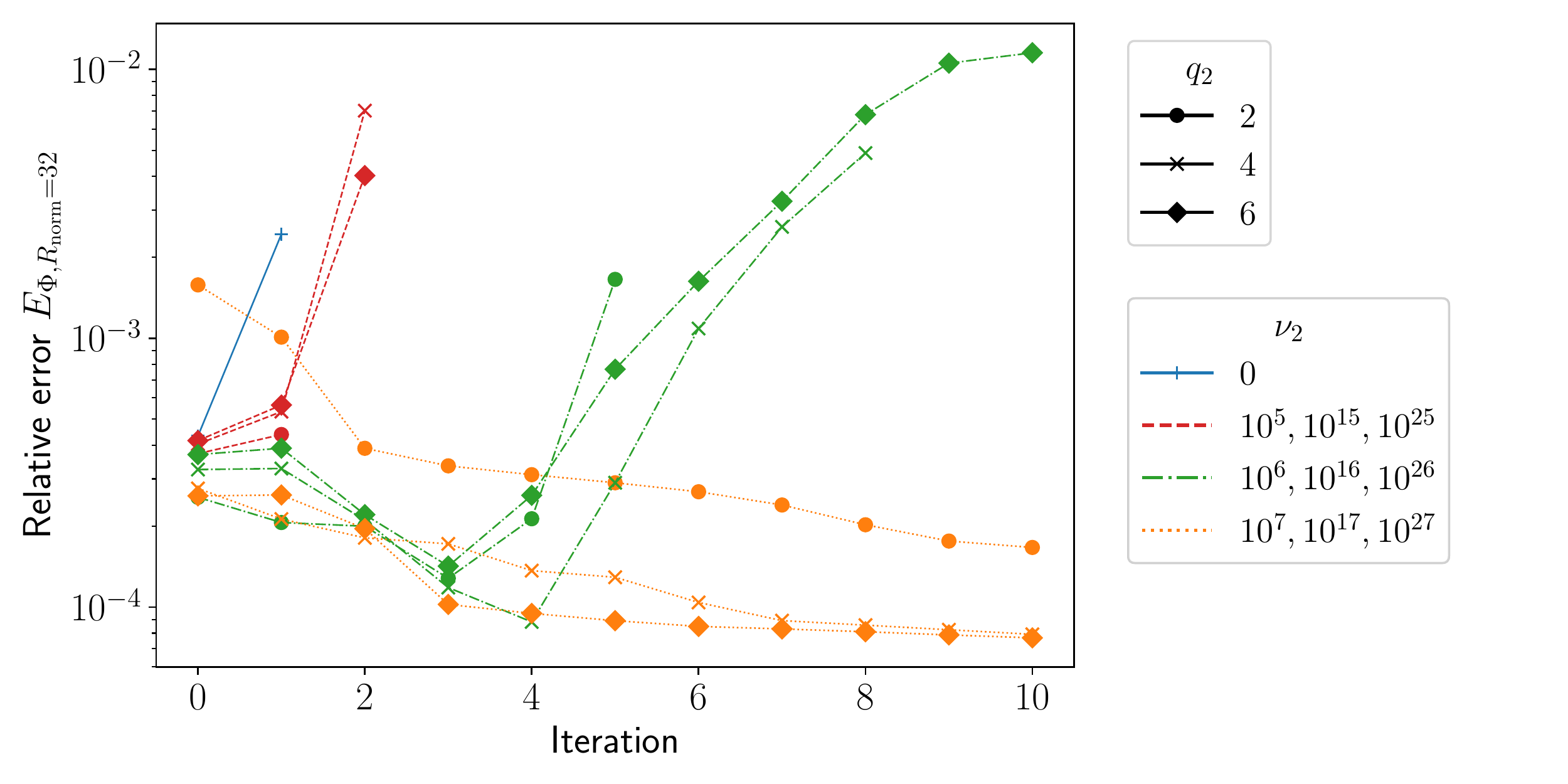}
        \caption{IMEX}
        \label{fig:gaussian_bumps_errors_viscosity_IMEX}
    \end{subfigure}
    \begin{subfigure}{.33\linewidth}
        \centering
        \includegraphics[scale=.35]{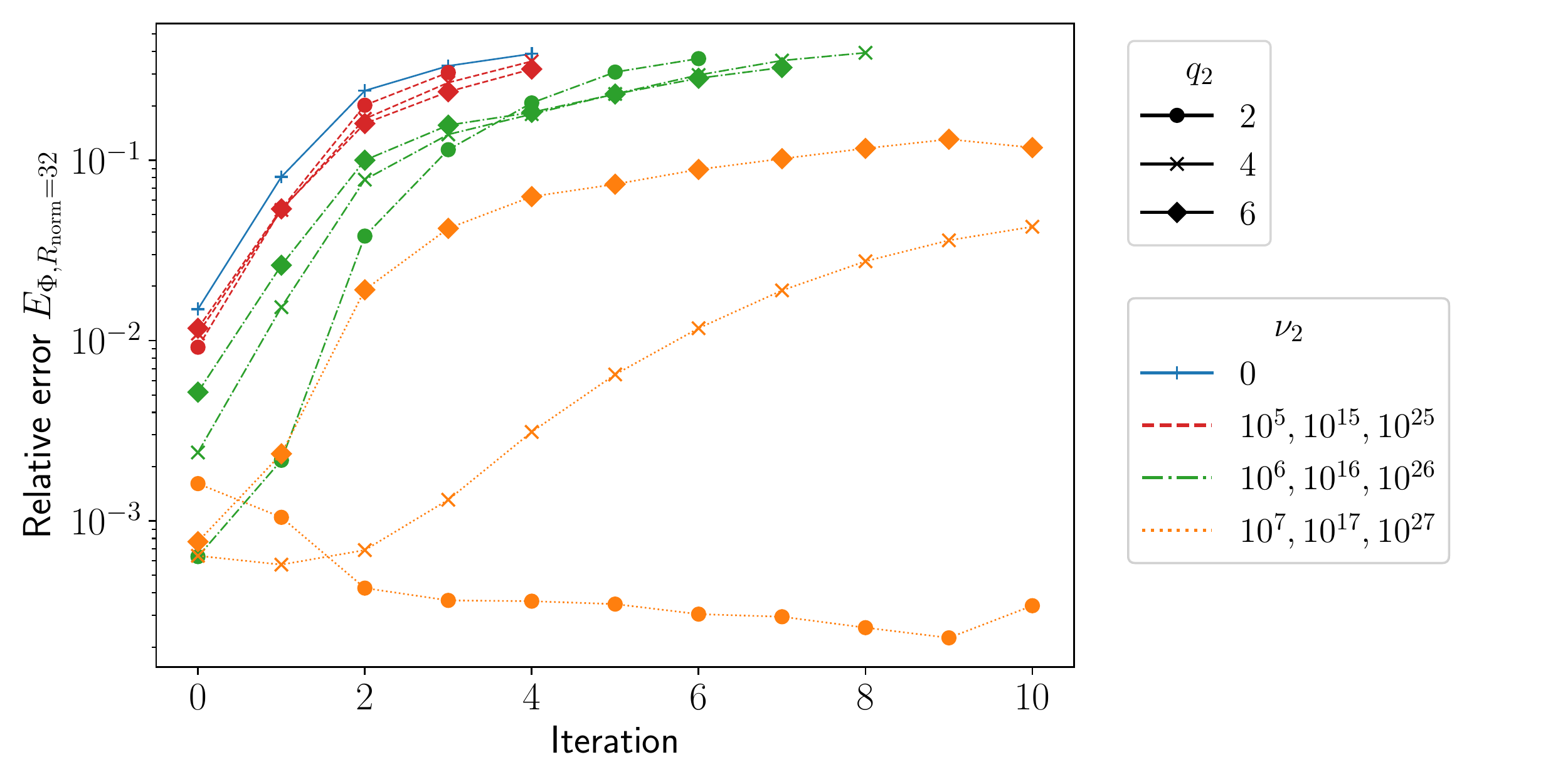}
        \caption{SL-SI-SETTLS}
        \label{fig:gaussian_bumps_errors_viscosity_SL_SI_SETTLS}
    \end{subfigure}
    \caption{Gaussian bumps test case: relative error $E_{\Phi, \rnorm}$ between the Pint and fine solutions at $t = T$ along iterations for $\rnorm = 32$, with IMEX (left) or SL-SI-SETTLS (right) being used as a coarse time-stepping scheme, as a function of the viscosity coefficients applied on the coarsest level. All simulations use configuration $(\nlevels, \cfactor, \nrelax, \Mcoarse) = (3,2,0,128)$. No viscosity is applied on the finest and intermediate levels. Viscosity coefficients expressed in $\visc{q}$.} \label{fig:gaussian_bumps_errors_viscosity}
\end{figure}

\indent \replaced[id=R3]{It has been observed that the convergence and stability behaviors of the simulations, both in the cases with IMEX and SL-SI-SETTLS as a coarse scheme, are determined mainly by the viscosity order and coefficient applied on the coarsest level, \ie $q_2$ and $\nu_2$, with only little influence of the values applied on the intermediate level. This observation indicates that it is critical to damp amplifications produced by the coarsest discretization, whose stability constraints are the most critical ones. Therefore, in order to avoid the presentation of a too large number of simulations (with both $(q_1,\nu_1)$ and $(q_2,\nu_2)$ varying), we fix $\nu_1 = 0$ (\ie, no viscosity on the intermediate level) and present the results for each pair $(q_2,\nu_2)$. Moreover, similar results were obtained in the cases $\rnorm = 32$ and $\rnorm = 128$, and we only present the former.}{Since these ranges of parameters lead to a large number of simulations, we compare the obtained relative errors at given iterations in Figures \ref{fig:gaussian_bumps_errors_viscosity_IMEX} and \ref{fig:gaussian_bumps_errors_viscosity_SL_SI_SETTLS}. The nine colored regions separated by black lines correspond to each pair $(q_1, q_2)$, and an additional row and column take into account the cases $\nu_1 = 0$ and $\nu_2 = 0$. Red regions indicate unavailable results due to instabilities that stopped the simulation before the respective iteration. We present the errors only for $\rnorm = 128$; similar results are obtained under $\rnorm = 32$, with slightly smaller errors in some configurations.}

\indent \replaced[id=R3]{In the simulations using IMEX as a coarse scheme (Figure \ref{fig:gaussian_bumps_errors_viscosity_IMEX}), we observe that large second-, fourth- or six-order viscosities are required to ensure stability: indeed, zero or too small viscosities lead to strong unstable behaviors from the initial iteration, and with intermediate values the simulations initially converge but instabilities are triggered after few iterations. The results also reveal that a compromise between stability and accuracy needs to be fulfilled when choosing the viscosity order: among the stable simulations, a faster convergence is obtained when the viscosity order is higher, indicating that only the largest wavenumbers need to be damped. Notably, with a large second-order viscosity, which damps a large region of the wavenumber spectrum, the error is much larger compared to higher orders already in the initial iteration. On the other hand, in the case where SL-SI-SETTLS is used as a coarse scheme (Figure \ref{fig:gaussian_bumps_errors_viscosity_SL_SI_SETTLS}), a much more critical dependence of the stability on the viscosity approach is observed. Only the simulation using a large second-order viscosity presents a convergent behavior, indicating that a large region of the wavenumber spectrum needs to be damped. Finally, the results for both coarse schemes reveal that, compared to reported values in serial simulations, larger viscosity coefficients (by between one and three orders of magnitude) are required to ensure stability in the PinT context.}{In the simulations using IMEX as coarse scheme (Figure \ref{fig:gaussian_bumps_errors_viscosity_IMEX}), we observe that a large second-order viscosity, either on the intermediate or on the coarsest discretization level, lead to a larger error in the initial iteration ($k = 0$; see the errors corresponding to $\nu_1 = 10^7 \visc{2}$ or $\nu_2 = 10^7 \visc{2}$), indicating that it produces an over-damping that affects the quality of the solution. It is more remarkable when $\nu_1$ is large since it is applied on the intermediate level, which is supposed to provide a better approximation to the fine solution than the coarsest level. However, we notice that large enough second-order viscosity coefficients are required to ensure stability along iterations; indeed, only simulations with $\nu_2 = 10^7 \visc{2}$ can reach iteration $k = 10$, independently of $\nu_1$; it indicates that a strong viscosity needs to be applied on the coarsest level to damp instabilities created by the large time step, and since these instabilities have already been damped, no viscosity is required on the intermediate level. On the other hand, large coefficient values of higher-order viscosities do not increase the initial error and are even able to slightly reduce it; moreover, higher-order viscosity approaches seem more effective in improving stability and convergence of the PinT simulation, with more simulations being able to reach advanced iterations, notably with sixth-order viscosities. As in the second-order case, we still observe better results when using large enough values for $\nu_2$ (\ie to damp instabilities on the coarsest level). Still, simulations using moderate values of fourth- and sixth-order viscosity can significantly reduce the error.}

\deleted[id=R3]{\indent Concerning the simulations using SL-SI-SETTLS on the coarse levels (Figure \ref{fig:gaussian_bumps_errors_viscosity_SL_SI_SETTLS}), we first notice that contrary to the configurations using IMEX, large viscosities usually reduce the error at $k = 0$, indicating that instabilities are already present in the initial iteration. The only exception is when a second-order $\nu_1$ is applied. Then, there is an over-damping of the solution on the intermediate level. However, the most noticeable conclusion is that only very few viscosity configurations can successfully control the more critical stability properties of the PinT methods. Indeed, only simulations with large enough viscosity values of second-, fourth- and sixth-order on the coarsest level can complete ten iterations. Almost all of them can reduce the error compared to the initial iteration using a second-order approach. As discussed above, instabilities of the PinT methods affect larger ranges of the wavenumber spectrum when using SL-SI-SETTLS instead of IMEX on the coarse levels; therefore, a less selective viscosity approach, damping both small and high wavenumbers seem to be required to make the parallel-in-time simulation stable and convergent. Finally, as observed in the simulations using IMEX, applying artificial viscosity is more critical on the coarsest level, in which the largest time step is used, with the results having little dependence on $\nu_1$.}

\begin{figure}[!htbp]
    \begin{subfigure}{.5\linewidth}
        \centering
        \includegraphics[scale=.35]{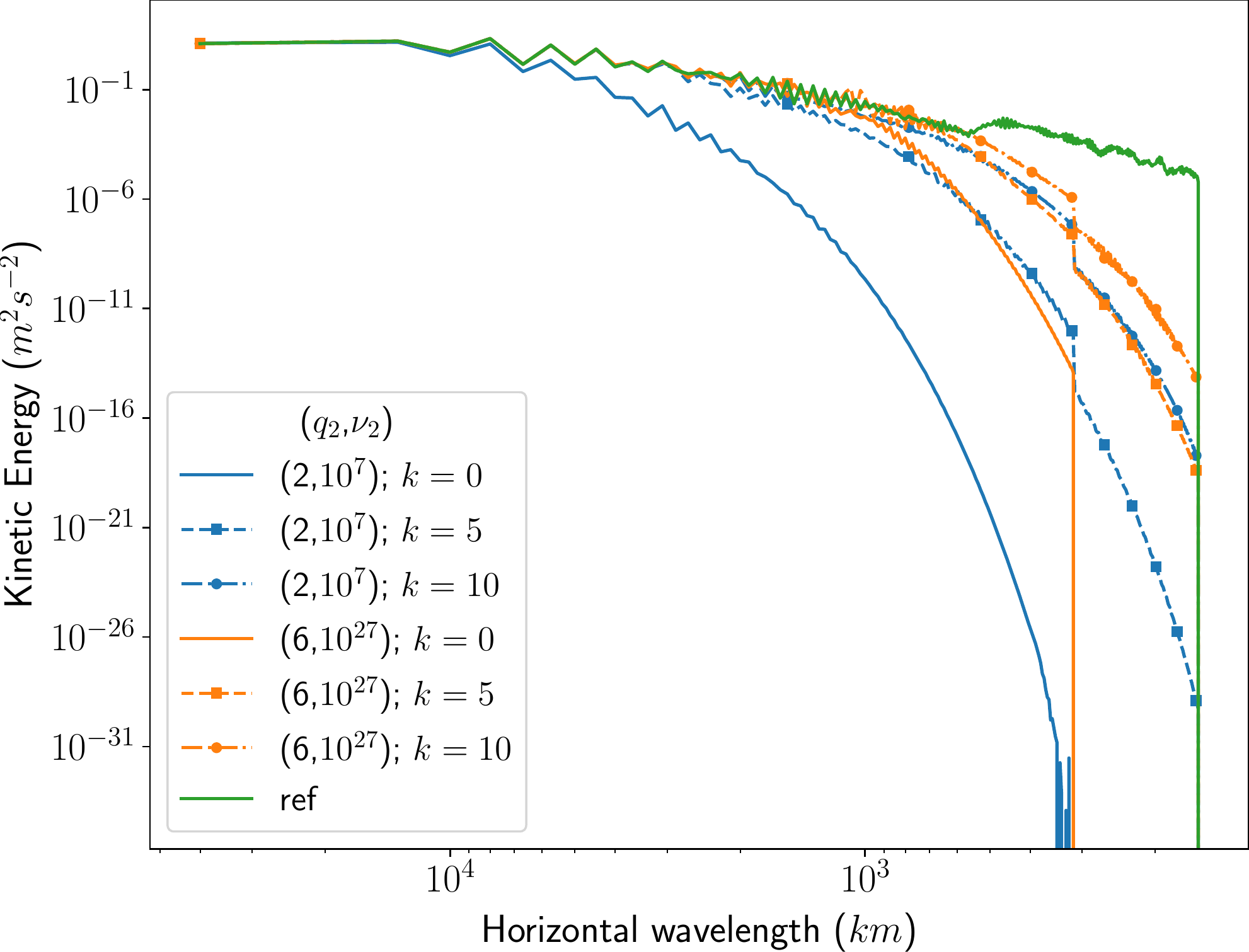}
        \caption{IMEX}
    \end{subfigure}
    \begin{subfigure}{.5\linewidth}
        \centering
        \includegraphics[scale=.35]{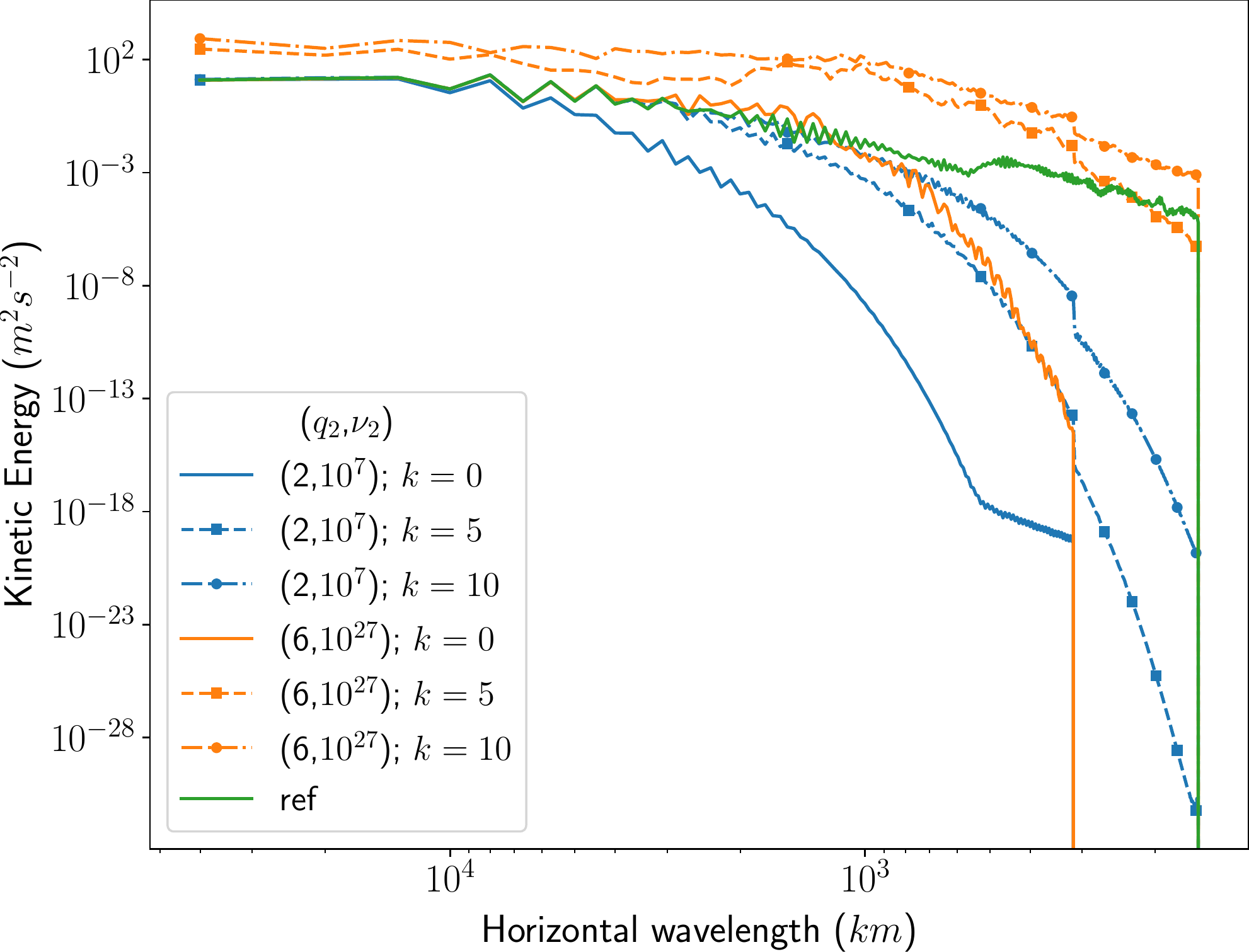}
        \caption{SL-SI-SETTLS}
    \end{subfigure}
    \caption{Gaussians bumps test case: kinetic energy spectra at $t = T$ of the reference (fine) and PinT solutions at iterations 0, 5 and 10 for chosen viscosity order and coefficient $(q_2, \nu_2)$ applied on the coarsest level. All simulations use configuration $(\nlevels, \cfactor, \nrelax, \Mcoarse) = (3,2,0,128)$, with no viscosity applied on the finest and intermediate levels. Left and right: IMEX and SL-SI-SETTLS as coarse time-stepping schemes, respectively.} 
    \label{fig:gaussian_bumps_spectrum_viscosity}
\end{figure}

\indent We illustrate in Figure \ref{fig:gaussian_bumps_spectrum_viscosity} the evolution of the kinetic energy spectrum for two chosen viscosity configurations, namely $(q_1,\nu_1,q_2,\nu_2) = (-,0,2,10^7)$ and $(q_1,\nu_1,q_2,\nu_2) = (-,0,6,10^{27})$ (\ie with no viscosity applied on the intermediate level and high second- or sixth-order viscosity on the coarsest none), for both time-stepping schemes on the coarse level. As shown in Figure \ref{fig:gaussian_bumps_errors_viscosity_IMEX}, both viscosity configurations provide convergence in the simulations using IMEX, mainly the sixth-order one. We indeed observe that the second-order viscosity strongly damps the initial PinT solution almost in the entire wavenumber spectrum, and the sixth-order one provides a better approximation to the reference spectrum. On the other hand, in the simulations using SL-SI-SETTLS on the coarse levels, the large damping due to the second-order viscosity is required for keeping the simulation stable along iterations; a higher-order approach, even with a very large coefficient, damps only the largest wavenumbers, and we observe an overamplification of medium wavenumbers already at iteration $k = 0$, which propagates to the entire spectrum after a few iterations.

\subsubsection{Evaluation of computing times and speedups}

\indent We now evaluate the computing times and respective speedups of some chosen PinT configurations presenting relatively stable and convergent behavior. A larger number of configurations using IMEX on the coarse levels is chosen since it provides better behavior in general, which allows studying the speedups as a  function of various MGRIT parameters, such as the number of levels, the relaxation strategy, the spectral resolution on the coarse levels and the artificial viscosity approach. In the case of SL-SI-SETTLS, a smaller set of configurations is chosen.

\begin{figure}[!htbp]
    \begin{subfigure}{.5\linewidth}
        \centering
        \includegraphics[scale=.475]{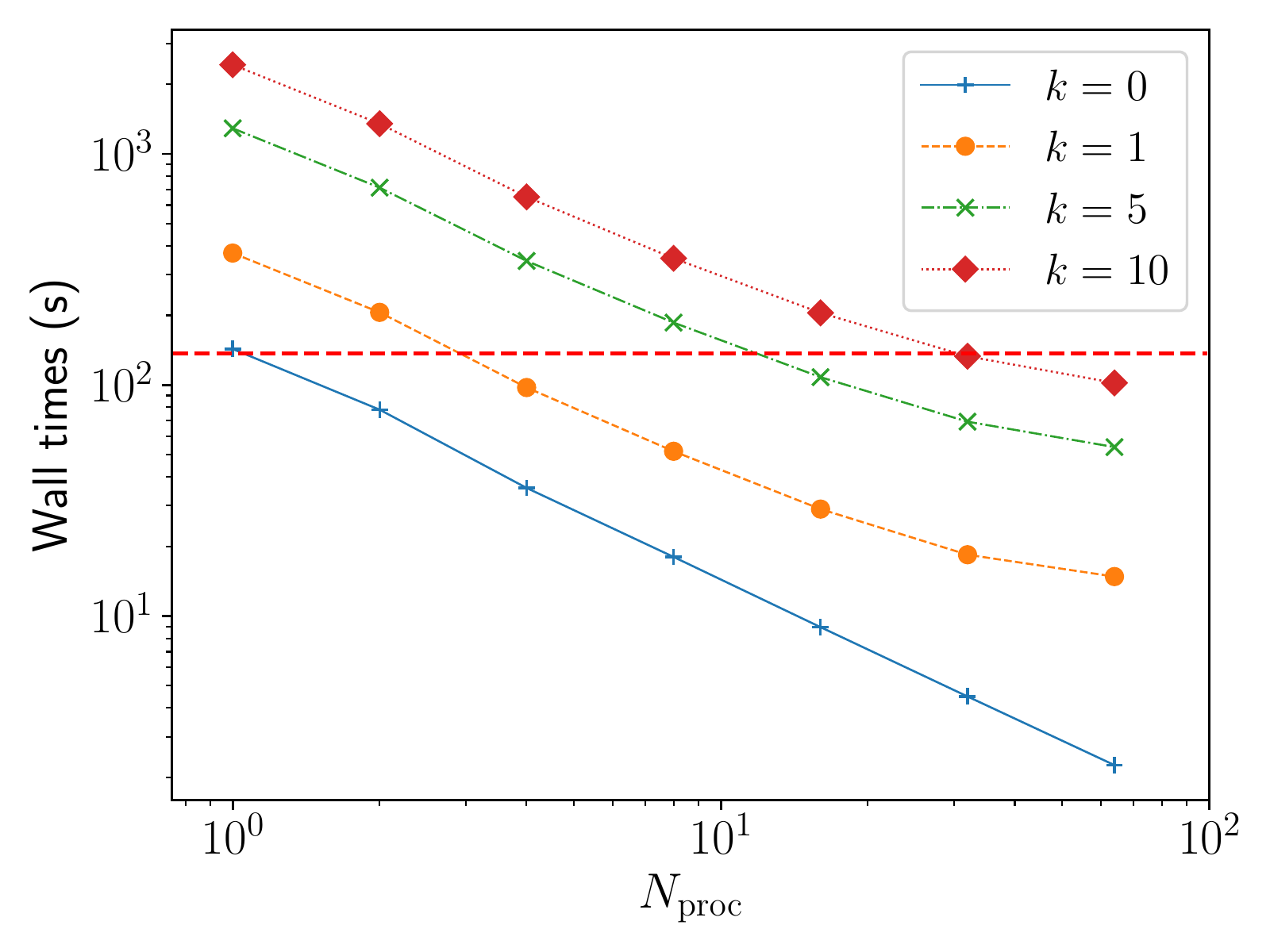}
        \caption{IMEX}
    \end{subfigure}
    \begin{subfigure}{.5\linewidth}
        \centering
        \includegraphics[scale=.475]{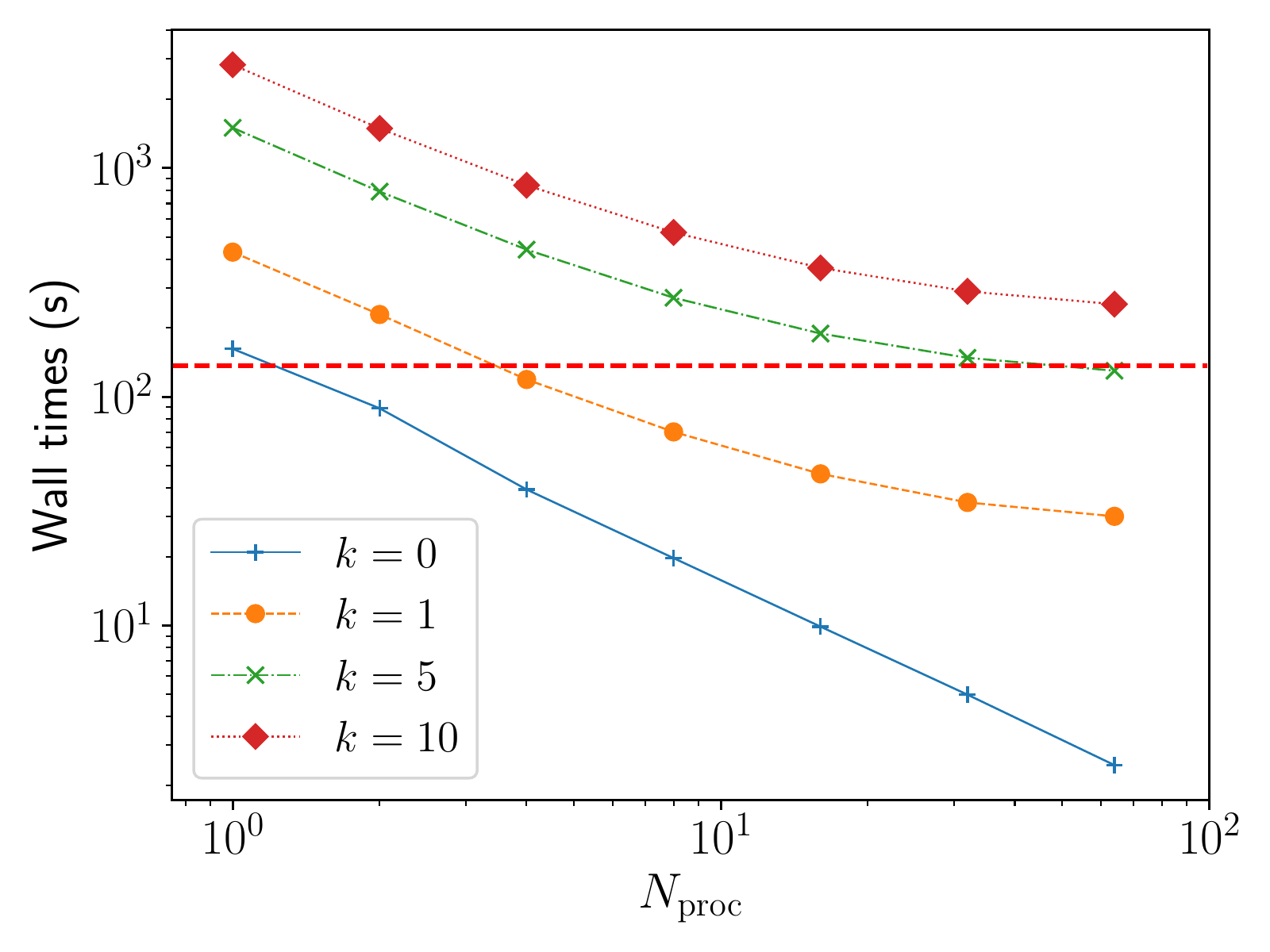}
        \caption{SL-SI-SETTLS}
    \end{subfigure}
    \caption{Gaussian bumps test case: wall times for completing $k$ iterations of PinT simulations for the configuration $(\nlevels, \cfactor, \nrelax, \Mcoarse, q_1, \nu_1) = (2,2,0,51,2,10^6)$ using IMEX (left) or SL-SI-SETTLS (right) on the coarse level as a function of the number of parallel processors in time. The horizontal, dashed line indicates the computing time of the reference solution, computed using $\Dt_{\text{ref}} = \Dt_0 = 60\text{s}$.}
    \label{fig:gaussian_bumps_strong_scaling}
\end{figure}

\indent We are mainly interested in the relation between the speedups and the errors provided by each PinT simulation. However, we first present some strong scaling results in Figure \ref{fig:gaussian_bumps_strong_scaling} for the configuration $(\nlevels, \cfactor,\\ \nrelax, \Mcoarse, q_1, \nu_1) = (2,2,0,51,2,10^6)$ using IMEX or SL-SI-SETTLS as coarse time-stepping schemes. The wall times for reaching given iterations are compared to the reference one ($T_{\text{ref}} \approx 137 \text{s}$) for simulations using $\nproc \in [1,64]$. In the initial iteration, we observe a good scaling in the entire range of $\nproc$; in the following ones, we begin to notice a saturation of the speedup for $\nproc \geq 32$, but there seems to still exist some room for further improvements if more processors are used. In all results presented hereafter, we consider $\nproc = 64$. We notice that the number of processors considered here is smaller than the theoretical maximum value of $\nproc$, which would still provide speedups, which depends on the number of fine time steps, the coarsening factor and the relaxation strategy, as described in Section \ref{subsec:MGRIT}; for instance, in the simulation depicted in Figure \ref{fig:gaussian_bumps_strong_scaling}, up to $\nproc = 1080$ could be used. This limit is smaller in other configurations considered in this work but remains larger than $\nproc = 64$ when $\nrelax = 0$ is used. As a last remark, we notice in Figure \ref{fig:gaussian_bumps_strong_scaling} that, under the same PinT configurations, the simulations using SL-SI-SETTLS as a coarse scheme are more expensive than those using IMEX, which may be due to additional costs linked to the semi-Lagrangian approach (trajectory estimations and interpolations to departure points).

\begin{figure}[!htbp]
    \begin{subfigure}{.5\linewidth}
        \centering
        \includegraphics[scale=.5]{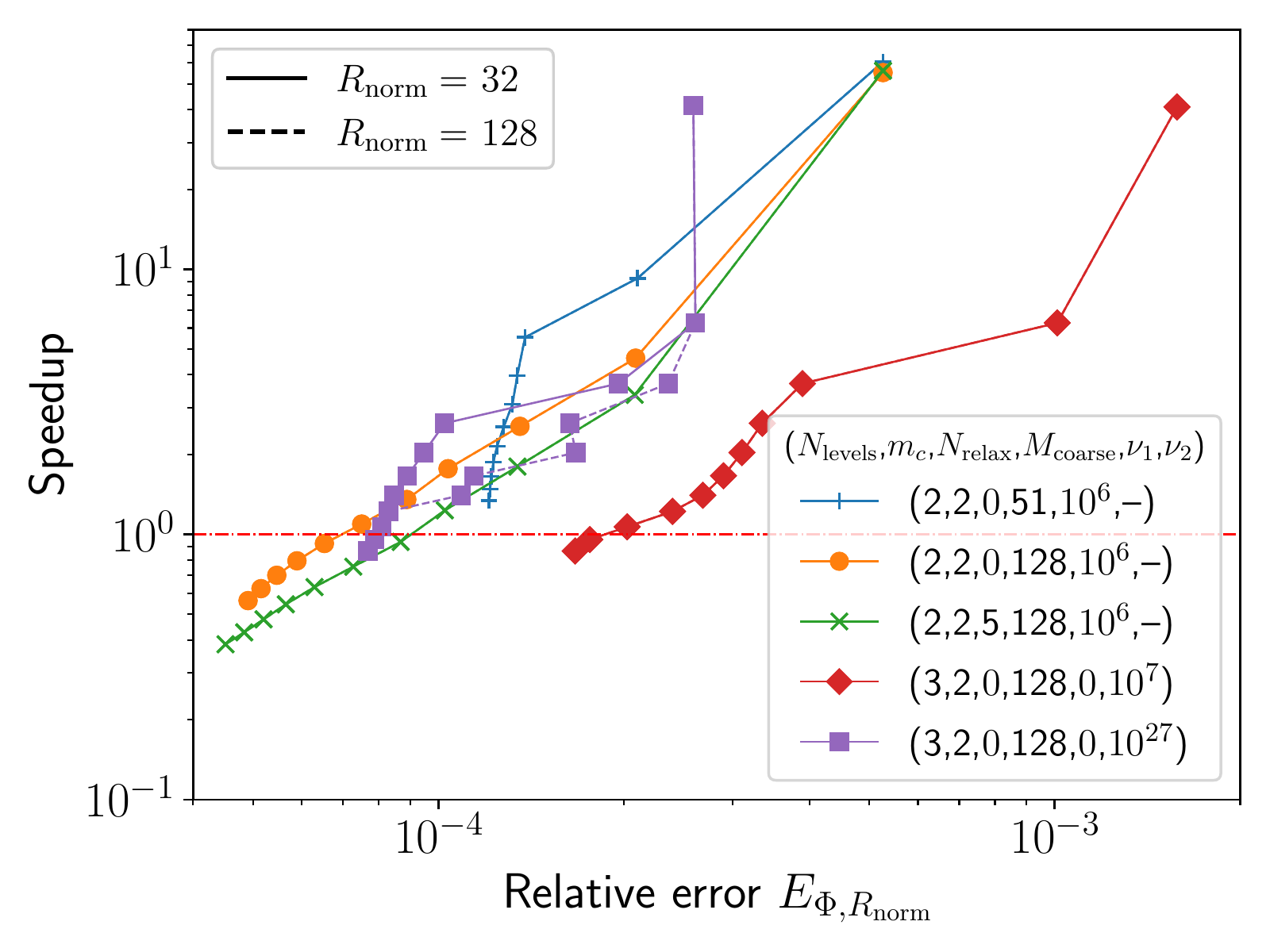}
        \caption{IMEX\label{fig:gaussian_bumps_speedup_IMEX}}
    \end{subfigure}
    \begin{subfigure}{.5\linewidth}
        \centering
        \includegraphics[scale=.5]{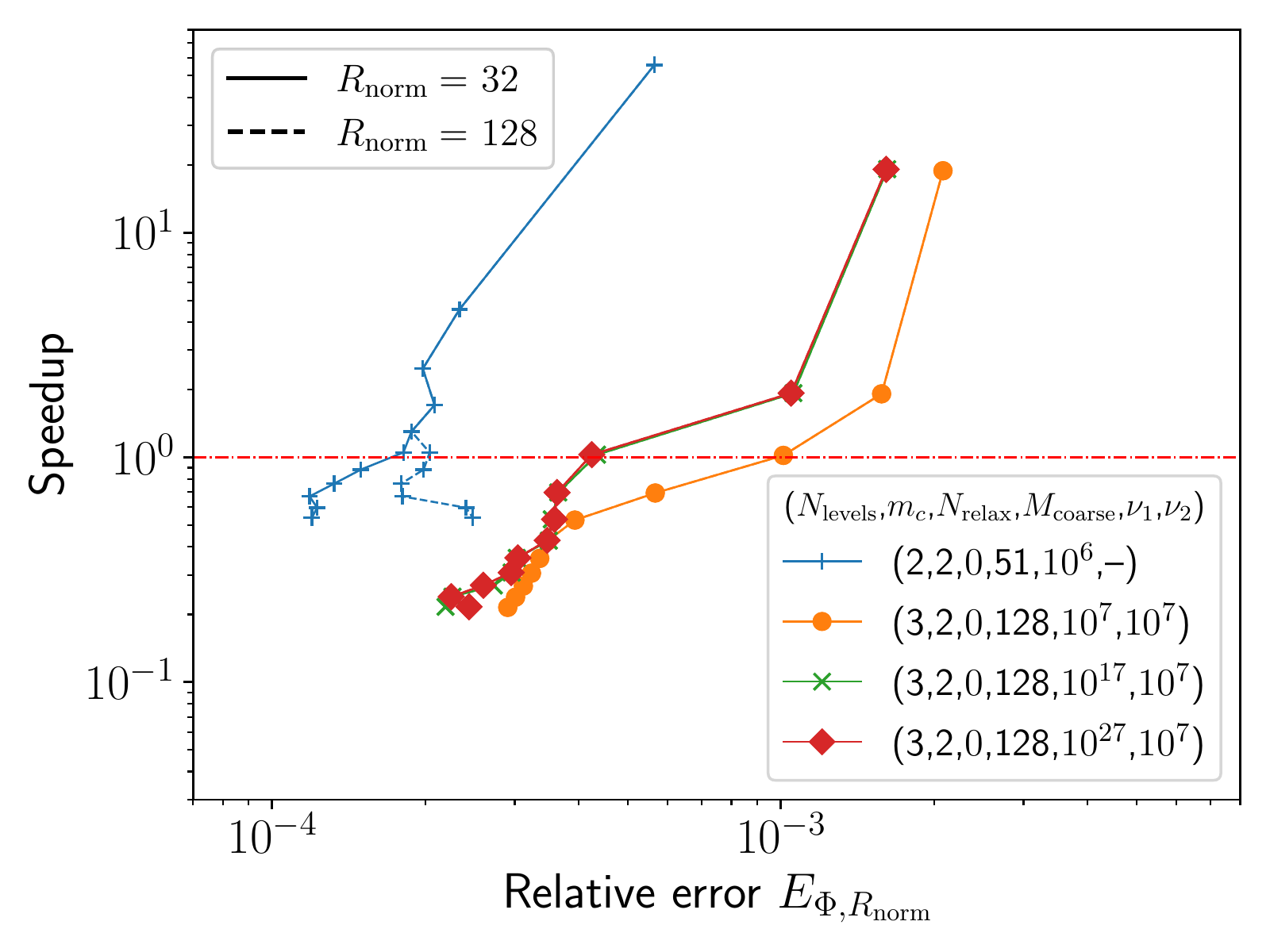}
        \caption{SL-SI-SETTLS\label{fig:gaussian_bumps_speedup_error_SL_SI_SETTLS}}
    \end{subfigure}
    \caption{Gaussian bumps test case: speedup as a function of the relative geopotential error in the spectral space for chosen configurations using IMEX (left) or SL-SI-SETTLS (right) as coarse time-stepping scheme. Simulations are identified by $(\nlevels, \cfactor, \nrelax, \Mcoarse, \nu_1, \nu_2)$ and each data point corresponds to an iteration. The horizontal, dashed-dotted line indicates a unitary speedup. All simulations use $\nproc = 64$ parallel processors in time.}
    \label{fig:gaussian_bumps_speedup_error}
\end{figure}

We now study the speedups \wrt the relative error in the geopotential field for simulations using IMEX in Figure \ref{fig:gaussian_bumps_speedup_IMEX}. The presented configurations are chosen to compare the influence of the PinT parameters individually. First, we observe that the use of a less refined resolution $\Mcoarse = 51$ reduces the cost of the time integration on the coarse level, leading to a speedup of approximately 5.5 for reducing the error by a factor close to 5; however, the convergence rapidly stagnates, with increasing computational times not leading to significant error reductions; if one wants to obtain smaller errors, it is necessary to use a larger coarse spectral resolution, which increases the computational cost of the temporal parallelization but leads to speedups still larger than one if enough processors are used. Concerning the relaxation strategy, it has been seen in Figure \ref{fig:gaussian_bumps_errors_params_mgrit_IMEX} that only slight improvements in convergence are obtained by increasing $\nrelax$. A more expensive relaxation, however, strongly increases the computational cost of the MGRIT simulation since a larger number of time steps needs to be computed per iteration; therefore, increasing $\nrelax$ negatively impacts the speedup for obtaining approximately the same errors, and $\nrelax = 0$ provides the better compromise. Finally, we compare the influence of the viscosity orders and coefficients applied on the coarse levels: the worst compromise between speedup and convergence is obtained by using a too-large second-order viscosity, which overdamps the entire wavenumber spectrum and leads to larger errors despite the convergent behavior. Better results are obtained using a large higher-order viscosity, mainly in the beginning of the wavenumber spectrum, with a speedup of approximately 2.6 to obtain relative errors close to $10^{-4}$ under $\rnorm = 32$. Still, the overall trade-off between convergence and computational cost is similar to moderate second-order viscosity configurations.

\indent Figure \ref{fig:gaussian_bumps_speedup_error_SL_SI_SETTLS} presents the speedups as a function of the errors of some of the few relatively stable simulations using SL-SI-SETTLS on the coarse levels, both with the default moderate second-order viscosity or other viscosity configurations. In the former case, only the less aggressive configuration $(\nlevels, \cfactor) = (2,2)$ remains stable and can improve the convergence with speedups larger than the unity. Still, after two iterations, the convergence deteriorates at the end of the wavenumber spectrum. In the other simulations, a better convergence behavior is obtained due to the larger second-order viscosity on the coarsest level; however, due to this same reason, mainly when a large second-order viscosity is also used on the intermediate level, the initial errors are larger than with $\nu_1 = 10^6 \visc{2}$ and even after ten iterations, they do not provide smaller errors than this case; moreover, the speedups are considerably smaller due to the larger spectral resolution on the coarse levels and go below the unity after two iterations.

\subsection{Unstable jet test case}

\begin{figure}[!htbp]
    \centering
    \includegraphics[scale=.5]{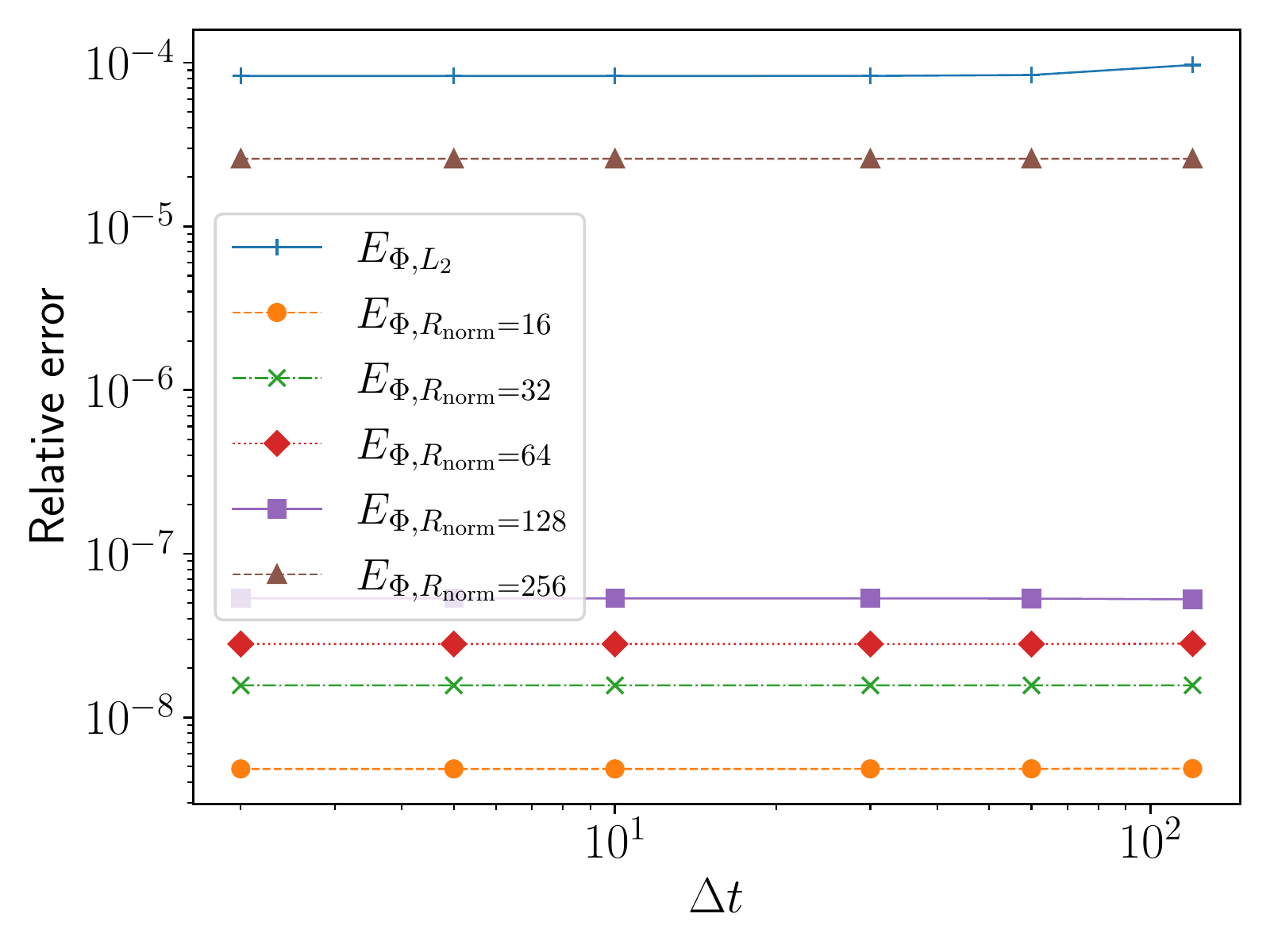}
    \caption{Unstable jet test case: relative $L_2$ error (computed in the physical space) and spectral errors (for various values of $\rnorm$) between a solution obtained with spectral resolution $M = 512$ and time step $\Dt = 2$ and solutions obtained with $M = M_0 = 256$ and various time steps. IMEX is used in all cases.}
    \label{fig:err_ref_unstable_jet}
\end{figure}

\begin{figure}[!htbp]
    \begin{subfigure}{.5\linewidth}
        \centering
        \includegraphics[scale=.375]{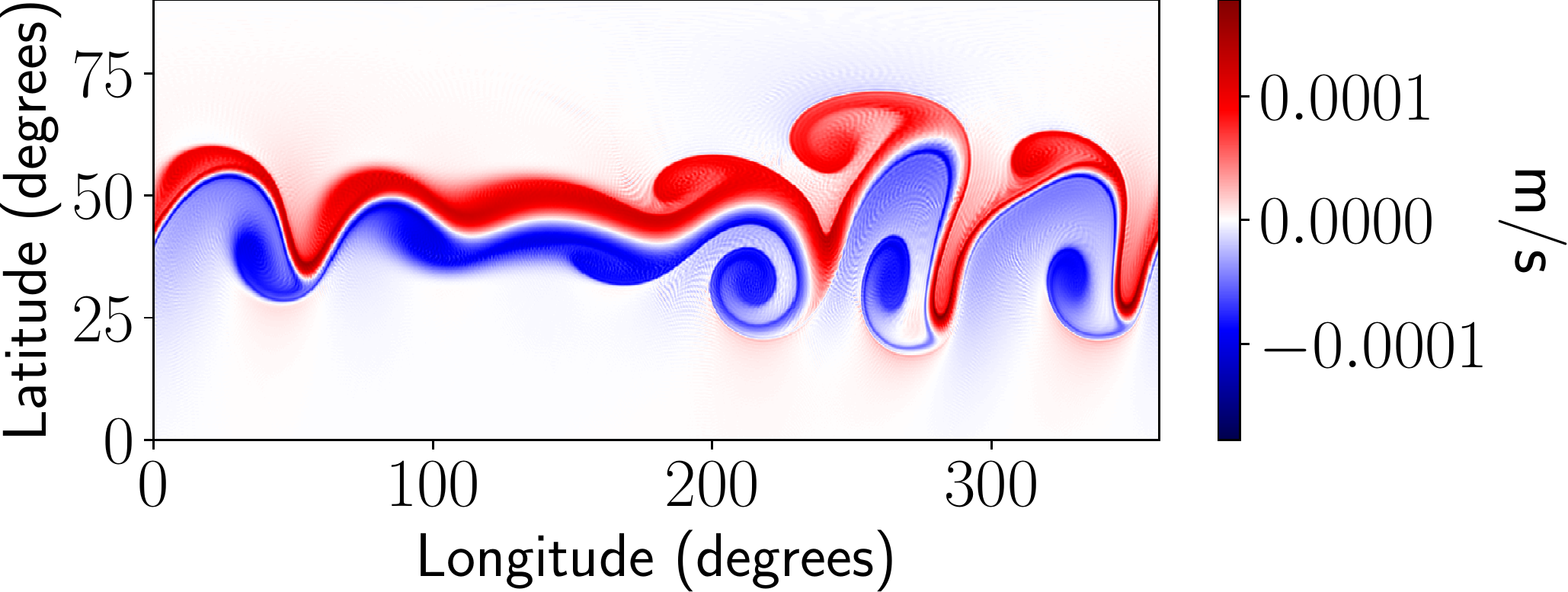}
        \caption{$M = 512$, $\Dt = 2$}
    \end{subfigure}
    \begin{subfigure}{.5\linewidth}
        \centering
        \includegraphics[scale=.375]{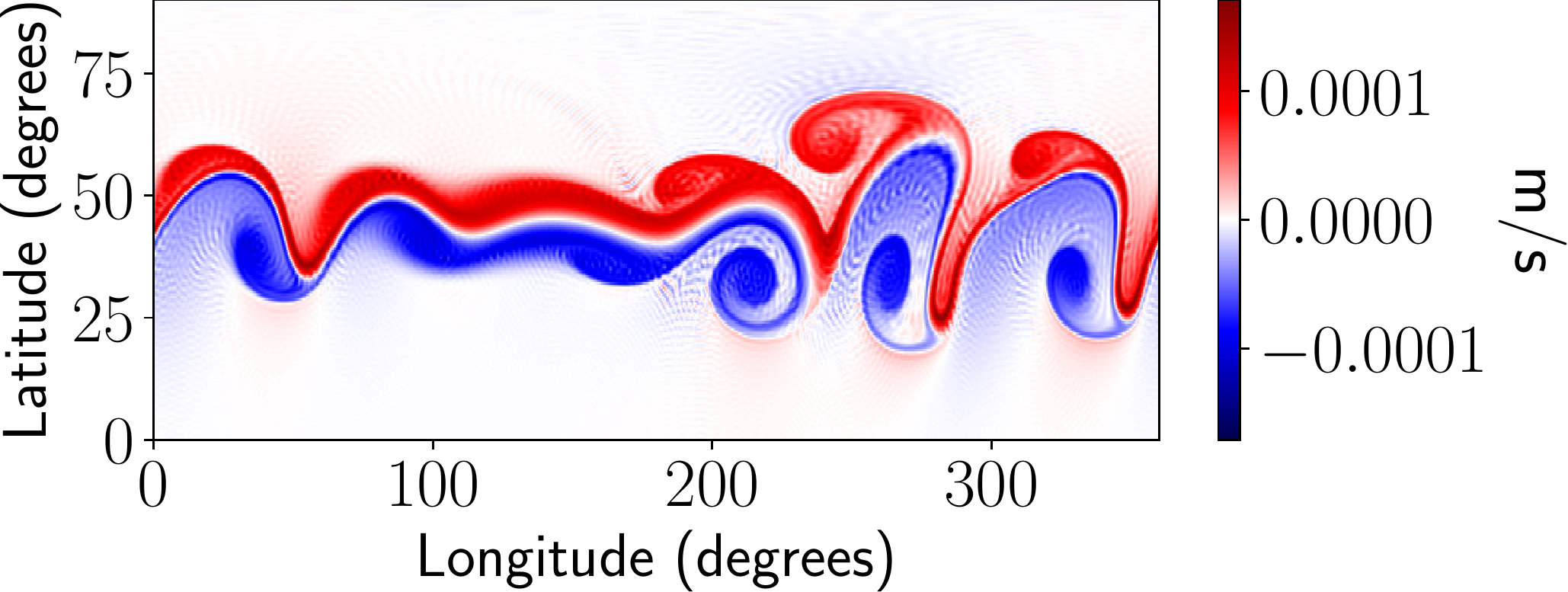}
        \caption{$M = 256$, $\Dt = 60$}
    \end{subfigure}
    \caption{Unstable jet test case: solution at $t = T = 144\text{h}$ computed with IMEX. Zoom on the north hemisphere, the vorticity field being negligible in the south one.}
    \label{fig:ref_unstable_jet}
\end{figure}

\indent We now consider the test case presented by \cite{galewski_al:2004}. In this test, a stationary zonal jet is perturbed by a Gaussian bump in the geopotential field, leading to the formation of vortices and a rapid energy transfer from low to high wavenumbers, which may be especially challenging in the context of PinT methods. This test case was used by \cite{hamon_al:2020} for studying the temporal parallelization of the SWE on the rotating sphere using PFASST, and we consider the same simulation length in time, namely $T = 144\text{h}$. As in the Gaussian bumps test case, we begin by choosing a temporal discretization size to be used as a reference and fine solution for the study performed here. We perform simulations using a spectral resolution $M_0 = 256$, timestep $\Dt_0 \in [2,960] \text(s)$ and no artificial viscosity ($\nu_0 = 0$) and evaluate the errors of the geopotential field both in physical and spectral spaces \wrt to a solution computed using $M = 512$ and $\Dt = 2$, also computed without viscosity. In all cases, the integration is performed using IMEX. Figure \ref{fig:err_ref_unstable_jet} shows that the errors for $\Dt_0 \leq 120 \text{s}$ are mainly due to the spatial discretization, with no visible dependence of the spectral errors on $\Dt_0$, for all $\rnorm$ values, and a slight increase of the physical error between $\Dt_0 = 60 \text{s}$ and $\Dt_0 = 120 \text{s}$. For all tested values $\Dt_0 > 120 \text{s}$, instabilities develop and the simulations are not able to reach $t = T$. Thus, we choose $\Dt_0 = 60 \text{s}$, allowing, in the context of Parareal and MGRIT, a stable simulation on the finest level and the use of not-too-large time steps on the coarse ones. Figure \ref{fig:ref_unstable_jet} presents the final vorticity field produced by the chosen fine discretization and by the simulation using $M = 512$ and $\Dt = 2$.

\subsubsection{Convergence study}

\indent For the study of convergence and stability of Parareal and MGRIT, we consider the same set of parameters as in the Gaussian bumps test case, namely with $\nlevels \in \{2,3\}$, $\cfactor \in \{2,4\}$, $\nrelax \in \{0,1,5\}$ and $\Mcoarse \in \{51, 128\}$. However, instead of using a fixed second-order artificial viscosity approach, with the same coefficient applied on all coarse levels, we use the results obtained in Section \ref{subsec:gaussian_bump_viscosity} in order to improve the performance of the temporal parallelization. In the simulations using IMEX on the coarse levels, we consider a higher-order viscosity; however, due to the higher complexity of the unstable jet test case compared to the Gaussian bumps one, with a larger temporal domain and possibly more challenging stability, we choose here to use fourth instead of sixth-order viscosity; moreover, the viscosity coefficients are chosen as a function of the time step size used on each level, such that larger coefficients are applied when stability constraints are more restrictive. In the case of SL-SI-SETTLS, we apply large second-order viscosity coefficients on all coarse levels, independently of their discretizations, due to more critical stability behavior. Table \ref{tab:unstable_jet_viscosity_coefficients} summarizes the adopted viscosity coefficients.

\begin{table}[!htbp]
    \centering
    \begin{tabular}{|c|c|c|}
        \cline{1-3}
        Coarse scheme & \multicolumn{1}{c|}{IMEX} & \multicolumn{1}{c|}{SL-SI-SETTLS} \\
        \cline{1-3}
        Viscosity order $q$ & \multicolumn{1}{c|}{4} & \multicolumn{1}{c|}{2}\\
        \cline{1-3}
        $\Dt = 2\Dt_0$ & $10^{16}$ & $10^7$ \\
        \cline{1-3}
        $\Dt = 4\Dt_0$ & $10^{17}$ & $10^7$ \\
        \cline{1-3}
        $\Dt = 16\Dt_0$  & $10^{17}$ & $10^7$ \\
        \cline{1-3}
    \end{tabular}
    \caption{Unstable jet test case: viscosity coefficients (in $\visc{q}$) applied on the coarse levels of the PinT simulations as a function of the time step and coarse time-stepping scheme.}
    \label{tab:unstable_jet_viscosity_coefficients}
\end{table}

\begin{figure}[!htbp]
    \begin{subfigure}{.5\linewidth}
        \centering
        \includegraphics[scale=.425]{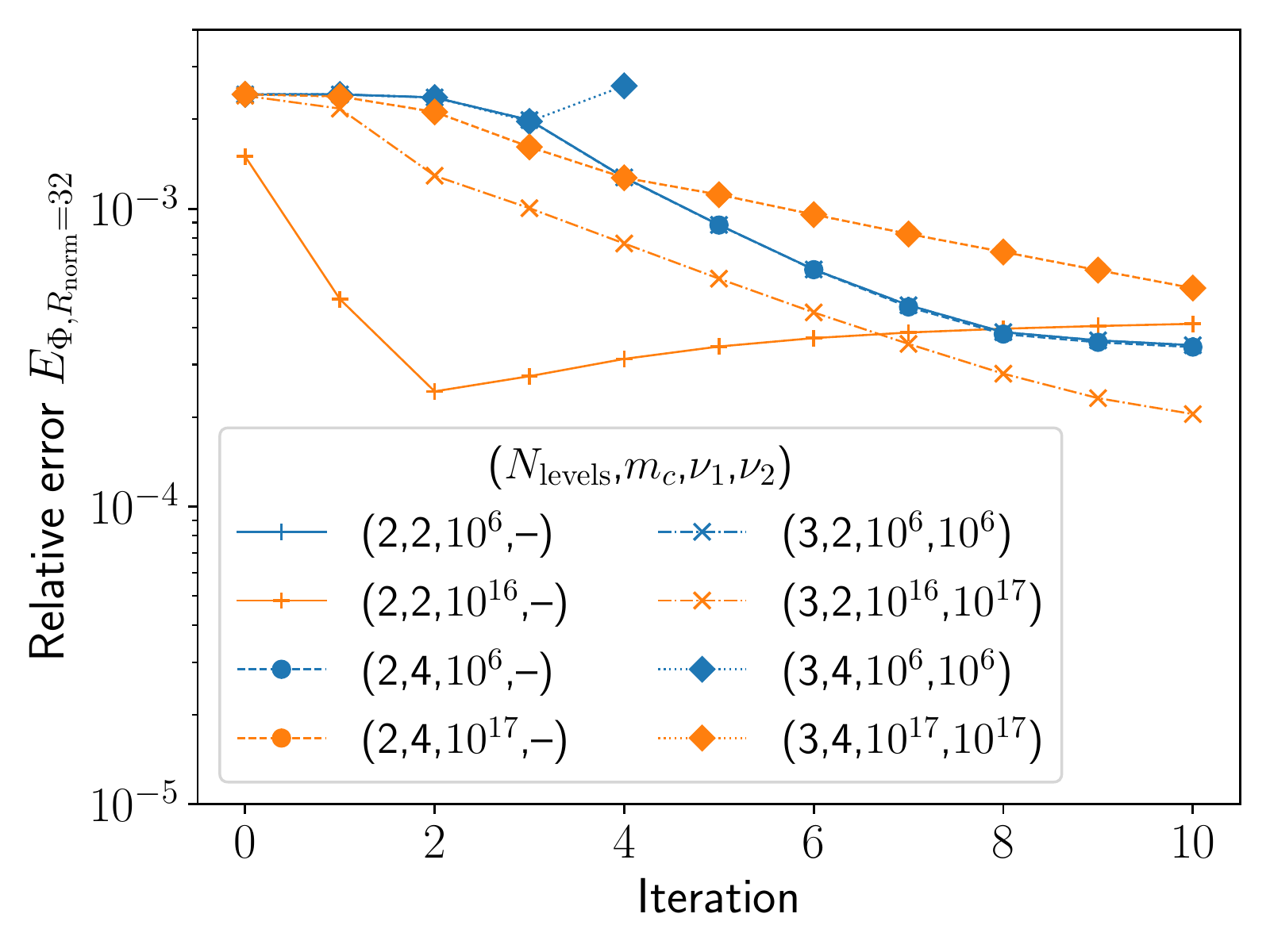}
        \caption{$M_{\text{coarse}} = 51$, $\rnorm = 32$}
    \end{subfigure}
    \begin{subfigure}{.5\linewidth}
        \centering
        \includegraphics[scale=.425]{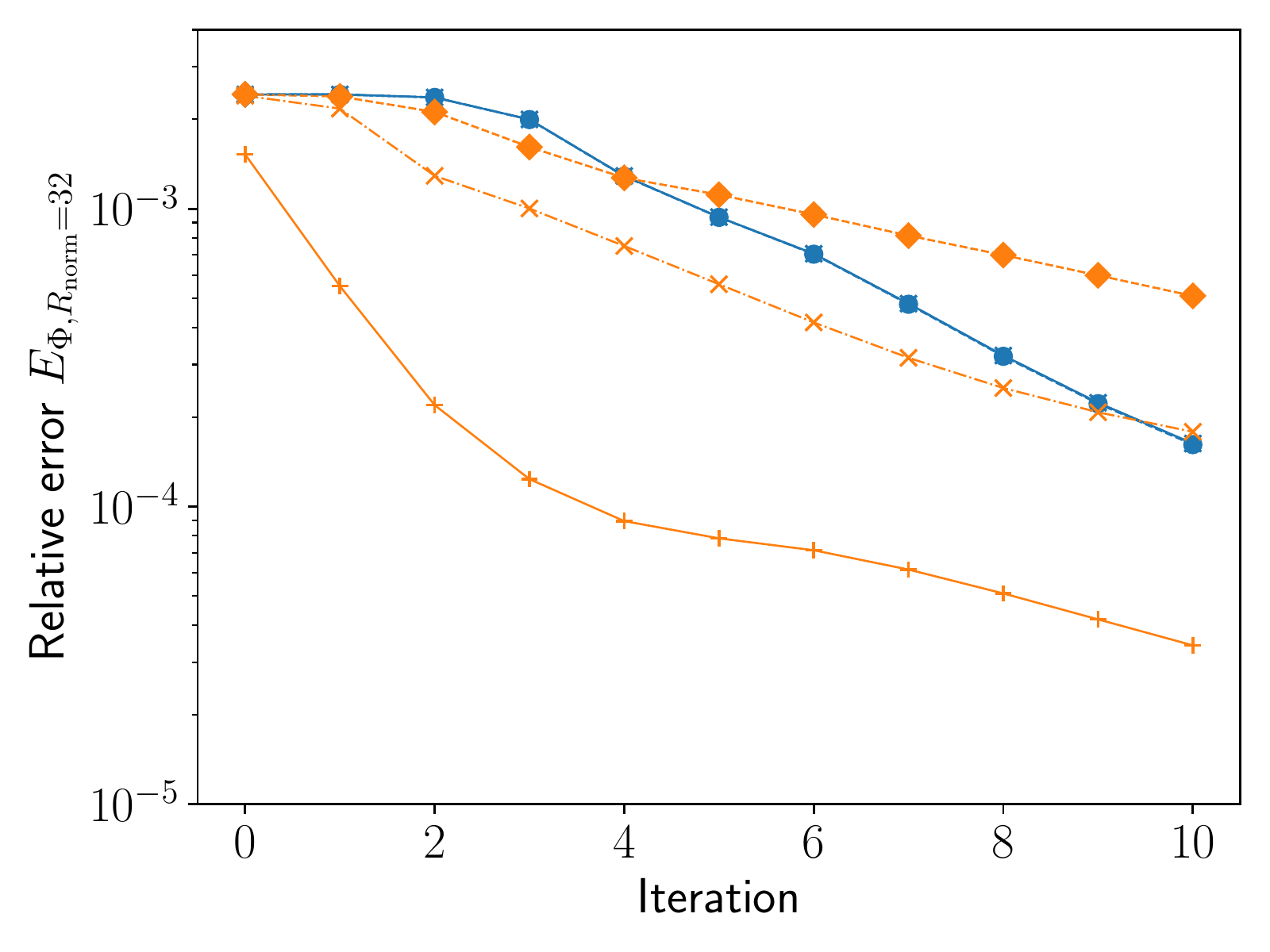}
        \caption{$M_{\text{coarse}} = 128$, $\rnorm = 32$}
    \end{subfigure}
    \caption{Unstable jet test case: relative error $E_{\Phi,\rnorm = 32}$ between the PinT and fine solutions at $t = T$ along iterations for $\Mcoarse = 51$ (left) and $\Mcoarse = 128$ (right), with IMEX used on the coarse levels. Results are identical under $\rnorm = 128$. Simulations are identified by $(\nlevels, \cfactor, \nu_1, \nu_2)$. The curves corresponding to the second-order viscosity visually coincide, as well as the curves $(\nlevels, \cfactor) = (2,4)$ and $(\nlevels, \cfactor) = (3,4)$ in the fourth-order viscosity case.}
    \label{fig:unstable_jet_errors_params_mgrit_IMEX}
\end{figure}

\indent Figure \ref{fig:unstable_jet_errors_params_mgrit_IMEX} presents the evolution of the PinT errors along iterations when IMEX is used on the coarse levels. As in the Gaussian bumps test case, only slight improvements are observed by using more expensive relaxation strategies; therefore, we only present the results for $\nrelax = 0$. Moreover, the observed stability and convergence behaviors are the same in all spatial scales, with the large ones dominating the error magnitudes, such that the convergence curves are identical under $\rnorm = 32$ and $\rnorm = 128$; we, therefore, present only the former case (these same remarks being valid for the results using SL-SI-SETTLS). Finally, in order to evaluate the proposed viscosity approach, we also present the results using the same viscosity order and coefficients initially adopted in the previous test case ($q = 2$ and $\nu = 10^6 \visc{2}$ on all coarse levels).

\indent Under the second-order viscosity approach, convergence, and stability is effectively obtained, except for the most aggressive configuration ($(\nlevels,\cfactor) = (3,4)$), which stops after three and one iteration, respectively in the cases $\Mcoarse = 51$ and $\Mcoarse = 128$, due to instabilities; however, all other simulations present an almost identical convergence behavior: it probably indicates that the second-order viscosity causes too large damping on the coarse levels and only a few contributions to the fine solution come from them. Indeed, proper choices of fourth-order viscosity coefficients provide better results: a much faster convergence is obtained in the configuration $(\nlevels,\cfactor,\Mcoarse) = (2,2,128)$, in which a moderate viscosity coefficient ($\nu_1 = 10^{16} \visc{4}$) is applied on the coarse level; a slower one (but still faster compared to the second-order viscosity) is obtained under $(\nlevels,\cfactor) = (3,2)$, in which a large viscosity ($\nu_2 = 10^{17} \visc{4}$) is applied only on the coarsest level; however, all simulations using this large viscosity on all coarse levels ($(\nlevels,\cfactor) = (2,4)$ and $(\nlevels,\cfactor) = (3,4)$) converge almost identically and slower than the second-order viscosity case, indicating an excessive damping. An unexpected result concerns the less aggressive configuration ($(\nlevels,\cfactor, \Mcoarse) = (2,2,51)$), which diverges after a very fast convergence in the first two iterations, probably indicating the triggering of numerical instabilities.

\begin{figure}[!htbp]
    \begin{subfigure}{.5\linewidth}
        \centering
        \includegraphics[scale=.425]{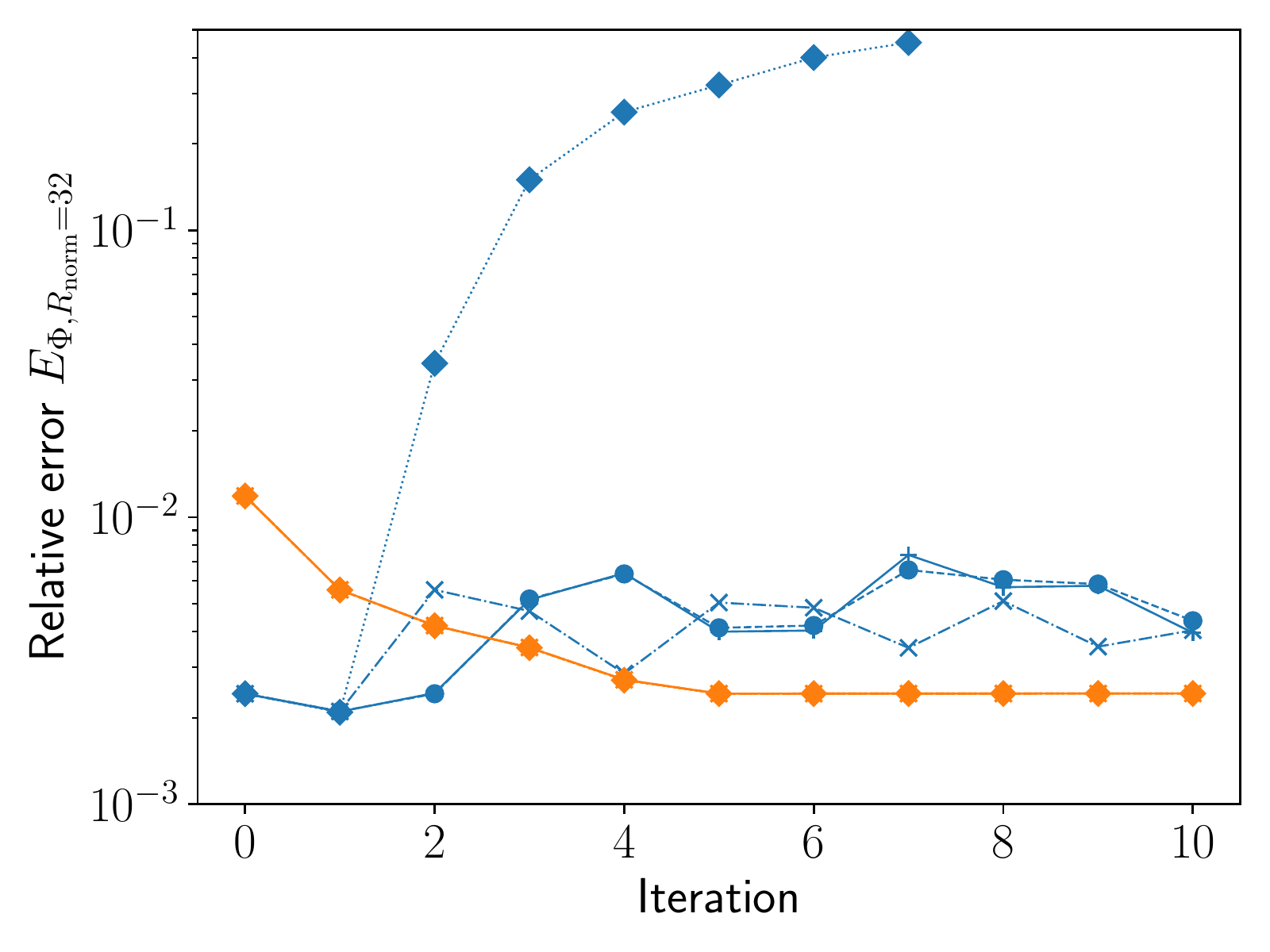}
        \caption{$M_{\text{coarse}} = 51$, $\rnorm = 32$}
    \end{subfigure}
    \begin{subfigure}{.5\linewidth}
        \centering
        \includegraphics[scale=.425]{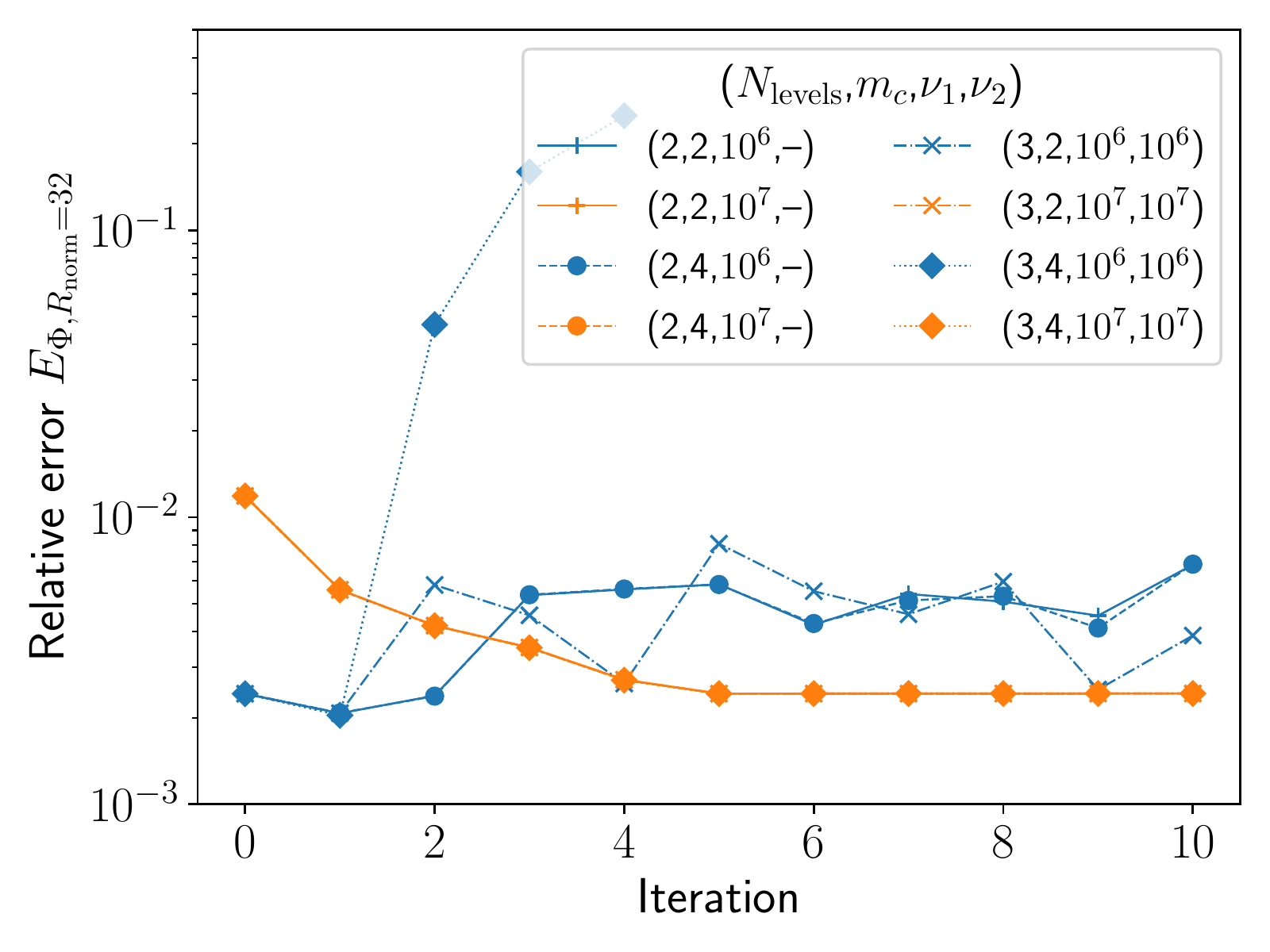}
        \caption{$M_{\text{coarse}} = 128$, $\rnorm = 32$}
    \end{subfigure}
    \caption{Unstable jet test case: relative error $E_{\Phi,\rnorm = 32}$ between the PinT and fine solutions at $t = T$ along iterations, for $\Mcoarse = 51$ (left) and $\Mcoarse = 128$ (right), with SL-SI-SETTLS used on the coarse levels. Results are identical under $\rnorm = 128$. Simulations are identified by $(\nlevels, \cfactor, \nu_1, \nu_2)$. The curves corresponding to $\nu_1 = 10^7 \visc{2}$ and $\nu_2 = 10^7 \visc{2}$ visually coincide.}
    \label{fig:unstable_jet_errors_params_mgrit_SL_SI_SETTLS}
\end{figure}

\indent The importance of the choice of the viscosity approach is even clearer in the simulations using SL-SI-SETTLS (Figure \ref{fig:unstable_jet_errors_params_mgrit_SL_SI_SETTLS}). The second-order viscosity with a moderate coefficient ($\nu_{\text{coarse}} = 10^6 \visc{2}$) is not able to ensure convergence and stability in any of the PinT configurations. For that, it is necessary to increase the viscosity coefficient to $\nu_{\text{coarse}} = 10^7 \visc{2}$. However, it leads to excessive damping on the coarse levels, and all configurations present the same convergence behavior, which stagnates around errors one to two orders of magnitude larger than those obtained using IMEX as a coarse scheme. Therefore, the viscosity setting does allow improving the PinT simulations using SL-SI-SETTLS, but it seems quite challenging to find a sweet spot between accuracy on the coarse levels and convergence.

\begin{figure}[!htbp]
    \begin{subfigure}{.5\linewidth}
        \centering
        \includegraphics[scale=.35]{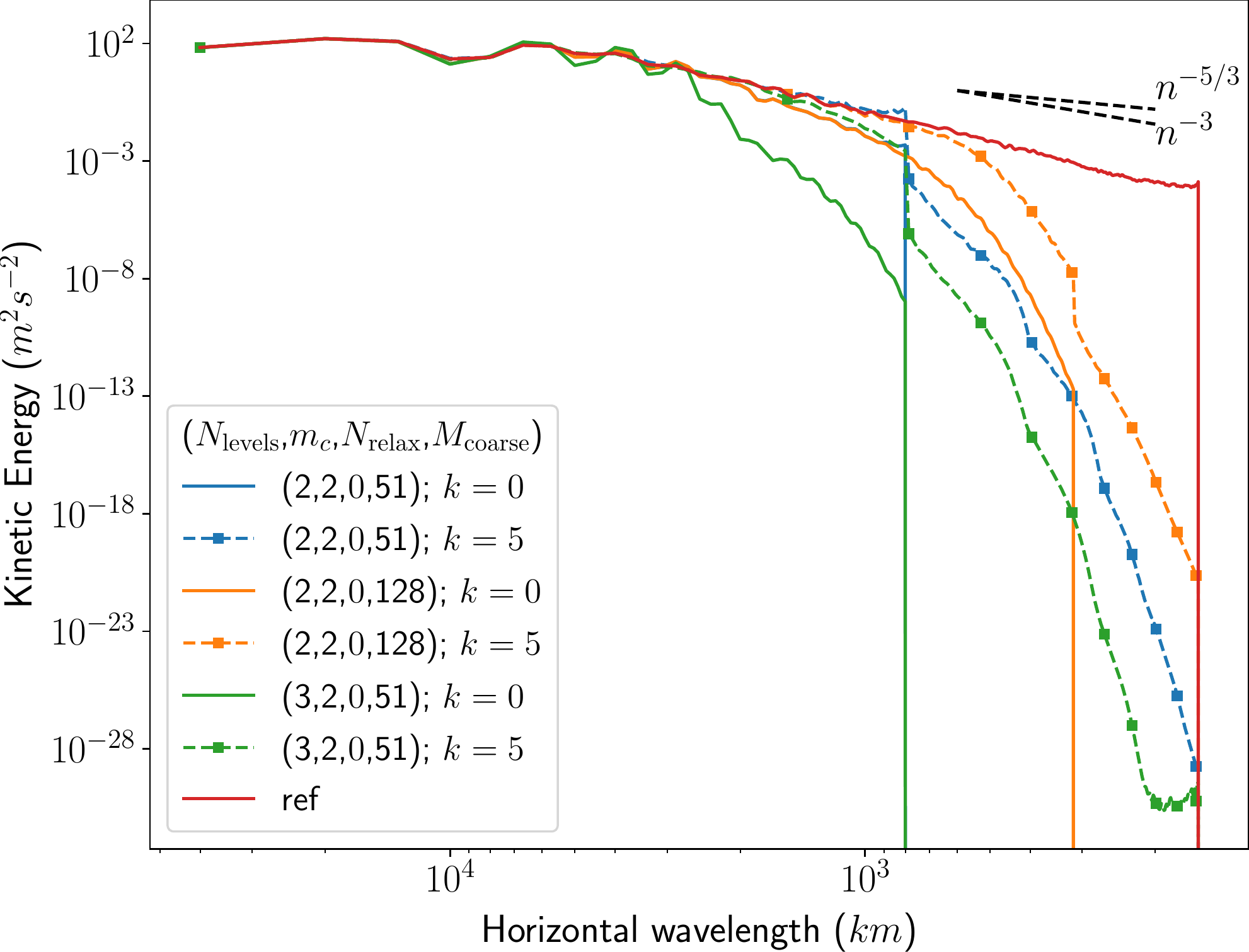}
        \caption{IMEX}
    \end{subfigure}
    \begin{subfigure}{.5\linewidth}
        \centering
        \includegraphics[scale=.35]{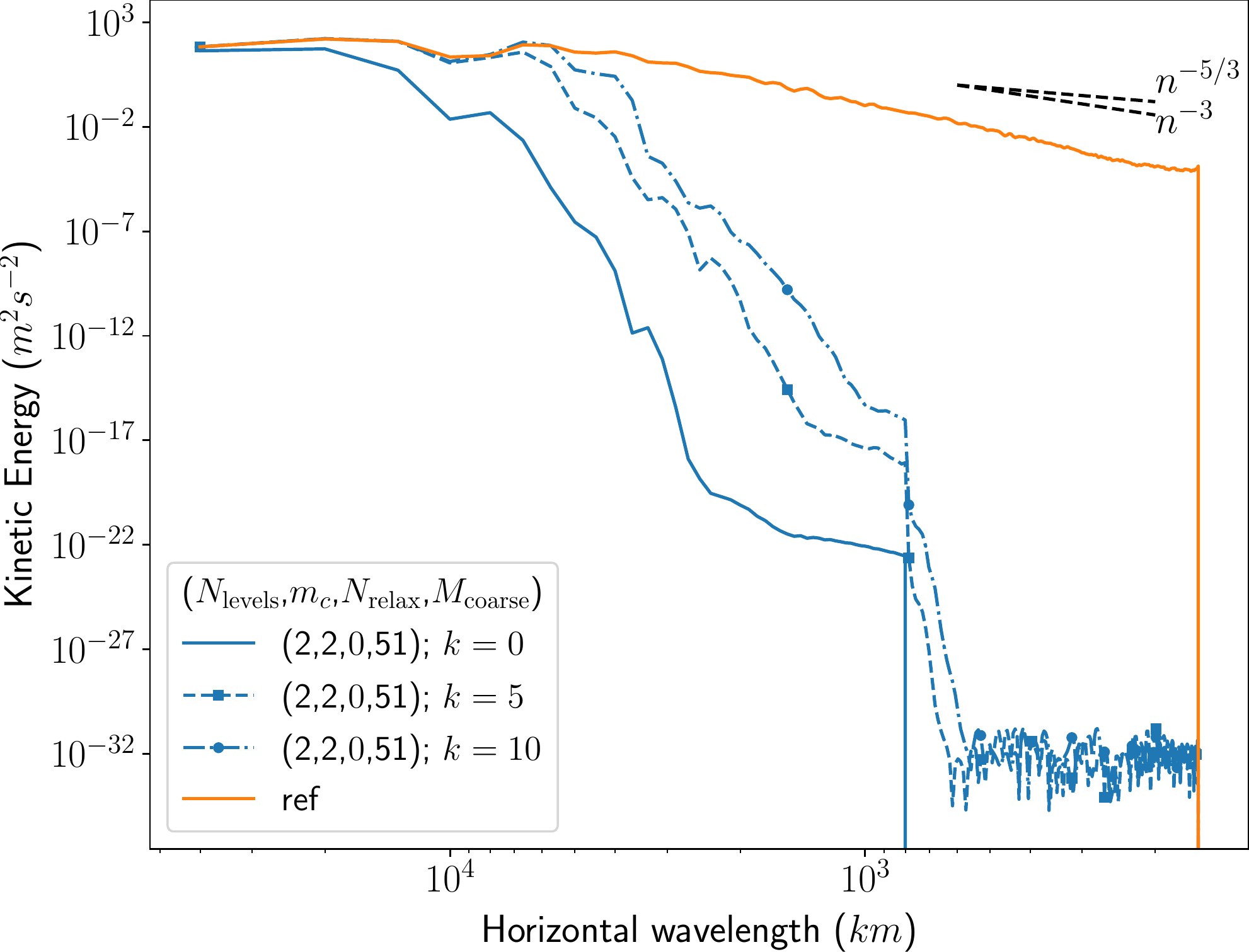}
        \caption{SL-SI-SETTLS}
    \end{subfigure}
    \caption{Unstable jet test case: kinetic energy spectra at $t = T$ of the reference (fine) and MGRIT solutions at given iterations $k$ for the configurations illustrated in Figure \ref{fig:unstable het_solution_params_mgrit}. Left and right: IMEX and SL-SI-SETTLS as coarse time-stepping schemes.} 
    \label{fig:unstable_jet_spectrum}
\end{figure}

\begin{figure}[!htbp]
    \begin{subfigure}{.5\linewidth}
        \centering
        \includegraphics[scale=.375]{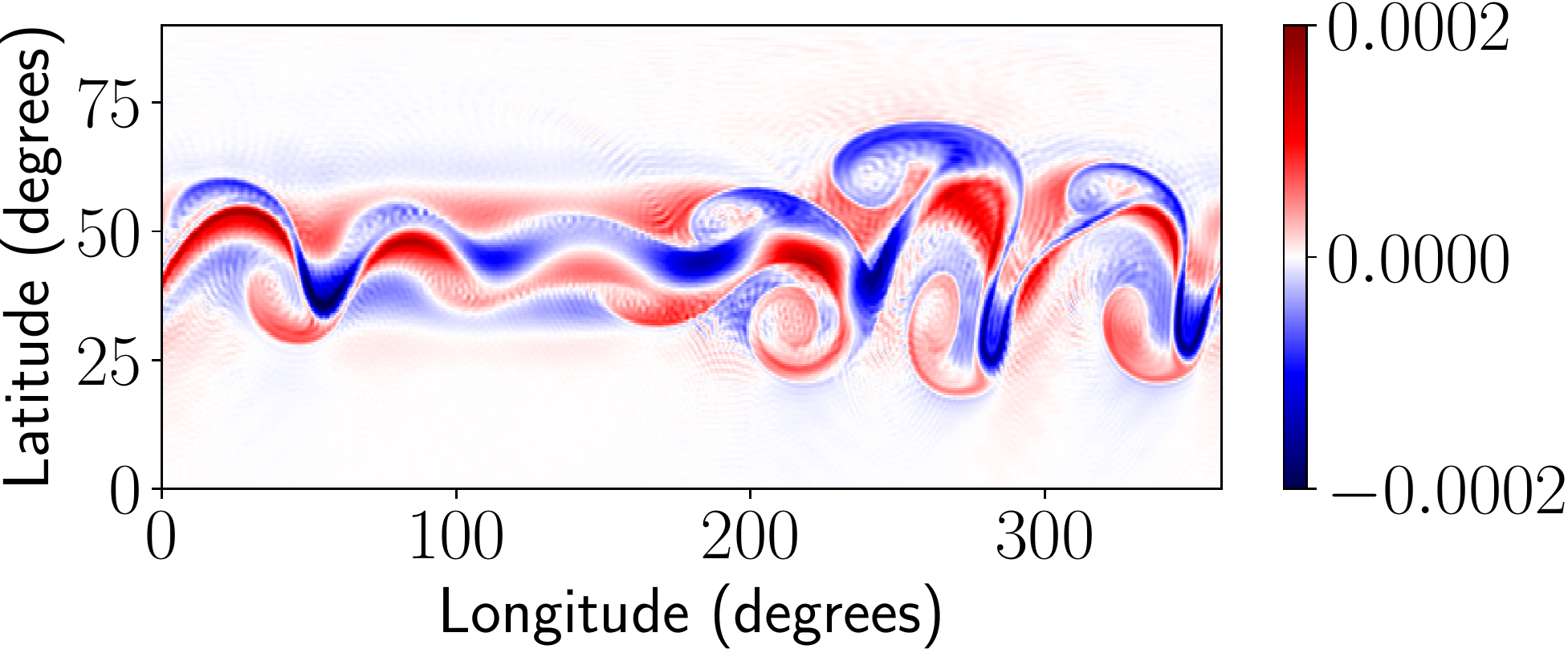}
        \caption{$(2,2,51,\text{IMEX})$; Iteration 0}
    \end{subfigure}
    \begin{subfigure}{.5\linewidth}
        \centering
        \includegraphics[scale=.375]{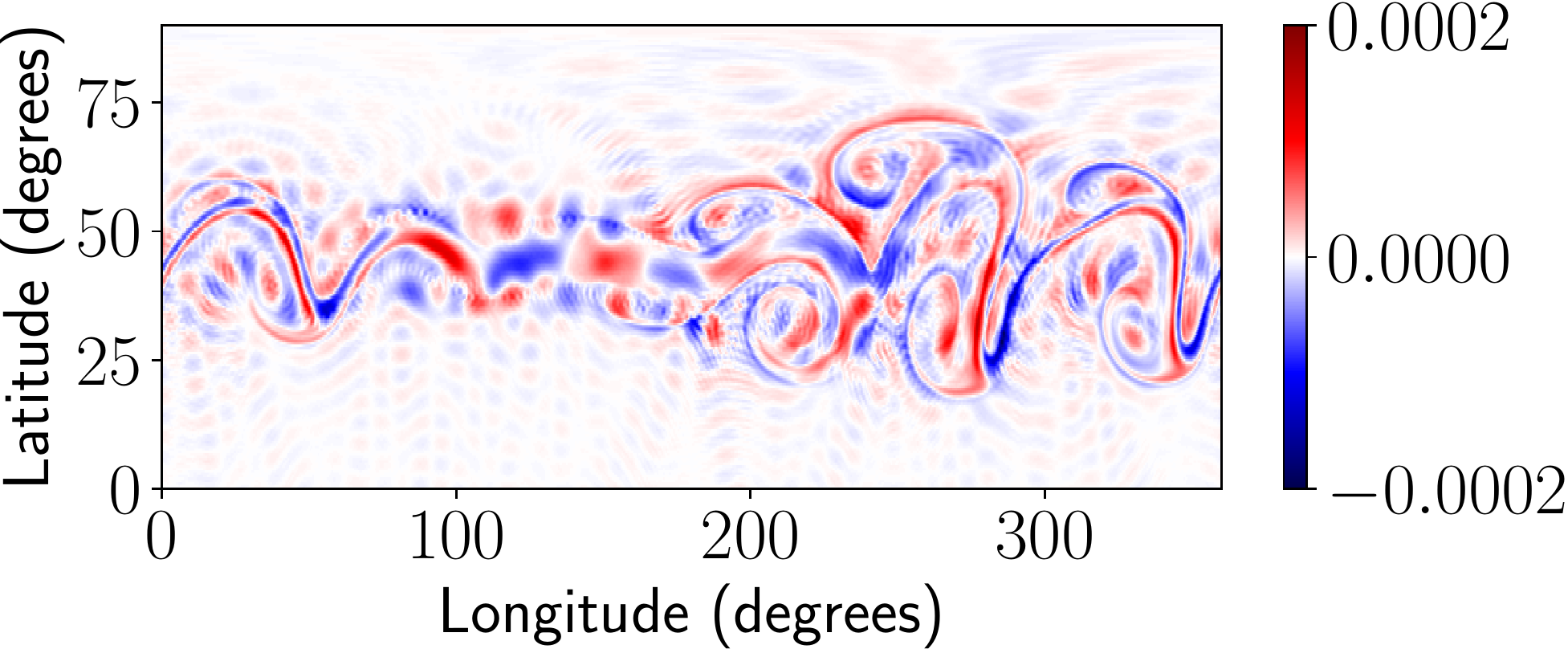}
        \caption{$(2,2,51,\text{IMEX})$; Iteration 5}
    \end{subfigure}
    \begin{subfigure}{.5\linewidth}
        \centering
        \includegraphics[scale=.375]{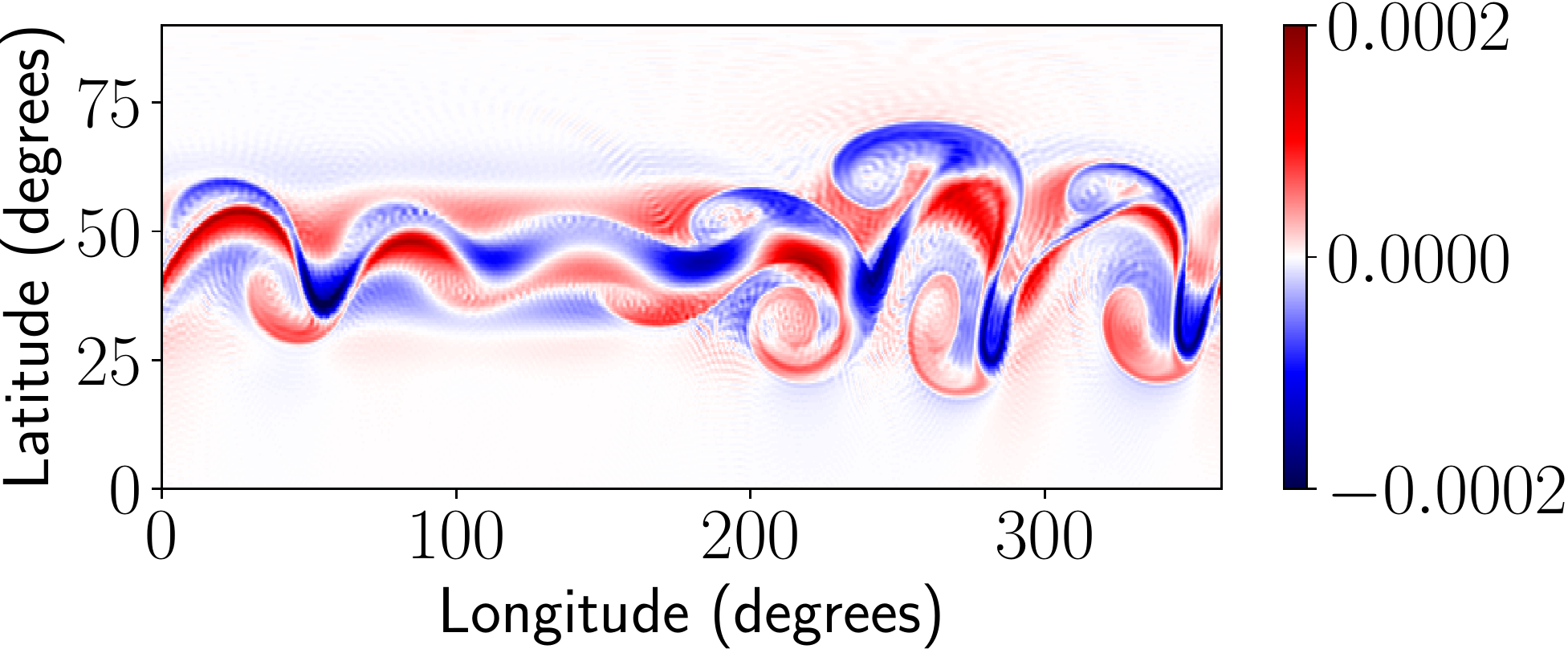}
        \caption{$(2,2,128,\text{IMEX})$; Iteration 0}
    \end{subfigure}
    \begin{subfigure}{.5\linewidth}
        \centering
        \includegraphics[scale=.375]{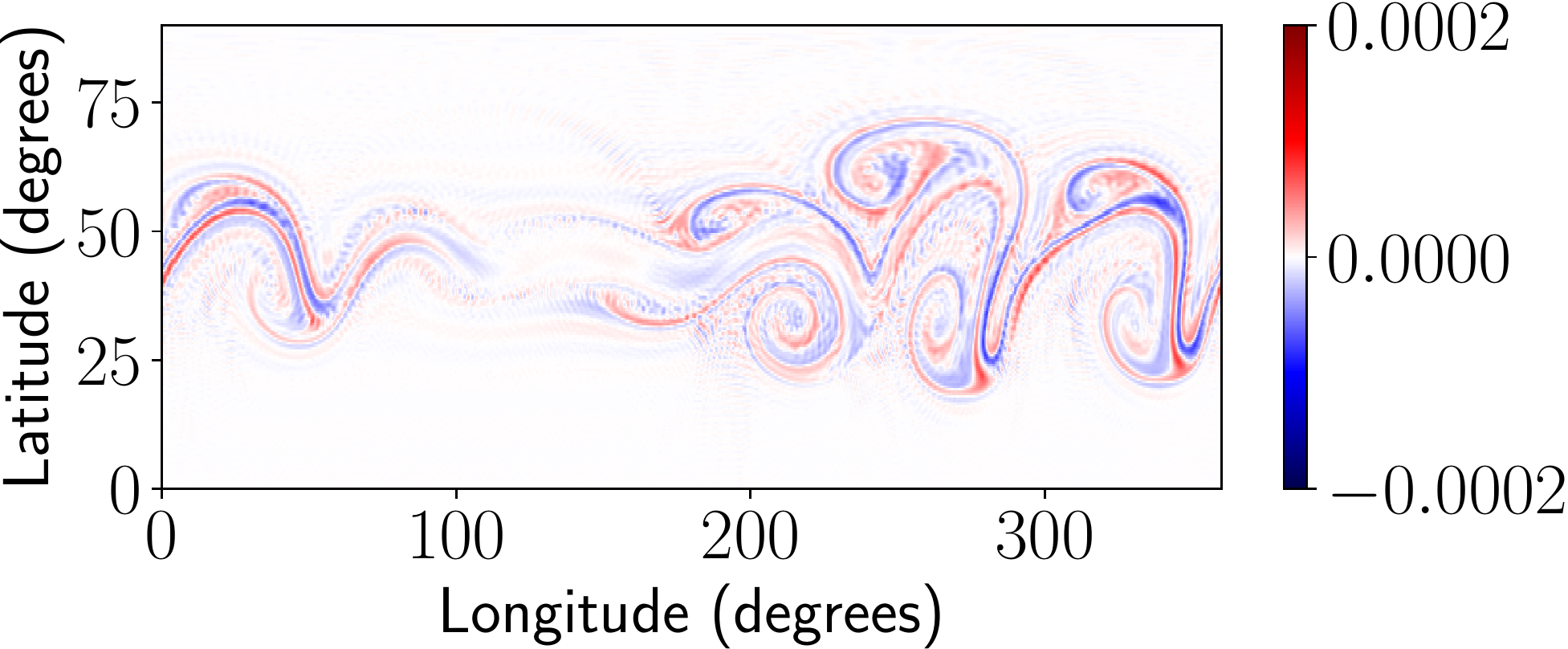}
        \caption{$(2,2,128,\text{IMEX})$; Iteration 5}
    \end{subfigure}
    \begin{subfigure}{.5\linewidth}
        \centering
        \includegraphics[scale=.375]{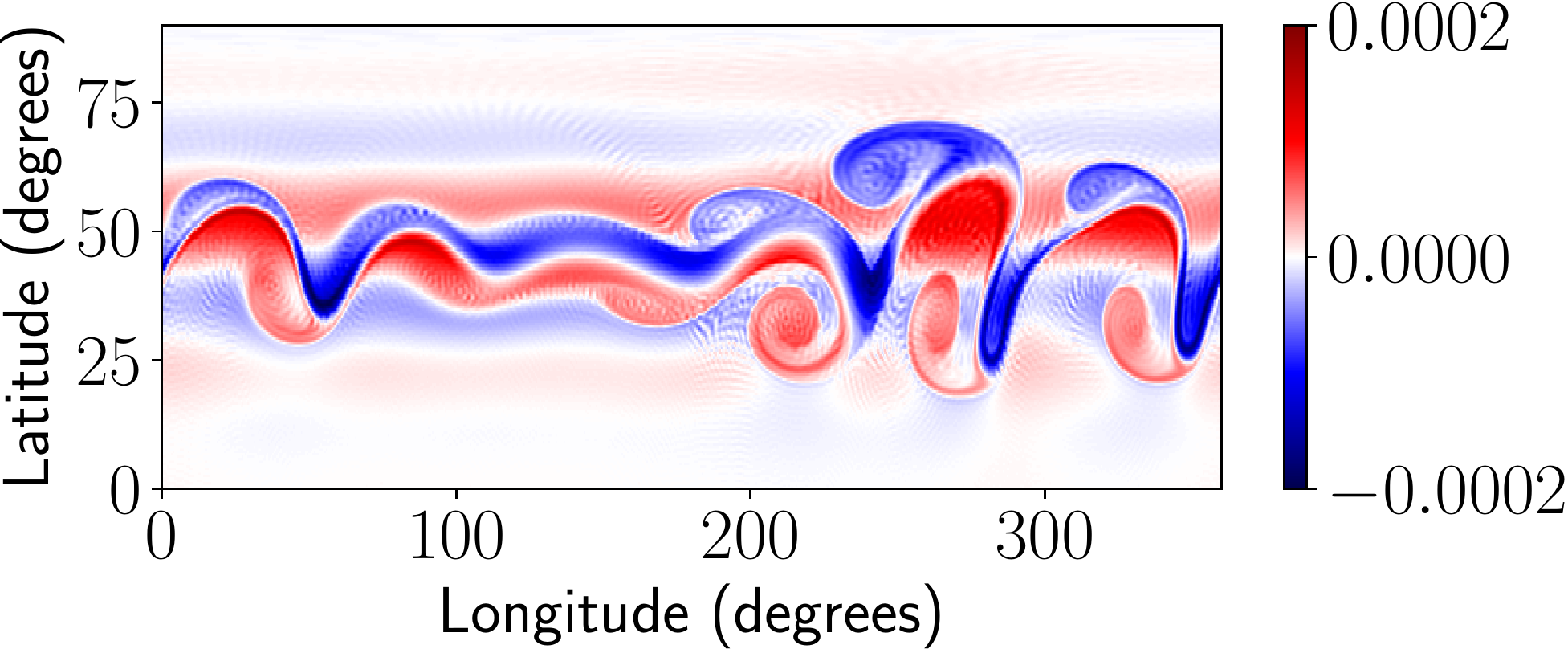}
        \caption{$(3,2,51,\text{IMEX})$; Iteration 0}
    \end{subfigure}
    \begin{subfigure}{.5\linewidth}
        \centering
        \includegraphics[scale=.375]{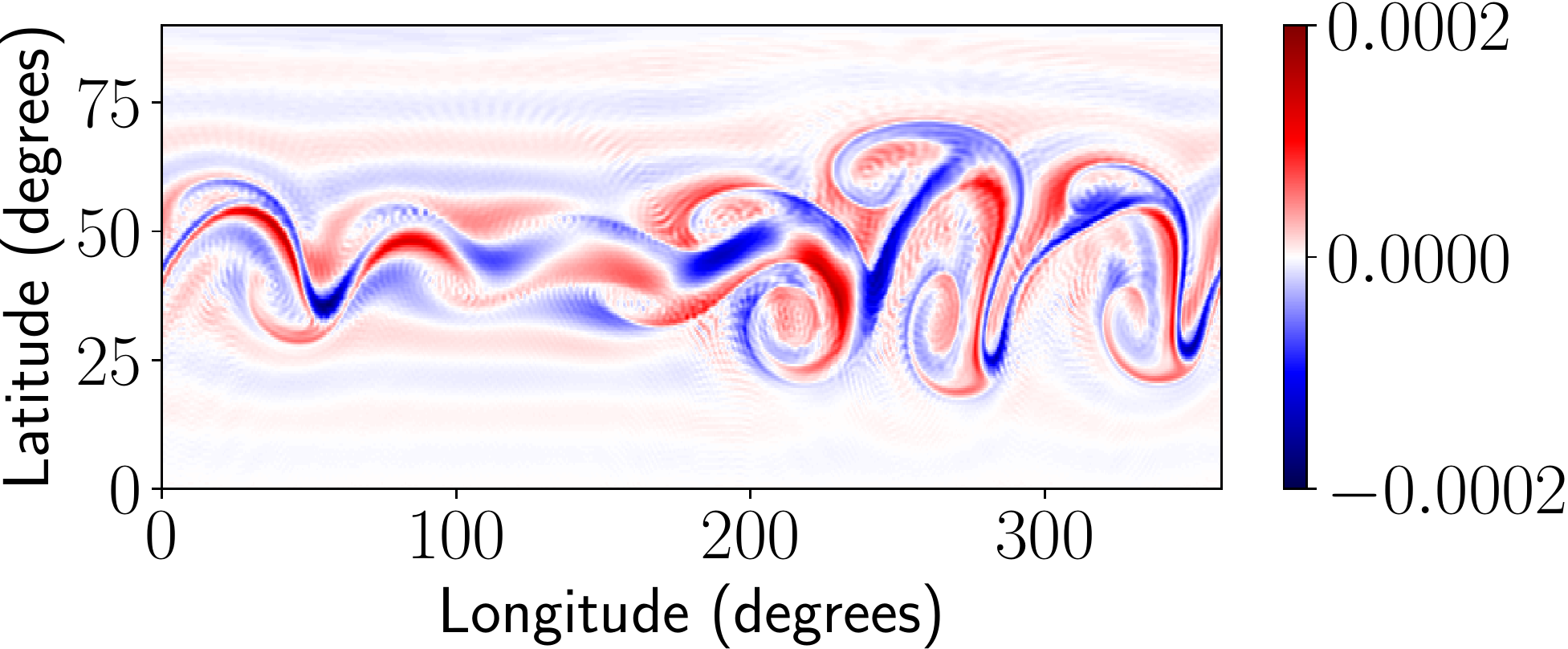}
        \caption{$(3,2,51,\text{IMEX})$; Iteration 5}
    \end{subfigure}
    \begin{subfigure}{.5\linewidth}
        \centering
        \includegraphics[scale=.375]{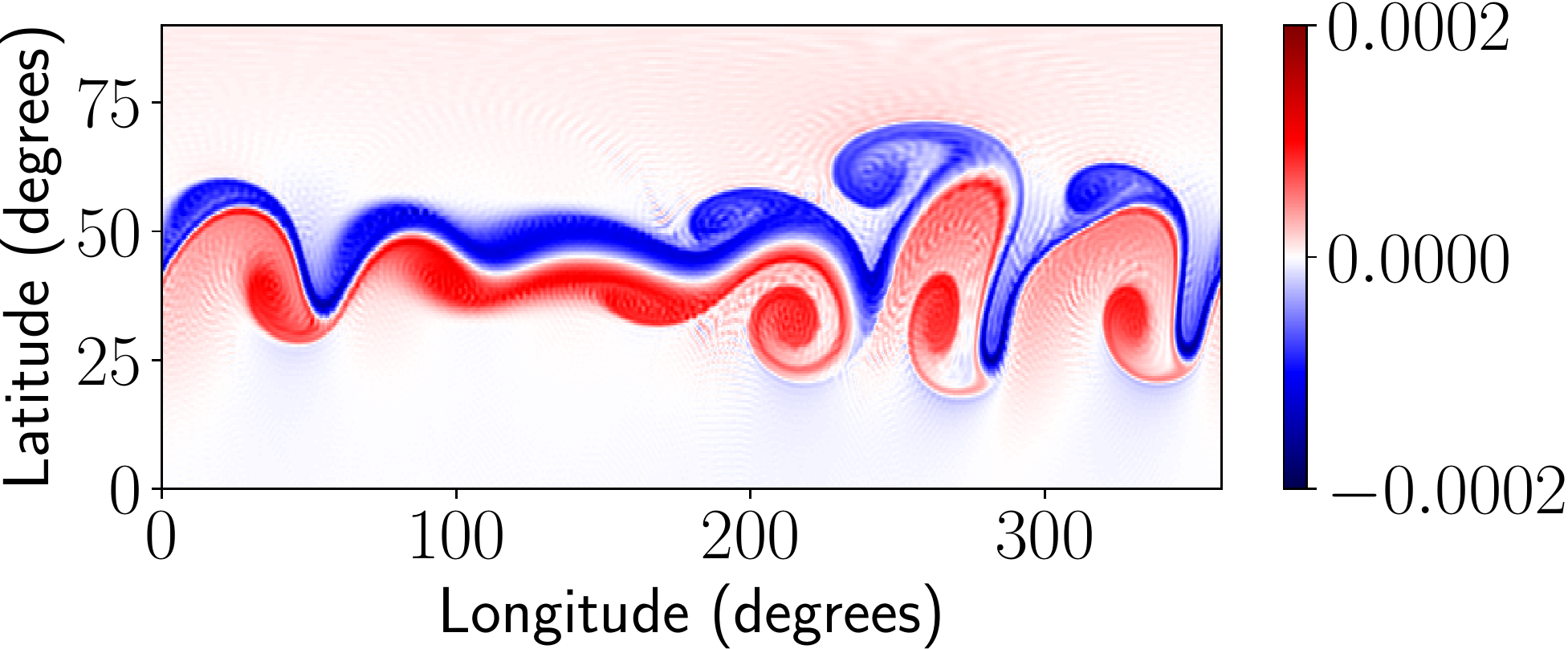}
        \caption{$(2,2,51,\text{SL-SI-SETTLS})$; Iteration 0}
    \end{subfigure}
    \begin{subfigure}{.5\linewidth}
        \centering
        \includegraphics[scale=.375]{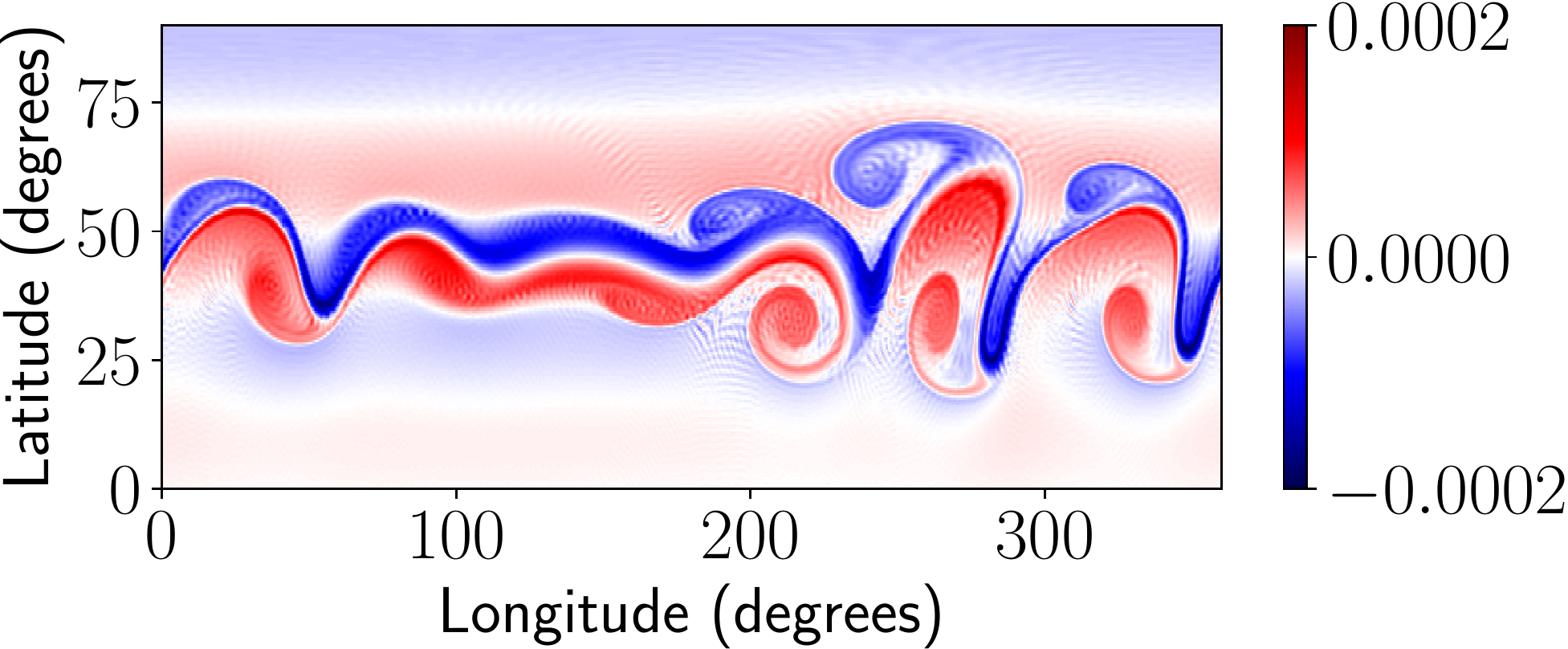}
        \caption{$(2,2,51,\text{SL-SI-SETTLS})$; Iteration 10}
    \end{subfigure}
    \caption{Unstable jet test case: the absolute difference between the PinT and the fine vorticity fields at $t = T$ under chosen configurations identified by $(\nlevels,\cfactor,\Mcoarse, \text{coarse time-stepping scheme})$ at iterations 0 (left) and 5 or 10 (right). All simulations use $\nrelax = 0$ and the viscosity configurations depicted in Table \ref{tab:unstable_jet_viscosity_coefficients}.}
    \label{fig:unstable het_solution_params_mgrit}
\end{figure}

\indent Figure \ref{fig:unstable het_solution_params_mgrit} presents the absolute difference between the PinT and reference vorticity fields at iterations $k = 0$ and $k = 5$ or $k = 10$ for chosen configurations, namely $(\nlevels, \cfactor, \nrelax, \Mcoarse) \in \{(2,2,0,51), \\(2,2,0,128), (3,2,0,51)\}$ and $(\nlevels, \cfactor, \nrelax, \Mcoarse) = (2,2,0,51)$, respectively with IMEX and SL-SI-SETTLS as coarse scheme, all of them using the viscosity configurations depicted in Table \ref{tab:unstable_jet_viscosity_coefficients}. The evolution of the respective kinetic energy spectra is illustrated in Figure \ref{fig:unstable_jet_spectrum}.
Concerning the configurations $(\nlevels, \cfactor, \nrelax, \Mcoarse) = (2,2,0,51)$ and $(\nlevels, \cfactor, \nrelax, \Mcoarse) = (2,2,0,128)$ using IMEX, which use a moderate viscosity coefficient on the coarse level, the kinetic energy spectrum is relatively close to the reference one already in the initial iteration; after five iterations, a slight amplification of medium wavenumbers is observed in the case $\Mcoarse = 51$, which translates to small-scale oscillations in the physical plot, despite a global reduction of the error magnitudes, whereas a much more important error decrease, without amplification of the spectrum, is obtained with $\Mcoarse = 128$. In the configuration $(\nlevels, \cfactor, \nrelax, \Mcoarse) = (3,2,0,51)$ using IMEX, the larger fourth-order viscosity coefficient applied on the coarsest level leads to a larger damping of the spectrum at iteration $k=0$, including small wavenumbers. It is clearly seen as large-scale oscillations in the respective physical plot. After five iterations, only smaller-scale errors are observed, but the spectrum is still outperformed by the initial one of the simulations using a smaller viscosity coefficient. Finally, in the simulation using SL-SI-SETTLS, in which a very aggressive second-order viscosity is applied, drastic damping is observed along the entire spectrum in the initial iteration, leading to very large-scale physical errors. The spectrum converges slowly: after ten iterations, it remains highly damped, and large-scale errors are still observed in the physical plot.

\subsubsection{Evaluation of computing times and speedups}

\indent We now proceed to the evaluation of the computing times and speedups provided by the PinT simulations in the unstable jet test case, whose reference solution takes a computing time $T_{\text{ref}} \approx 546 \text{s}$. As in the Gaussian bumps test case, we consider $\nproc = 64$ in all simulations.

\begin{figure}[!htbp]
    \begin{subfigure}{.5\linewidth}
        \centering
        \includegraphics[scale=.5]{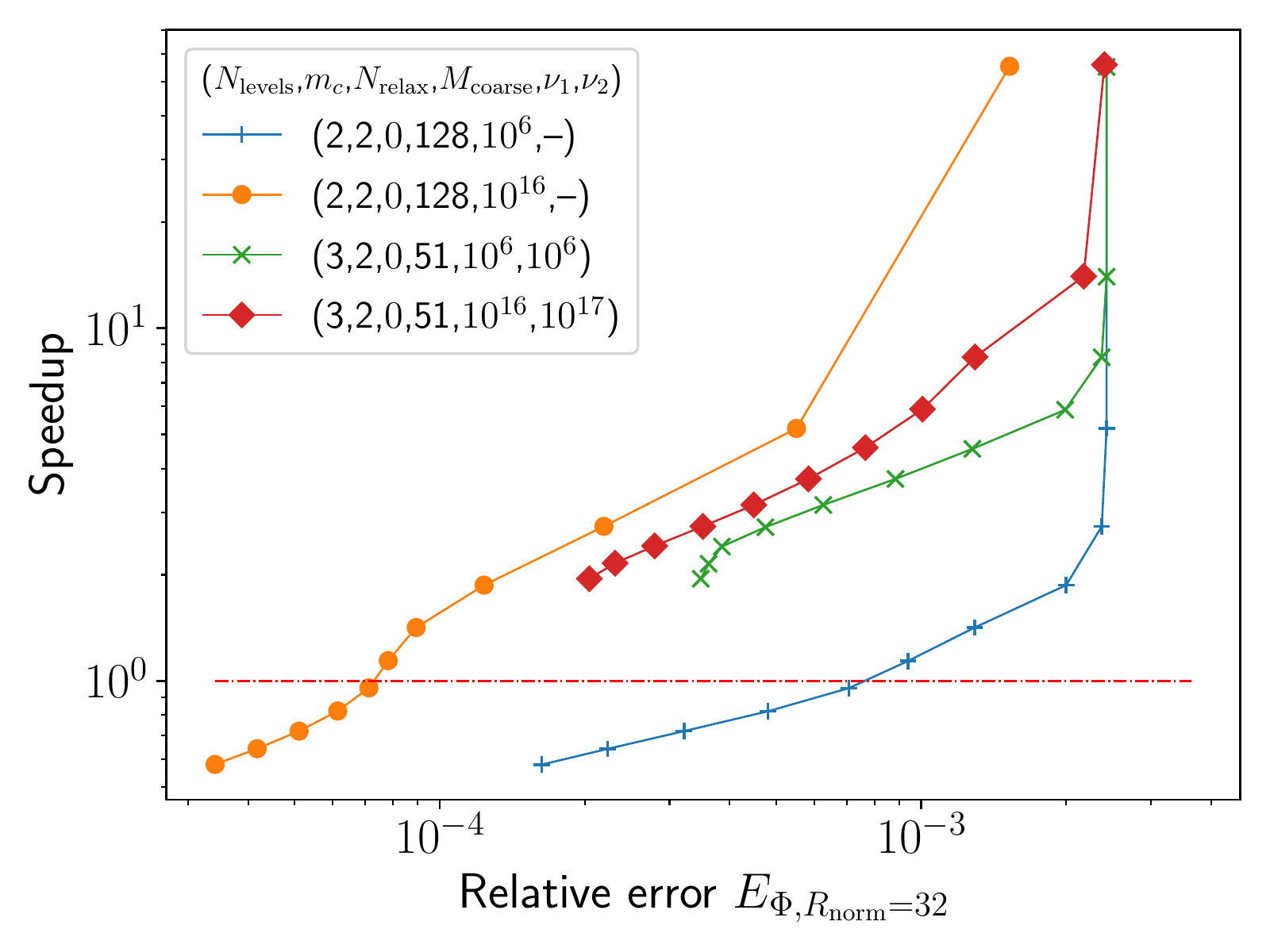}
        \caption{IMEX\label{fig:unstable_jet_speedup_error_IMEX}}
    \end{subfigure}
    \begin{subfigure}{.5\linewidth}
        \centering
        \includegraphics[scale=.5]{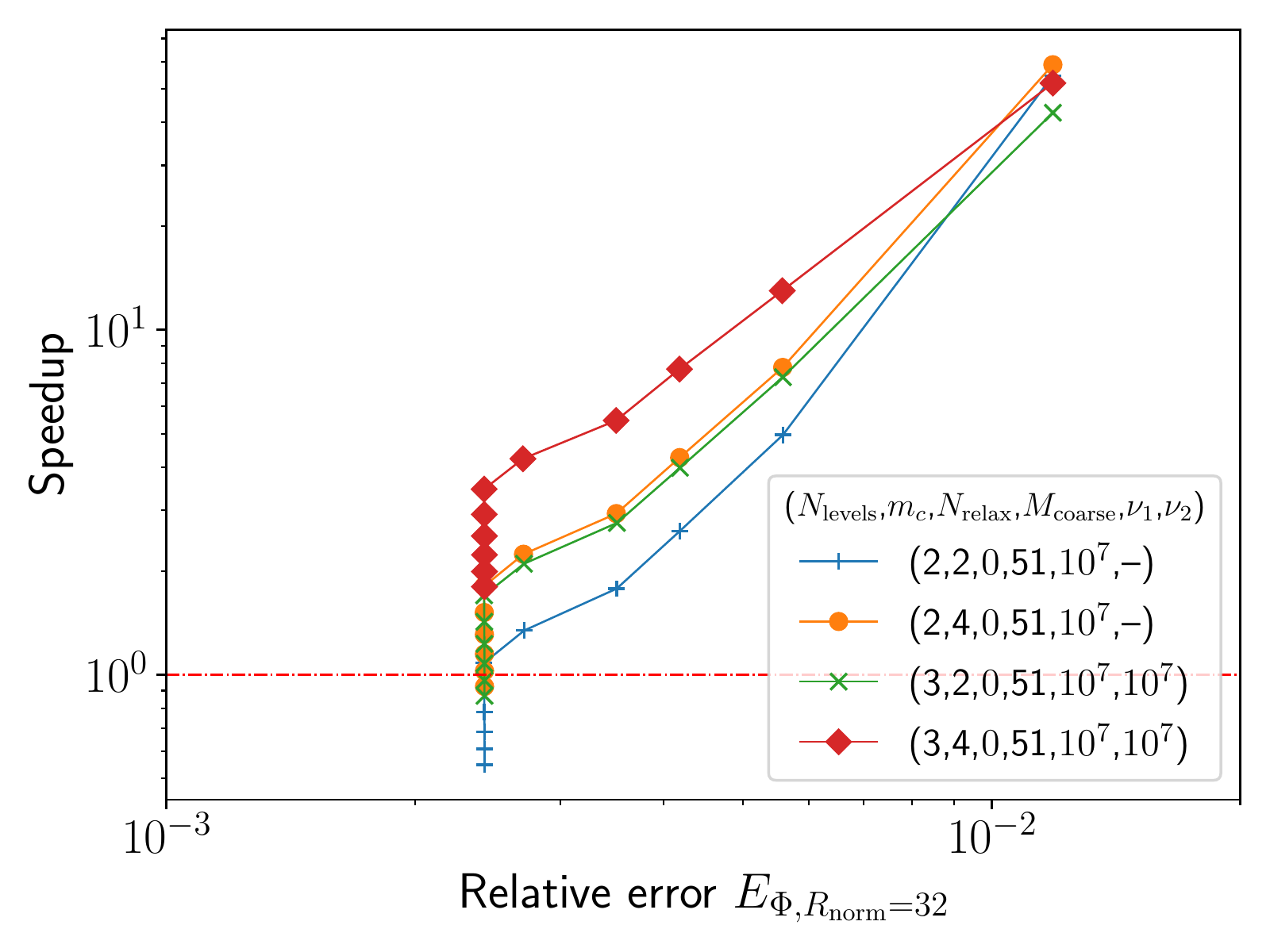}
        \caption{SL-SI-SETTLS\label{fig:unstable_jet_speedup_error_SL_SI_SETTLS}}
    \end{subfigure}
    \caption{Unstable jet test case: speedup with respect to the relative geopotential error in spectral space for chosen configurations using IMEX or SL-SI-SETTLS as coarse time-stepping scheme. Simulations are identified by $(\nlevels, \cfactor, \nrelax, \Mcoarse, \nu_1, \nu_2)$ and each data point corresponds to an iteration. The horizontal, dashed-dotted line indicates a unitary speedup. All simulations use $\nproc = 64$ parallel processors in time.}
    \label{fig:unstable_jet_speedup_error}
\end{figure}

\indent Figure \ref{fig:unstable_jet_speedup_error} presents the speedups as a function of the relative error on the geopotential field for some chosen configurations using IMEX or SL-SI-SETTLS on the coarse levels. In the former case, we consider two pairs of configurations, namely $(\nlevels,\cfactor,\Mcoarse) = (2,2,128)$ and $(\nlevels,\cfactor,\Mcoarse) = (3,2,51)$; in each pair, one configuration uses a second-order viscosity and the other a fourth-order one. Contrary to the Gaussian bumps simulations, we now observe, in this more complex test case, characterized by more important interactions between wavenumber modes, clearer differences between the viscosity approaches in terms of compromise between convergence and numerical acceleration. We easily see that the simulations using the second-order approach provide less interesting results in terms of speedup since almost no convergence is observed in the first three iterations. Better results are obtained using the fourth-order viscosity, mainly the configuration $(\nlevels,\cfactor,\Mcoarse) = (2,2,128)$, despite its larger spectral resolution on the coarse levels. Indeed, two iterations provide nearly the same error as ten iterations of the configuration $(\nlevels,\cfactor,\Mcoarse) = (3,2,51)$ using fourth-order viscosity, with approximate speedup factors of 2.7 and 1.9, respectively. Finally, for these two configurations, the solutions whose error plots are depicted in Figure \ref{fig:unstable het_solution_params_mgrit}, corresponding to the iteration $k=5$, are obtained with respective speedups of 1.1 and 3.7.

\indent In the case where SL-SI-SETTLS is used as a coarse time-stepping scheme, we recall that all stable simulations present approximately the same convergence behavior, due to the overdamping induced by the large second-order viscosity approach. Therefore, the comparison of the speedups allows us to study the influence of the choice of Parareal and MGRIT parameters on the computational cost. In Figure \ref{fig:unstable_jet_speedup_error_SL_SI_SETTLS}, we only present simulations using the coarse spectral resolution $\Mcoarse = 51$ and the errors under $\rnorm = 32$. The largest and smallest speedups are provided respectively by $(\nlevels,\cfactor) = (3,4)$ and $(\nlevels,\cfactor) = (2,2)$, whose coarsest levels, on which the time integration is serial, use the largest and smallest time step sizes, with speedups of approximately 3.5 and 1.1 before the stagnation of the convergence. The simulations $(\nlevels,\cfactor) = (3,2)$ and $(\nlevels,\cfactor) = (2,4)$, whose coarsest levels use the same time step size, present intermediate and very similar speedup results.

\subsection{Discussion}

\indent \added[id = R3]{The study presented in this section aimed to evaluate the performance, in terms of stability, convergence and computational cost, of Parareal and MGRIT applied to the numerical integration of the SWE on the rotating sphere. More specifically, the main goal was to study the influence of the choice of time stepping scheme used on the coarse levels, which is known to be a crucial factor for the PinT performance, and how the results can be improved by proper parametric and discretization choices. The obtained results provide indeed indications on which characteristics a coarse scheme should have in order to lead to a successful temporal parallelization, \ie with stable and relatively accurate solutions being provided within smaller computational times than serial integrations.}

\indent \added[id=R3]{ From the numerical simulations, and also from the analytical stability study developed in Section \ref{subsec:stability}, it is clear that it is necessary but not sufficient for the coarse scheme to be stable: poor results are obtained when the temporal parallelization uses SL-SI-SETTLS on the coarse levels, despite the stability and popularity of this scheme for the serial integration of atmospheric models and recent works indicating that the use of semi-Lagrangian approaches on the coarse levels improves the performance of PinT methods applied to simpler problems.}

\indent \added[id=R3]{Moreover, the coarse discretization should be able to accurately represent large spatial scales and ensure stability on the fine ones. Indeed, the best results using IMEX as a coarse scheme are obtained with large higher-order viscosity coefficients on the coarse levels, damping the largest wavenumbers but preserving the intermediate and smallest ones. In the configurations using SL-SI-SETTLS on the coarse levels, the very restrictive stability constraints require the use of large second-order viscosities, which damps a large range of the spectrum, providing stability at the expense of a great accuracy loss.}

\indent \added[id=R3]{The detailed and level-dependent study conducted on the influence of the adopted viscosity approach provides other insightful conclusions. The tests were performed considering viscosity values around typical values adopted on atmospheric modeling, and it was found that relatively large values are required for ensuring stability in the PinT framework. However, it was verified that the artificial viscosity is required only on the coarsest levels, on which the stability constraints are more severe due to the large time step sizes adopted. It implies that the fine dynamics do not need to be damped on the fine discretization levels, and more accurate and faster converging PinT solutions can be expected.}

\indent \added[id=R3]{Finally, these conclusions imply that a compromise between accuracy and stability of the coarse discretization has to be found in order to achieve speedup using PinT methods. Our results indicate that, if the conditions of the coarse discretization described above are met, and by using enough parallel resources, Parareal and MGRIT using IMEX are able to provide relatively accurate solutions within shorter computational times compared to fine reference simulations.}

\section{Conclusion and perspectives}
\label{sec:conclusion}

\indent In this work, we have studied the temporal parallelization of the shallow water equations on the rotating sphere using Parareal and MGRIT. \replaced[id=R2]{The development and application of PinT methods to hyperbolic problems still advance relatively slowly due to the well-known stability and convergence issues when applied to simple problems such as the advection equation}{It is well-known that PinT methods present stability and convergence issues when applied to hyperbolic problems}, with the choice of coarse temporal discretization being a crucial aspect in this context. Therefore, the focus here was to investigate if popular and well-established time-stepping schemes in the atmospheric modeling community can provide good performance results when used on the coarse levels of the PinT methods. The two considered schemes, IMEX and SL-SI-SETTLS, allow the use of relatively large time steps in serial integrations since they avoid stability issues linked to the stiff terms of the governing equations.

\indent Two approaches were considered in this study. First, we conducted an analytical stability investigation of the PinT methods applied to a linearized ODE. A notable result is that Parareal and MGRIT using SL-SI-SETTLS as a coarse scheme present very poor stability properties compared to IMEX. Second, we performed numerical simulations of two test cases with increasing complexity.
This confirmed the poor stability using SL-SI-SETTLS, with only very restrictive PinT configurations presenting a stable behavior. Better results were obtained with IMEX, but still relatively limited in terms of choice of several levels and coarsening factors in time. A further investigation indicated that better choices of artificial viscosity parameters only on the coarse levels allow for improved stability and convergence by using second- and higher-order viscosities for SL-SI-SETTLS and IMEX. In the former case, however, the stability is obtained at the cost of a strong damping of the entire wavenumber spectrum, which reduces the accuracy provided by the temporal parallelization. To the best of our knowledge, this
is the first time that a level-selective viscosity has been
investigated in PinT methods. Finally, studies on parallel performance indicated that the best compromises between accuracy and computational cost, with speedups larger than one, are obtained using IMEX as a coarse scheme, mainly due to its better stability and convergence properties in the PinT framework.

\indent In summary, the study presented in this work indicates that inferior results are obtained when the temporal parallelization uses SL-SI-SETTLS on the coarse levels, despite the stability and popularity of this scheme for the serial integration of atmospheric model and recent works indicating that the use of semi-Lagrangian approaches on the coarse levels improves the performance of PinT methods. Better results are obtained using IMEX, but depending on the discretization applied on each level and, consequently, on appropriate choices of viscosity and hyperviscosity parameters, whose influence is especially remarkable in the more complex unstable jet test case. Under these conditions, and by using enough parallel resources, Parareal and MGRIT using IMEX can provide relatively accurate solutions within shorter computational times compared to fine reference simulations.

\indent \added[id=R3]{An important and open challenge concerns the gap to be filled between the results obtained in this work and the application of parallel-in-time methods for real weather and climate problems. Through parametrization and discretization choices coherent with practical applications, we tried, to the extent possible, to provide insights in this direction, \eg with the parameters considered in the analytical stability study, the spectral resolutions adopted, the ranges of tested viscosity and hyperviscosity coefficients and the simulation of quite challenging test cases, mainly the unstable jet one, which is a standard test in atmospheric modeling due to its complex dynamics. However, more refined investigations would need to be conducted to reach practical applications in operational models. }

\indent A natural future work consists of investigating if modifications of SL-SI-SETTLS and the use of other time-stepping schemes, eventually not still used operationally in atmospheric models or even similar schemes but with higher discretization orders, could provide better stability and convergence properties for Parareal and MGRIT. For instance, exponential integration methods, which can integrate precisely the linear terms of the governing equations, have aroused a growing interest in the context of atmospheric modeling, \eg in \cite{gaudreault_pudykiewicz:2016, schreiber_al:2019, gaudreault_al:2022}. It includes a semi-Lagrangian variant of  this family of method, which have been proposed and used for solving the SWE on the plane by \cite{peixoto_schreiber:2019}, and to the SWE on the rotating sphere with a cubed sphere spatial discretization by \cite{shashkin_goyman:2020}. Ongoing studies indicate that the use of exponential schemes \added[id=R2]{and their semi-Lagrangian versions} on the coarse discretization levels indeed provides better stability and convergence properties for Parareal and MGRIT, which will be presented in the following article.

\indent Finally, the results presented in this work illustrate different behaviors of the wavenumber spectra of the solution along the PinT iterations, with the adopted artificial viscosity approach being crucial for ensuring stability and convergence\added[id=R3]{, and notably being required only on the coarse discretization levels. }Since the SWEs are a nonlinear model, with nonlinear interactions between modes transferring energy to the highest wavenumbers, future works should focus on understanding how each time-stepping scheme treats these interactions and how they influence the performance of the temporal parallelization.

\section*{CRediT authorship contribution statement}
\indent \textbf{João G. Caldas Steinstraesser:} Conceptualization, Methodology, Software, Validation, Formal analysis, Investigation, Writing - Original Draft, Visualization.
\textbf{Pedro da Silva Peixoto:} Conceptualization, Methodology, Resources, Writing - Review \& Editing, Supervision, Project
administration, Funding acquisition.
\textbf{Martin Schreiber:} Conceptualization, Methodology, Software, Validation, Formal analysis,
Resources, Writing - Review \& Editing, Supervision, Project
administration, Funding acquisition.

\section*{Declaration of competing interest}

\indent The authors declare that they have no known competing financial interests or personal relationships that could have appeared to influence the work reported in this paper.

\section*{Data availability}

% \indent The source code and numerical tests presented here are available for reproducibility in the repository \cite{schreiber:2016}.

\indent The source code and numerical tests presented here are available for reproducibility in the repository \begin{tiny}\url{https://gitlab.inria.fr/sweet/sweet/-/tree/58a4c051fe1e928573284d9d94aaf571104d9e05/benchmarks_sphere/paper_jcp_pint_imex_slsi}\end{tiny}

\section*{Acknowledgements}

\indent \textbf{Funding:} This work was supported by the São Paulo Research Foundation (FAPESP) grants 2021/03777-2 and 2021/06176-0, as well as the Brazilian National Council for Scientific and Technological Development (CNPq), Grant 303436/2022-0. This project also received funding from the Federal Ministry of Education and Research and the European High-Performance Computing Joint Undertaking (JU) under grant agreement No 955701, Time-X. The JU receives support from the European Union’s Horizon 2020 research and innovation programme and Belgium, France, Germany, Switzerland.

\indent Most of the computations presented in this paper were performed using the GRICAD infrastructure (\url{https://gricad.univ-grenoble-alpes.fr}), which is supported by Grenoble research communities.

\appendix
\section{Viscosity and hyperviscosity in spectral methods}
\label{app:viscosity}

\indent We make a brief overview of the theory of viscosity approaches in spectral methods and some considerations about the orders of magnitude to be chosen for the viscosity coefficients, following the presentation in \cite{jablonowski_williamson:2011}, to which we refer the reader for further details.

\indent As mentioned in Section \ref{sec:swe}, the viscosity approach with even order $q \geq 2$ applied to a field $\psi$ consists in

\begin{equation}
    \label{eq:viscosity}
    \pdert{\psi} = (-1)^{\frac{q}{2} + 1}\nu \nabla^q \psi
\end{equation}

\noindent where the $(-1)^{\frac{q}{2}+1}$ allows to define $\nu \geq 0$. In the spectral space of spherical harmonics, \eqref{eq:viscosity} reads

\begin{equation*}
    \pdert{\psi_{m,n}} = (-1)^{\frac{q}{2}+1} \nu \left( \frac{-n(n+1)}{a^2} \right)^{\frac{q}{2}}  \psi_{m,n} = -\nu \left( \frac{n(n+1)}{a^2} \right)^{\frac{q}{2}}  \psi_{m,n}
\end{equation*}

\noindent whose exact solution is

\begin{equation*}
    \psi_{m,n}(t) = \psi_{m,n}(0) \exp{ \left(-\nu \left( \frac{n(n+1)}{a^2} \right)^{\frac{q}{2}} t \right) } 
\end{equation*}

\indent The viscosity coefficient can be determined by setting a time $\tau$ in which a given wavenumber $n_0$ is damped to a fraction $b_{\tau, n_0}$ of $\psi_{m,n}(0)$. Typically, $\nu$ is set such that the largest wavenumber $n_0 = M$ damps to a fraction $b_{\tau, n_0} = 1/e$, \ie

\begin{equation}
    \label{eq:viscosity_coefficient}
    \nu = \frac{1}{\tau} \left( \frac{M(M+1)}{a^2} \right)^{-\frac{q}{2}}
\end{equation}

\indent With this viscosity parameter being applied to the entire spectrum, a given mode $\psi_{m,n}$ is damped after a time step $\Dt$ by a factor

\begin{equation*}
    \begin{aligned}
    b_{\Dt,n} & := \exp{\left( - \nu \left( \frac{n(n+1)}{a^2} \right)^{\frac{q}{2}} \Dt   \right)} = \exp{  \left( -\frac{\Dt}{\tau} \left( \frac{n(n+1)}{M(M+1)} \right)^{\frac{q}{2}} \right) }\\
    & \approx \left[ 1 +  \Dt \nu \left( \frac{n(n+1)}{a^2} \right)^{\frac{q}{2}} \right]^{-1} =: \hat{b}_{\Dt,n}
    \end{aligned}
\end{equation*}

\noindent where $\hat{b}_{\Dt, n}$ corresponds to a backward Euler discretization of \eqref{eq:viscosity}, as considered in this work.

\begin{figure}[!htbp]
    \begin{subfigure}{.5\linewidth}
        \centering
        \includegraphics[scale=0.45]{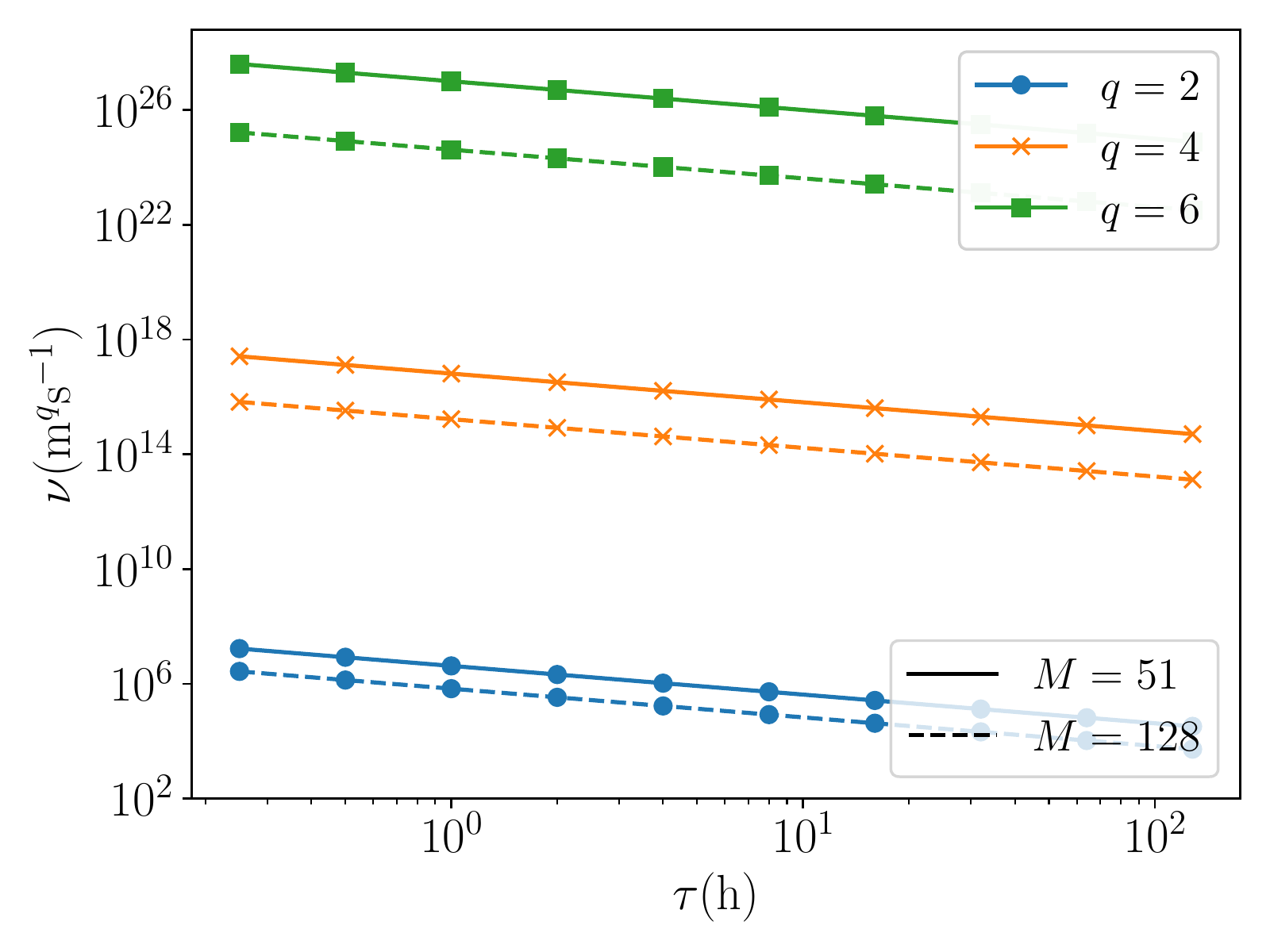}
        \caption{\label{fig:viscosity_time}}
    \end{subfigure}
    \begin{subfigure}{.5\linewidth}
        \centering
        \includegraphics[scale=0.45]{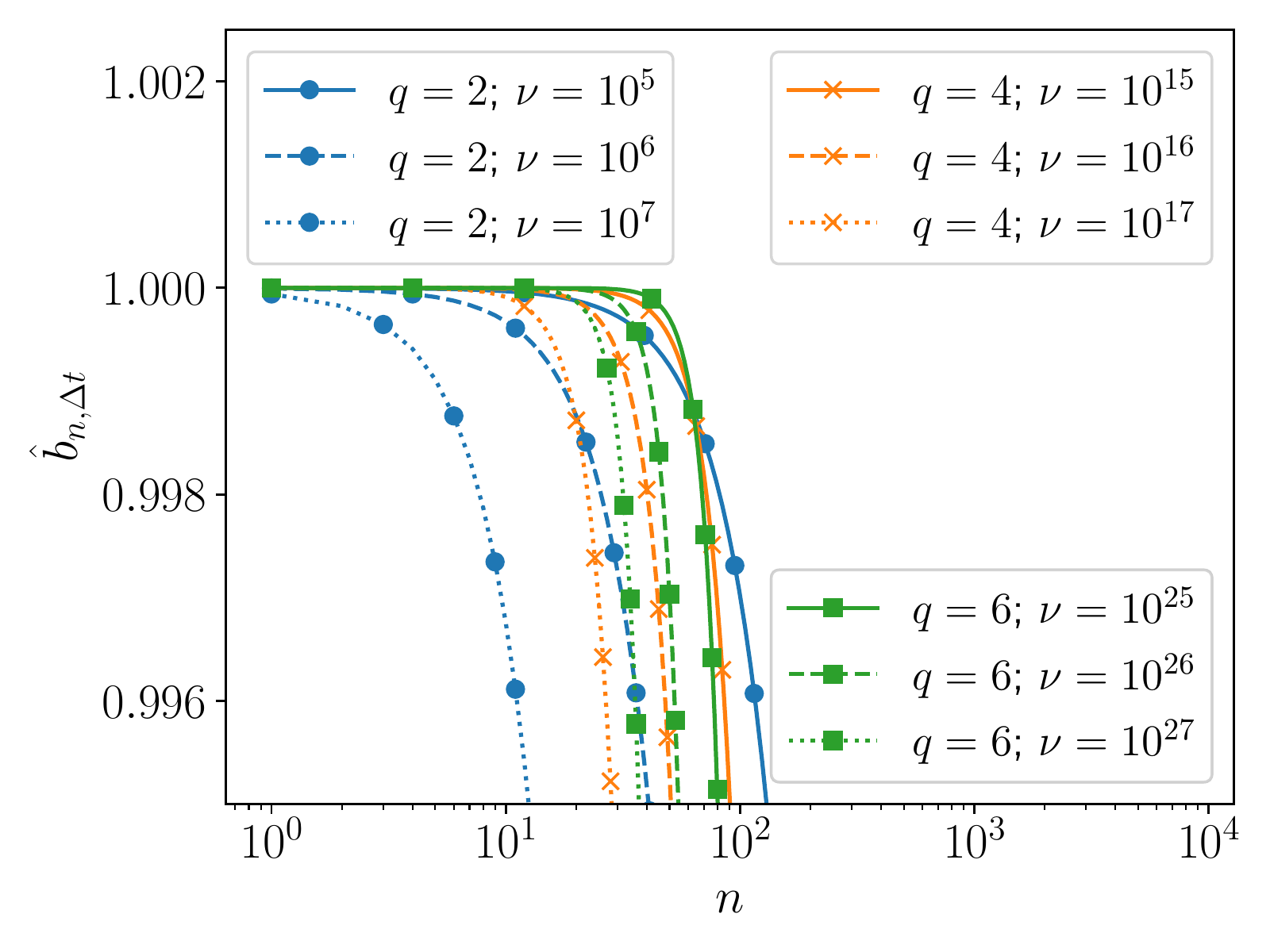}
        \caption{\label{fig:viscosity_damping_factor}}
    \end{subfigure}
    \caption{(a): viscosity coefficient of orders $q \in \{2,4,6\}$ as a function of the damping time for chosen spectral resolutions $M$. Computed using \eqref{eq:viscosity_coefficient}. (b): discrete damping factor $\hat{b}_{\Dt,n}$ as a function of the spectral mode for viscosity orders 2, 4 and 6 and chosen viscosity coefficients (in $\text{m}^q\text{s}^{-1}$). Computed considering $\Dt = 120 \text{s}$.}
    \label{fig:viscosity}
\end{figure}

\indent Figure \ref{fig:viscosity_time} illustrates the orders of magnitude of the viscosity coefficient values for viscosity orders $q \in \{2,4,6\}$ as a function of the damping time $\tau$ and for the two spectral resolutions considered on the coarse discretization levels of the PinT numerical simulations performed here, namely $M \in \{51, 128\}$. Damping timescales reported and suggested in the literature for spectral models usually range around units or tenths of hours, depending on the spectral resolution (see \eg \cite{boville:1991,hamilton_al:2008,williamson_al:2008}). In Figure \ref{fig:viscosity_damping_factor} we plot the discrete damping factor $\hat{b}_{\Dt, n}$ \wrt the wavenumber modes $n$ for some viscosity coefficients around the approximate average values depicted in Figure \ref{fig:viscosity_time} and time step $\Dt = 120\text{s}$. We observe that smaller-order viscosities lead to a faster decay of smaller wavenumbers and a larger sensitivity on the coefficient value (compare the spacing between the curves corresponding to the same viscosity order). It indicates that a higher-order viscosity produces a more refined stability filter, damping only the largest wavenumbers and better preserving the large-scale spatial features of the solution.

\bibliography{biblio}

\begin{thebibliography}{10}

\bibitem{xbraid-package}
{XB}raid: Parallel multigrid in time.
\newblock \url{http://llnl.gov/casc/xbraid}.

\bibitem{abel_al:2020}
N.~Abel, J.~Chaudhry, R.~D. Falgout, and J.~Schroder.
\newblock Multigrid-reduction-in-time for the rotating shallow water equations,
  8 2020.

\bibitem{astorino_al:2012}
Matteo Astorino, Franz Chouly, and Alfio Quarteroni.
\newblock {Multiscale coupling of finite element and lattice Boltzmann methods
  for time dependent problems}.
\newblock Technical report, October 2012.

\bibitem{baffico_al:2002}
L.~Baffico, S.~Bernard, Y.~Maday, G.~Turinici, and G.~Z{\'{e}}rah.
\newblock Parallel-in-time molecular-dynamics simulations.
\newblock {\em Physical Review E}, 66(5), November 2002.

\bibitem{bal:2005}
Guillaume Bal.
\newblock On the convergence and the stability of the parareal algorithm to
  solve partial differential equations.
\newblock In Timothy~J. Barth, Michael Griebel, David~E. Keyes, Risto~M.
  Nieminen, Dirk Roose, Tamar Schlick, Ralf Kornhuber, Ronald Hoppe, Jacques
  P{\'e}riaux, Olivier Pironneau, Olof Widlund, and Jinchao Xu, editors, {\em
  Domain Decomposition Methods in Science and Engineering}, pages 425--432,
  Berlin, Heidelberg, 2005. Springer Berlin Heidelberg.

\bibitem{bal_maday:2002}
Guillaume Bal and Yvon Maday.
\newblock A {\textquotedblleft}parareal{\textquotedblright} time discretization
  for non-linear {PDE}'s with application to the pricing of an american put.
\newblock In {\em Lecture Notes in Computational Science and Engineering},
  pages 189--202. Springer Berlin Heidelberg, 2002.

\bibitem{bellen_zennaro:1989}
Alfredo Bellen and Marino Zennaro.
\newblock Parallel algorithms for initial-value problems for difference and
  differential equations.
\newblock {\em Journal of Computational and Applied Mathematics},
  25(3):341--350, May 1989.

\bibitem{berry_al:2012}
L.A. Berry, W.~Elwasif, J.M. Reynolds-Barredo, D.~Samaddar, R.~Sanchez, and
  D.E. Newman.
\newblock Event-based parareal: A data-flow based implementation of parareal.
\newblock {\em Journal of Computational Physics}, 231(17):5945--5954, July
  2012.

\bibitem{boville:1991}
Byron~A. Boville.
\newblock Sensitivity of simulated climate to model resolution.
\newblock {\em Journal of Climate}, 4(5):469--485, May 1991.

\bibitem{brandt:1977}
Achi Brandt.
\newblock Multi-level adaptive solutions to boundary-value problems.
\newblock {\em Mathematics of Computation}, 31(138):333--390, 1977.

\bibitem{carpenter_al:2005}
M.~H. Carpenter, C.~A. Kennedy, Hester Bijl, S.~A. Viken, and Veer~N. Vatsa.
\newblock Fourth-order runge{\textendash}kutta schemes for fluid mechanics
  applications.
\newblock {\em Journal of Scientific Computing}, 25(1):157--194, October 2005.

\bibitem{chen_al:2014}
Feng Chen, Jan~S. Hesthaven, and Xueyu Zhu.
\newblock {\em On the Use of Reduced Basis Methods to Accelerate and Stabilize
  the Parareal Method}, pages 187--214.
\newblock Springer International Publishing, Cham, 2014.

\bibitem{cox_matthews:2002}
S.~M. Cox and P.~C. Matthews.
\newblock Exponential time differencing for stiff systems.
\newblock {\em Journal of Computational Physics}, 176(2):430--455, 2002.

\bibitem{dai_maday:2011}
Xiaoying Dai and Yvon Maday.
\newblock Stable parareal in time method for first and second order hyperbolic
  system.
\newblock Technical report, December 2011.

\bibitem{durran:2010}
Dale Durran.
\newblock {\em Numerical methods for fluid dynamics : with applications to
  geophysics}.
\newblock Springer, New York, 2010.

\bibitem{ECMWF:2003}
{ECMWF} and P.~White.
\newblock {IFS Documentation CY25R1 - Part III: Dynamics and Numerical
  Procedures}.
\newblock 2003.

\bibitem{emmett_minion:2012}
Matthew Emmett and Michael~L. Minion.
\newblock {Toward an Efficient Parallel in Time Method for Partial Differential
  Equations}.
\newblock {\em Communications in Applied Mathematics and Computational
  Science}, 7:105--132, 2012.

\bibitem{falgout_al:2014}
R.~D. Falgout, S.~Friedhoff, T.~V. Kolev, Scott~P. MacLachlan, and Jacob~B.
  Schroder.
\newblock {Parallel time integration with multigrid}.
\newblock {\em SIAM Journal on Scientific Computing}, 36:C635--C661, 2014.

\bibitem{friedhoff_al:2013}
S.~Friedhoff, R.~D. Falgout, T.~V. Kolev, Scott~P. MacLachlan, and Jacob~B.
  Schroder.
\newblock {A Multigrid-in-Time Algorithm for Solving Evolution Equations in
  Parallel}.
\newblock In {\em {Presented at: Sixteenth Copper Mountain Conference on
  Multigrid Methods, Copper Mountain, CO, United States, Mar 17 - Mar 22,
  2013}}, 2013.

\bibitem{galewski_al:2004}
Joseph Galewsky, Richard~K. Scott, and Lorenzo~M. Polvani.
\newblock An initial-value problem for testing numerical models of the global
  shallow-water equations.
\newblock {\em Tellus A: Dynamic Meteorology and Oceanography}, 56(5):429--440,
  January 2004.

\bibitem{gander:2008}
Martin~J. Gander.
\newblock {Analysis of the Parareal Algorithm Applied to Hyperbolic Problems
  using Characteristics}.
\newblock {\em Bol. Soc. Esp. Mat. Apl.}, 42:21--35, 2008.

\bibitem{gander:2015}
Martin~J. Gander.
\newblock 50 years of time parallel time integration.
\newblock In {\em Contributions in Mathematical and Computational Sciences},
  pages 69--113. Springer International Publishing, 2015.

\bibitem{gander_al:2018}
Martin~J. Gander, Felix Kwok, and Hui Zhang.
\newblock Multigrid interpretations of the parareal algorithm leading to an
  overlapping variant and {MGRIT}.
\newblock {\em Computing and Visualization in Science}, 19(3-4):59--74, June
  2018.

\bibitem{gander_vanderwalle:2007}
Martin~J. Gander and Stefan Vandewalle.
\newblock Analysis of the parareal time‐parallel time‐integration method.
\newblock {\em SIAM J. Scientific Computing}, 29:556--578, 01 2007.

\bibitem{gaudreault_al:2022}
St{\'{e}}phane Gaudreault, Martin Charron, Valentin Dallerit, and Mayya Tokman.
\newblock High-order numerical solutions to the shallow-water equations on the
  rotated cubed-sphere grid.
\newblock {\em Journal of Computational Physics}, 449:110792, January 2022.

\bibitem{gaudreault_pudykiewicz:2016}
St{\'{e}}phane Gaudreault and Janusz~A. Pudykiewicz.
\newblock An efficient exponential time integration method for the numerical
  solution of the shallow water equations on the sphere.
\newblock {\em Journal of Computational Physics}, 322:827--848, October 2016.

\bibitem{geiser_guttel:2012}
J\"{u}rgen Geiser and Stefan G\"{u}ttel.
\newblock Coupling methods for heat transfer and heat flow: Operator splitting
  and the parareal algorithm.
\newblock {\em Journal of Mathematical Analysis and Applications},
  388(2):873--887, April 2012.

\bibitem{gotschel_al:2021}
Sebastian G\"{o}tschel, Michael Minion, Daniel Ruprecht, and Robert Speck.
\newblock Twelve ways to fool the masses when giving parallel-in-time results.
\newblock In {\em Springer Proceedings in Mathematics {\&} Statistics}, pages
  81--94. Springer International Publishing, 2021.

\bibitem{gotschel:2018}
Sebastian G\"{o}tschel and Michael~L. Minion.
\newblock Parallel-in-time for parabolic optimal control problems using
  {PFASST}.
\newblock In {\em Lecture Notes in Computational Science and Engineering},
  pages 363--371. Springer International Publishing, 2018.

\bibitem{hack_ruediger:1992}
James Hack and Ruediger Jakob.
\newblock Description of a global shallow water model based on the spectral
  transform method.
\newblock Technical report, 1992.

\bibitem{hamilton_al:2008}
Kevin Hamilton, Yoshiyuki~O. Takahashi, and Wataru Ohfuchi.
\newblock Mesoscale spectrum of atmospheric motions investigated in a very fine
  resolution global general circulation model.
\newblock {\em Journal of Geophysical Research}, 113(D18), September 2008.

\bibitem{hamon_schreiber:2020}
Fran{\c{c}}ois~P. Hamon, Martin Schreiber, and Michael~L. Minion.
\newblock Parallel-in-time multi-level integration of the shallow-water
  equations on the rotating sphere.
\newblock {\em Journal of Computational Physics}, 407:109210, April 2020.

\bibitem{hamon_al:2020}
François~P. Hamon, Martin Schreiber, and Michael~L. Minion.
\newblock Parallel-in-time multi-level integration of the shallow-water
  equations on the rotating sphere.
\newblock {\em Journal of Computational Physics}, 407:109210, 2020.

\bibitem{haut_wingate:2014}
Terry Haut and Beth Wingate.
\newblock An asymptotic parallel-in-time method for highly oscillatory {PDEs}.
\newblock {\em SIAM Journal on Scientific Computing}, 36(2):A693--A713, 2014.

\bibitem{hortal:2002}
Mariano Hortal.
\newblock The development and testing of a new two-time-level semi-lagrangian
  scheme ({SETTLS}) in the {ECMWF} forecast model.
\newblock {\em Quarterly Journal of the Royal Meteorological Society},
  128(583):1671--1687, July 2002.

\bibitem{jablonowski_williamson:2011}
Christiane Jablonowski and David~L. Williamson.
\newblock The pros and cons of diffusion, filters and fixers in atmospheric
  general circulation models.
\newblock In Peter Lauritzen, Christiane Jablonowski, Mark Taylor, and
  Ramachandran Nair, editors, {\em Numerical Techniques for Global Atmospheric
  Models}, chapter~13, pages 381--493. Springer Berlin Heidelberg, 2011.

\bibitem{kanamitsu_al:1983}
M.~Kanamitsu, K.~Tada, T.~Kudo, N.~Sato, and S.~Isa.
\newblock Description of the {JMA} operational spectral model.
\newblock {\em Journal of the Meteorological Society of Japan. Ser. {II}},
  61(6):812--828, 1983.

\bibitem{kelly_giraldo:2012}
James~F. Kelly and Francis~X. Giraldo.
\newblock Continuous and discontinuous galerkin methods for a scalable
  three-dimensional nonhydrostatic atmospheric model: Limited-area mode.
\newblock {\em Journal of Computational Physics}, 231(24):7988--8008, October
  2012.

\bibitem{koshyk_hamilton:2001}
John~N. Koshyk and Kevin Hamilton.
\newblock The horizontal kinetic energy spectrum and spectral budget simulated
  by a high-resolution
  troposphere{\textendash}stratosphere{\textendash}mesosphere {GCM}.
\newblock {\em Journal of the Atmospheric Sciences}, 58(4):329--348, February
  2001.

\bibitem{lauritzen_al:2011}
Peter Lauritzen, Christiane Jablonowski, Mark Taylor, and Ramachandran Nair,
  editors.
\newblock {\em Numerical Techniques for Global Atmospheric Models}.
\newblock Springer Berlin Heidelberg, 2011.

\bibitem{lindborg:1999}
Erik Lindborg.
\newblock Can the atmospheric kinetic energy spectrum be explained by
  two-dimensional turbulence?
\newblock {\em Journal of Fluid Mechanics}, 388:259--288, June 1999.

\bibitem{lions_al:2001}
Jacques-Louis Lions, Yvon Maday, and Gabriel Turinici.
\newblock R{\'{e}}solution d{\textquotesingle}{EDP} par un sch{\'{e}}ma en
  temps \quotes{parar{\'{e}}el}.
\newblock {\em Comptes Rendus de l{\textquotesingle}Acad{\'{e}}mie des Sciences
  - Series I - Mathematics}, 332(7):661--668, April 2001.

\bibitem{mengaldo_al:2018}
Gianmarco Mengaldo, Andrzej Wyszogrodzki, Michail Diamantakis, Sarah-Jane Lock,
  Francis~X. Giraldo, and Nils~P. Wedi.
\newblock Current and emerging time-integration strategies in global numerical
  weather and climate prediction.
\newblock {\em Archives of Computational Methods in Engineering},
  26(3):663--684, February 2018.

\bibitem{libpfasst:2018}
Michael Minion, Brandon Krull, Mathew Emmett, and Sebastian Goetschel.
\newblock libpfasst v1.0.
\newblock [Computer Software] \url{https://doi.org/10.11578/dc.20180711.5}, jun
  2018.

\bibitem{NOAA_NCEP:2016}
U.S.~National Oceanic and Atmospheric Administration (NOAA) / National~Centers
  for Environmental Prediction~(NCEP).
\newblock {Global Forecast System - Global Spectral Model (GSM) - V13.0.2}.
\newblock https://vlab.noaa.gov/web/gfs/documentation, May 2016.

\bibitem{ong_schroder:2020}
Benjamin~W. Ong and Jacob~B. Schroder.
\newblock Applications of time parallelization.
\newblock {\em Computing and Visualization in Science}, 23(1-4), September
  2020.

\bibitem{peixoto_schreiber:2019}
Pedro~S. Peixoto and Martin Schreiber.
\newblock Semi-lagrangian exponential integration with application to the
  rotating shallow water equations.
\newblock {\em SIAM Journal on Scientific Computing}, 41(5):B903--B928, 2019.

\bibitem{philippi_slawig:2023}
B.~Philippi and T.~Slawig.
\newblock {A Micro-Macro Parareal Implementation for the Ocean-Circulation
  Model FESOM2}, 2023.

\bibitem{philippi_slawig:2022}
Benedict Philippi and Thomas Slawig.
\newblock {The Parareal Algorithm Applied to the FESOM 2 Ocean Circulation
  Model}, 2022.

\bibitem{ries_trottenberg:1979}
M.~Ries and U.~Trottenberg.
\newblock Mgr-ein blitzschneller elliptischer löser.
\newblock {\em Tech. Rep. Preprint}, 277(SFB 72), 1979.

\bibitem{ruprecht:2018}
Daniel Ruprecht.
\newblock {Wave propagation characteristics of Parareal}.
\newblock {\em {Computing and Visualization in Science}}, 19:1--17, june 2018.

\bibitem{ruprecht_krause:2012}
Daniel Ruprecht and Rolf Krause.
\newblock Explicit parallel-in-time integration of a linear acoustic-advection
  system.
\newblock {\em Computers Fluids}, 59:72–83, Apr 2012.

\bibitem{samaddar_al:2019}
D.~Samaddar, D.P. Coster, X.~Bonnin, L.A. Berry, W.R. Elwasif, and D.B.
  Batchelor.
\newblock Application of the parareal algorithm to simulations of {ELMs} in
  {ITER} plasma.
\newblock {\em Computer Physics Communications}, 235:246--257, February 2019.

\bibitem{schmitt_al:2018}
A.~Schmitt, M.~Schreiber, P.~Peixoto, and M.~Sch\"{a}fer.
\newblock A numerical study of a semi-lagrangian parareal method applied to the
  viscous burgers equation.
\newblock {\em Computing and Visualization in Science}, 19(1-2):45--57, June
  2018.

\bibitem{schreiber_al:2019}
M.~Schreiber, N.~Schaeffer, and R.~Loft.
\newblock Exponential integrators with parallel-in-time rational approximations
  for shallow-water equations on the rotating sphere.
\newblock {\em Parallel Computing}, 2019.

\bibitem{shashkin_goyman:2020}
Vladimir~V. Shashkin and Gordey~S. Goyman.
\newblock Semi-lagrangian exponential time-integration method for the shallow
  water equations on the cubed sphere grid.
\newblock {\em Russian Journal of Numerical Analysis and Mathematical
  Modelling}, 35(6):355--366, December 2020.

\bibitem{skamarock_al:2012}
William~C. Skamarock, Joseph~B. Klemp, Michael~G. Duda, Laura~D. Fowler,
  Sang-Hun Park, and Todd~D. Ringler.
\newblock A multiscale nonhydrostatic atmospheric model using centroidal
  voronoi tesselations and c-grid staggering.
\newblock {\em Monthly Weather Review}, 140(9):3090--3105, September 2012.

\bibitem{staff_ronquist:2005}
Gunnar~Andreas Staff and Einar~M. R{\o}nquist.
\newblock Stability of the parareal algorithm.
\newblock In {\em Lecture Notes in Computational Science and Engineering},
  pages 449--456. Springer-Verlag, 2005.

\bibitem{staniforth:1991}
Andrew Staniforth and Jean C{\^{o}}t{\'{e}}.
\newblock Semi-lagrangian integration schemes for atmospheric
  models{\textemdash}a review.
\newblock {\em Monthly Weather Review}, 119(9):2206--2223, September 1991.

\bibitem{sterck_al:2022}
H.~De Sterck, R.~D. Falgout, and O.~A. Krzysik.
\newblock Fast multigrid reduction-in-time for advection via modified
  semi-lagrangian coarse-grid operators.
\newblock 2022.

\bibitem{sterck_al:2021}
Hans~De Sterck, Robert~D. Falgout, Stephanie Friedhoff, Oliver~A. Krzysik, and
  Scott~P. MacLachlan.
\newblock Optimizing multigrid reduction-in-time and parareal coarse-grid
  operators for linear advection.
\newblock {\em Numerical Linear Algebra with Applications}, 28(4), March 2021.

\bibitem{de_sterck_al:2019}
Hans~De Sterck, Stephanie Friedhoff, Alexander J.~M. Howse, and Scott~P.
  MacLachlan.
\newblock Convergence analysis for parallel-in-time solution of hyperbolic
  systems.
\newblock {\em Numerical Linear Algebra with Applications}, 27(1), November
  2019.

\bibitem{swarztrauber:2004}
Paul~N. Swarztrauber.
\newblock Shallow water flow on the sphere.
\newblock {\em Monthly Weather Review}, 132(12):3010--3018, December 2004.

\bibitem{temperton:1997}
Clive Temperton.
\newblock Treatment of the coriolis terms in semi-lagrangian spectral models.
\newblock {\em Atmosphere-Ocean}, 35(sup1):293--302, January 1997.

\bibitem{trindade_pereira:2004}
Joana M. F.~da Trindade and José~F. Pereira.
\newblock Parallel-in-time simulation of the unsteady {Navier–Stokes}
  equations for incompressible flow.
\newblock {\em International Journal for Numerical Methods in Fluids},
  45(10):1123--1136, 2004.

\bibitem{washington:2008}
Warren~M Washington, Lawrence Buja, and Anthony Craig.
\newblock The computational future for climate and earth system models: on the
  path to petaflop and beyond.
\newblock {\em Philosophical Transactions of the Royal Society A: Mathematical,
  Physical and Engineering Sciences}, 367(1890):833--846, December 2008.

\bibitem{williamson:2007}
David~L. Williamson.
\newblock The evolution of dynamical cores for global atmospheric models.
\newblock {\em Journal of the Meteorological Society of Japan. Ser. {II}},
  85B:241--269, 2007.

\bibitem{williamson_al:2008}
David~L. Williamson.
\newblock Convergence of aqua-planet simulations with increasing resolution in
  the community atmospheric model, version 3.
\newblock {\em Tellus A: Dynamic Meteorology and Oceanography}, 60(5):848,
  January 2008.

\bibitem{williamson_al:1992}
David~L. Williamson, John~B. Drake, James~J. Hack, R\"{u}diger Jakob, and
  Paul~N. Swarztrauber.
\newblock A standard test set for numerical approximations to the shallow water
  equations in spherical geometry.
\newblock {\em Journal of Computational Physics}, 102(1):211--224, September
  1992.

\end{thebibliography}

\end{document}